\documentclass[11pt]{article}
\usepackage[utf8]{inputenc}
\usepackage{amssymb,amsmath,amsthm,graphicx,geometry,hyperref,xcolor,bm,mathtools,authblk,subfigure,tikz,multirow,tensor}
\usepackage[numbers]{natbib}

\newgeometry{left=1in,right=1in,top=1in,bottom=1in}

\hypersetup{
    colorlinks=true,
    linkcolor=blue,
    filecolor = blue
, urlcolor = blue,
citecolor = blue
}

\newtheorem{theorem}{Theorem}

\newtheorem{corollary}{Corollary}
\newtheorem{lemma}{Lemma}
\newtheorem{definition}{Definition}

\newtheorem{innercustomthm}{Theorem}
\newenvironment{customthm}[1]
  {\renewcommand\theinnercustomthm{#1}\innercustomthm}
  {\endinnercustomthm}

\newtheorem{innercustomlem}{Lemma}
\newenvironment{customlem}[1]
  {\renewcommand\theinnercustomlem{#1}\innercustomlem}
  {\endinnercustomlem}

\newtheorem{innercustomcor}{Corollary}
\newenvironment{customcor}[1]
  {\renewcommand\theinnercustomcor{#1}\innercustomcor}
  {\endinnercustomcor}

\newcommand{\tr}{\text{tr}}

\newcommand{\Sym}{\text{Sym}(p)}
\newcommand{\Sp}{\mathcal{S}^p_{++}}
\newcommand{\Skron}{\mathcal{S}^{p_1,p_2}_{++}}
\newcommand{\Glp}{\text{GL}_p}
\newcommand{\Glkron}{\text{GL}_{p_1,p_2}}
\newcommand{\Okron}{\text{O}_{p_1,p_2}}
\newcommand{\Wishp}{\text{Wishart}_p}

\begin{document}

\title{Information Geometry and Asymptotics for Kronecker Covariances}
\author{Andrew McCormack} 
\author{Peter Hoff}
\affil{Department of Statistical Science, Duke University} 
\date{\today}

\maketitle

\begin{abstract}
We explore the information geometry and asymptotic behaviour of estimators for Kronecker-structured covariances, in both growing-$n$ and growing-$p$ scenarios, with a focus towards examining the quadratic form or partial trace estimator proposed by Linton and Tang \cite{LintonQuadratic}. It is shown that the partial trace estimator is asymptotically inefficient  An explanation for this inefficiency is that the partial trace estimator does not scale sub-blocks of the sample covariance matrix optimally. To correct for this, an asymptotically efficient, rescaled partial trace estimator is proposed. Motivated by this rescaling, we introduce an orthogonal parameterization for the set of Kronecker covariances. High-dimensional consistency results using the partial trace estimator are obtained that demonstrate a blessing of dimensionality. In settings where an array has at least order three, it is shown that as the array dimensions jointly increase, it is possible to consistently estimate the Kronecker covariance matrix, even when the sample size is one. 
   
\smallskip
 \noindent\textit{Keywords:} exponential family; Fisher information metric; high-dimensional; Kronecker product; orthogonal parameterization; partial trace; tensor.   
\end{abstract}

\section{Introduction}
An order-$k$ tensor, or synonymously, an array with a 
$k$-tuple of indices, $Y$, is a collection of numbers $Y_{i_1 \ldots i_k} \in \mathbb{R}$ where the index $i_j$ ranges from $1$ to $p_j$ for each $j = 1,\ldots,k$ \cite{Kolda2009}. Vectors and matrices give familiar examples of order one and order two arrays respectively. Each index set $i_j$ is referred to as a mode of the tensor, where the dimension of mode $j$ is defined as $p_j$. In data analysis settings the $j$th mode can be interpreted as a factor with $p_j$ different levels. An observation of a random tensor is tantamount to obtaining real-valued measurements $Y_{i_1\ldots i_k}$ at every $k$-way combination of factors.  For instance, the international trade data  examined in Section \ref{Sec:SimulationandData}, consists of an order-3 tensor $Y_{ijt}$, representing the trade from country $i$ to country $j$ at time $t$ \cite{hoff2017amen}. Some other examples of tensor-valued data include microarray data \cite{allen2012inference}, financial time-series \cite{cheng2018modeling}, fMRI data \cite{shvartsman2018matrix}, and MIMO communication signals \cite{da2018tensor}.

In this work we study a covariance matrix estimation problem for tensor data. As the number of entries of a tensor $Y$ is $p \coloneqq \prod_{i = 1}^k p_i$, the covariance matrix of any vectorization of $Y$ is ${p + 1 \choose 2}$-dimensional. 
When the order of the tensor or the dimension of the modes are large, a substantial number of parameters have to be estimated, leading to unstable estimates of the covariance matrix. Alternatively, restrictions can be placed on the hypothesized covariance structure for which more stable estimators are available, assuming that the restricted covariance submodel is approximately valid. The Kronecker covariance structure is one such restriction that leverages the tensor structure of $Y$, and describes the influence of each of the index-factors on the covariance of $Y$ in an interpretable manner \cite{dawidmatrixnormal}.





A random tensor $Y$ is defined to have the Kronecker covariance matrix $\bigotimes_{i = 1}^k \Sigma_i = \Sigma_1 \otimes \cdots \otimes \Sigma_k$ when 
\begin{align}
\label{eqn:KronDef1}
    \text{Cov}(Y_{i_1\ldots i_k},Y_{j_1\ldots j_k}) = \prod_{l = 1}^k [\Sigma_l]_{i_l j_l},
\end{align}
and each $\Sigma_i$ is in the set  is a $\mathcal{S}^{p_i}_{++}$  of non-singular, $p_i \times p_i$ dimensional, positive-definite matrices. The notation 
\begin{align}
\label{eqn:KronCovNotation}
    \mathcal{S}^{p_1,\ldots,p_k}_{++} \coloneqq \bigg\{ \bigotimes_{i = 1}^k \Sigma_i: \Sigma_i \in \mathcal{S}^{p_i}_{++}\bigg\}
\end{align}
denotes the set of Kronecker-structured covariance matrices. The number of parameters needed to define a Kronecker covariance is on the order of $\sum_{i = 1}^k p_i^2$, a number that is typically much smaller than the order $\prod_{i = 1}^k p_i^2$ parameters needed to defined a covariance matrix with no constraints. 

An interpretation of the Kronecker factors $\Sigma_i$ is that the covariance along mode $i$ is proportional to $\Sigma_i$. Fixing all modes but mode $j$, equation \eqref{eqn:KronDef1} implies the proportionality equation
\begin{align}
\label{eqn:modecovproportion}
\text{Var}\left(
    \begin{bmatrix}
Y_{i_1\ldots 1 \ldots i_k}
\\
\vdots
\\
Y_{i_1 \ldots p_j \ldots i_k}
    \end{bmatrix} \right) = \left(\prod_{l \neq j} [\Sigma_l]_{i_l i_l}\right) \Sigma_j \propto \Sigma_j.
\end{align}
When $k = 2$ and $Y$ is a random matrix, the content of equation \eqref{eqn:modecovproportion} is that the covariance of every row of $Y$ is proportional to $\Sigma_2$ and the covariance of every column is proportional to $\Sigma_1$. Kronecker covariance structures arise naturally through matrix, and more generally multilinear, multiplication. If $Z$ is a $p_1 \times p_2$ matrix with i.i.d. standard normal entries then $\Sigma_1^{1/2} Z \Sigma_2^{1/2}$ has the Kronecker covariance matrix $\Sigma_1 \otimes \Sigma_2$. 
It should be noted that alternative notation is commonly used, where the covariance of $\Sigma_1^{1/2} Z \Sigma_2^{1/2}$ is denoted by $\Sigma_2 \otimes \Sigma_1$ instead of $\Sigma_1 \otimes \Sigma_2$. This distinction amounts to vectorizing a matrix row-wise for $\Sigma_1 \otimes \Sigma_2$ or column-wise for $\Sigma_2 \otimes \Sigma_1$. To remain consistent with the definition \eqref{eqn:KronDef1}, where $\Sigma_i$ represents the covariance of mode $i$, we use the notation $\text{Var}(\Sigma_1^{1/2}Z\Sigma_2^{1/2}) \coloneqq \Sigma_1 \otimes \Sigma_2$.

Formalizing the scenario we consider throughout this work, let $Y^{(1)},\ldots,Y^{(n)} \overset{i.i.d.}{\sim} N_p(0,\bigotimes_{i = 1}^k \Sigma_i)$ be $n$ independent observations of tensors that each have the Kronecker covariance $\bigotimes_{i = 1}^k \Sigma_i$, defined by \eqref{eqn:KronDef1}. 
Equivalently,
$Y^{(1)},\ldots,Y^{(n)}$ can also be reduced to the sample covariance matrix sufficient statistic, $n S^{(n)} \sim \text{Wishart}_p(\bigotimes_{i = 1}^k \Sigma_i,n)$. We consider $S^{(n)}$ to be indexed as a $(p_1,\ldots,p_k,p_1,\ldots,p_k)$-dimensional array where
\begin{align*}
    S^{(n)}_{i_1\ldots i_k, j_1 \ldots j_k} \coloneqq \frac{1}{n} \sum_{l = 1}^n Y^{(l)}_{i_1\ldots i_k} Y^{(l)}_{j_1\ldots j_k},
\end{align*}
and thus
\begin{align}
\label{eqn:SampleCovExpectation}
    E\big(S^{(n)}_{i_1\ldots i_k, j_1 \ldots j_k}\big) = \text{Cov}\big( Y_{i_1\ldots i_k}^{(r)}, Y_{j_1\ldots j_k}^{(r)}\big) = \prod_{l = 1}^k [\Sigma_l]_{i_l j_l}, \;\; r = 1,\ldots,n .
\end{align}
Asymptotic properties of various estimators for $\bigotimes_{i = 1}^k \Sigma_i$ are analyzed in this article. A wide variety of asymptotic regimes are relevant for tensor data, as the numbers $p_1,\ldots,p_k,k$, and $n$ can all be adjusted. Previous work on the efficiency of estimators for Kronecker covariances can be found in \cite{Werner2008KronAsymp}. In \cite{EfronMircoarray} it is shown that the effective sample size of partial trace estimates of correlation matrices in a $k = 2$ setting is small when there is strong between-row or between-column correlation. Our approach  emphasizes geometric aspects of efficiency that lend additional insights into the behaviour of the partial trace estimator. Specifically, we show how the dispersion of the eigenvalues of the Kronecker factors impacts the asymptotic performance of the partial trace estimator. Prior work on high-dimensional sampling properties of the MLE and the partial estimator that we build upon are provided in \cite{herographicalLasso,LintonQuadratic}. 

In the first portion of this work we explore large $n$ and fixed $p_1,\ldots,p_k$ asymptotics.  To aid comprehension, these results are formulated primarily for order-two arrays. Order-$k$ generalizations follow in a straightforward manner.  Section \ref{Sec:WishartInfoGeom} begins by describing the information geometry of the Kronecker submodel, contained in the full, Wishart exponential family. Two important estimators for Kronecker covariances, the maximum likelihood estimator and the partial trace estimator, are characterized in Section \ref{Sec:KronEstimators}. An analysis of the level sets of these estimators given in Section \ref{Sec:AsymptoticEfficiency} shows that the partial trace estimator is asymptotically inefficient, except at the identity matrix. A rescaled version of the partial trace estimator is proposed and shown to be asymptotically efficient. Motivated by the rescaled partial trace estimator, an orthogonal parameterization for Kronecker covariances is introduced in Section \ref{Sec:OrthogonalParam}, where we describe the implications of this parameterization to testing problems. In Section \ref{Sec:LargePAsymptotics} a fixed $n$ and growing $(p_1,\ldots,p_k)$ asymptotic regime is analyzed. We establish mild conditions on $n,p_1,\ldots,p_k$ and the eigenvalues of $\otimes_{i = 1}^k \Sigma_i$ for which the partial trace estimator is consistent in relative Frobenius norm. A simulation study and a high-dimensional data example that support our theoretical results appear in the last section.

\section{Wishart Fisher Information Geometry}
\label{Sec:WishartInfoGeom}


In this section we introduce the relevant Fisher information geometry of the Wishart model that will be used to contrast the asymptotic performance of estimators for Kronecker covariance matrices. Given a regular, parametric model $\mathcal{P} = \{P_\theta:\theta \in \Theta \subset \mathbb{R}^m\}$, the Fisher information matrix at the distribution $P_\theta$ is the positive definite matrix with entries 
\begin{align}
\label{eqn:FIMatrix}
    [\mathcal{I}(\theta)]_{ij} = E_\theta\left(\partial_i\ell(\theta)\partial_j\ell(\theta)\right),
\end{align}
 where $\ell(\theta) = \log(p(y|\theta))$ is the log-likelihood function of $P_\theta$, and $\partial_i$ is the partial derivative with respect to the $i$th component of $\theta$. The matrix \eqref{eqn:FIMatrix} has direct statistical implications in that $\mathcal{I}(\theta)^{-1}$ is the asymptotic variance matrix of an optimal, asymptotically efficient estimator of $\theta$.   
 
Rather than deal with $\mathcal{I}(\theta)$ directly, it can be more convenient to work with the quadratic form associated with this matrix.
An inner product over $\mathbb{R}^m$ can be defined according to the rule 
 \begin{align}
 \label{eqn:FisherInfoMetric}
     \langle v,w \rangle_\theta \coloneqq v^\top \mathcal{I}(\theta) w \;\;\; v,w \in \mathbb{R}^m,
 \end{align}
which is referred to as the Fisher information metric (FIM) at $\theta$ \cite{AmariMethods}. The Fisher information metric gives the parameter space $\Theta$ the geometric structure of a Riemannian manifold where the FIM is a Riemannian metric \cite{lee2018Riem}. More generally, a Riemannian metric $\{\langle \cdot,\cdot\rangle_x\}_{x \in \mathcal{M}}$ over a manifold $\mathcal{M}$ is a collection of smoothly varying inner products, where $\langle \cdot,\cdot\rangle_x$ is an inner product over the tangent space $T_x\mathcal{M}$ at the point $x \in \mathcal{M}$. The tangent space $T_x\mathcal{M}$ can be viewed as the affine plane that intersects the surface $\mathcal{M}$ tangentially at $x$. Tangent spaces are especially simple when $\mathcal{M}$ is an open set in $\mathbb{R}^m$, as the approximating tangent plane is all of $\mathbb{R}^m$. This is the case for the FIM on a regular parametric model, where $\Theta$ is assumed to be open in $\mathbb{R}^m$. Further information on tangent spaces and differential geometry can be found in \cite{lee2013smooth,lee2018Riem}.


The unconstrained Wishart model $\mathcal{P} = \{\text{Wishart}_p(\Sigma,1):\Sigma \in \Sp\}$ is a regular exponential family \cite{BrownExpFam} that can be parameterized by the mean parameter $\Sigma \in \Sp$. The set $\Sp$  can be viewed as an open set in the ${p + 1\choose 2}$-dimensional vector space $\Sym$ of $p \times p$ dimensional symmetric matrices. The vector space $\Sym$ has the standard inner product
\begin{align}
    \label{eqn:SymInnerProduct}
    \langle A,B\rangle \coloneqq \tr(A^\top B) = \tr(AB), \; \text{for all } A,B \in \Sym.
\end{align}
Applying the above Fisher information construction to $\mathcal{P}$ in mean parameterization, the FIM can be succinctly represented \cite{SkovgaardMVN} as
\begin{align}
  \label{eqn:AIMetricMean}
    \langle H_1,H_2 \rangle_\Sigma & = \tfrac{1}{2} \tr\big(\Sigma^{-1/2}H_1\Sigma^{-1} H_2 \Sigma^{-1/2}\big), \;\; \Sigma \in \Sp, \; H_1,H_2 \in \Sym.
\end{align}
As $\Sp$ is an open set contained in $\Sym$, the tangent space $T_\Sigma \Sp$ at the point $\Sigma \in \Sp$ can be identified with $\Sym$ itself. The Riemannian metric \eqref{eqn:AIMetricMean}, also referred to as the affine-invariant Riemannian metric (AIRM), gives $\Sp$ the structure of a complete Riemannian manifold \cite[Ch 12]{LangDiffGeo}.

Affine-invariance refers to the invariance of the AIM  with respect to the congruence group action $\Sigma \mapsto A\Sigma A^\top$ of the general linear group $\Glp$ of invertible matrices $A$ on $\Sp$. This transitive action preserves the metric, as it is easily checked that
\begin{align}
\label{eqn:AIMetricIsom}
    \langle AH_1A^\top,AH_2A^\top \rangle_{A\Sigma A^\top} = \langle H_1,H_2\rangle_\Sigma. 
\end{align}
From a statistical perspective, this invariance is a result of the Wishart model being a group transformation family under $\Glp$. Specifically, if $S \sim \Wishp(\Sigma,1)$ then $ASA^\top \sim \text{Wishart}_p(A\Sigma A^\top, 1)$. The existence of this transitive group action simplifies many arguments, as a statement typically only has to be shown at a single, conveniently chosen, $\Sigma$.

The relationship between the FIM and the asymptotic variance matrix for asymptotically efficient estimators takes an especially simple form for sufficient statistics and mean parameters in regular exponential families. Projecting the sufficient statistic $S \sim \text{Wishart}_p(\Sigma,1)$ onto the one-dimensional subspaces determined by $A,B \in \Sym$, the covariance of these projections is given by
\begin{align}
\label{eqn:CovMetricFormula}
    \text{Cov}\big(\langle A,S\rangle_\Sigma, \langle B, S\rangle_\Sigma\big) = 4\langle A,B\rangle_\Sigma, \; \text{for all } A,B \in \Sym.
\end{align}
Formula \eqref{eqn:CovMetricFormula} is a restatement of the fact that the variance of the sufficient statistic in an exponential family is equal to the Fisher information matrix with respect to the corresponding canonical parameter \cite{AmariMethods}. See Lemma \ref{Lem:CovInnerProdExpFam} in the  Appendix for a proof.

\subsection{Kronecker Fisher Information Geometry}
We define the Kronecker submodel $\mathcal{P}_0 \coloneqq \{\text{Wishart}_p(\Sigma,1): \Sigma = \Sigma_1 \otimes \Sigma_2 \in \Skron\}$ of $\mathcal{P}$. Formally, $\mathcal{P}_0$ is a curved exponential family that can be parameterized by the submanifold $\Skron$ contained in $\Sp$ \cite[Sec 2.3]{KassVos}. That $\Skron$ has a manifold structure follows from the map $(a,B,C) \mapsto a (B\otimes C)$, where $B \in \mathcal{S}^{p_1}_{++}$, $C \in \mathcal{S}^{p_2}_{++}$ are matrices with trace $1$ and $a > 0$. This map gives a global parameterization of $\Skron$ in terms of a product of relatively open \cite[Sec 6]{rockafellar1997convex} sets in Euclidean space. As the dimension of the set of matrices in $\mathcal{S}^{p_i}_{++}$ with trace $1$ is ${p_i + 1 \choose 2} - 1$, the dimension of $\Skron$ as a manifold is ${p_1 + 1\choose 2} + {p_2 + 1 \choose 2}- 1$.

The tangent plane $ T_{\Sigma_1 \otimes \Sigma_2}\Skron$  at $\Sigma_1 \otimes \Sigma_2$ in $\Skron$ is more complicated than $T_{\Sigma_1 \otimes \Sigma_2}\Sp$, since the former tangent plane to the surface $\Skron$ depends on $\Sigma_1 \otimes \Sigma_2$. The subspace $T_{\Sigma_1 \otimes \Sigma_2}\Skron$ can be found by computing all possible derivatives of smooth curves $\Sigma_1(t) \otimes \Sigma_2(t)$ passing through $\Sigma_1 \otimes \Sigma_2$ at $t = 0$. By the product rule
\begin{align}
    \frac{d}{dt}\bigg\vert_{t = 0} \big(\Sigma_1(t) \otimes \Sigma_2(t)\big) = \Sigma_1'(0) \otimes \Sigma_2 + \Sigma_1 \otimes \Sigma_2'(0),
\end{align}
where $\Sigma_1'(0) \in T_{\Sigma_1}\mathcal{S}_{++}^{p_1}$ and $\Sigma_2'(0) \in T_{\Sigma_2}\mathcal{S}_{++}^{p_2}$ are symmetric matrices. As a subspace of $  T_{\Sigma_1 \otimes \Sigma_2}\Sp \cong \Sym$, the tangent space  $  T_{\Sigma_1 \otimes \Sigma_2}\Skron$ is the vector space
\begin{align}
\label{eqn:KronTangentSpace}
    T_{\Sigma_1 \otimes \Sigma_2}\Skron \cong \text{span}\{H_1 \otimes \Sigma_2,  \Sigma_1 \otimes H_2: H_i \in \text{Sym}(p_i)\}.
\end{align}
The submodel of $\mathcal{P}_0$ inherits the same Fisher information metric \eqref{eqn:AIMetricMean} as $\mathcal{P}$, restricted to the tangent space \eqref{eqn:KronTangentSpace}.

The collection of Kronecker covariances also has a natural group action on it that results from restricting the congruence action defined previously.  Let $\Glkron$ denote the subgroup of Kronecker products of non-singular matrices, $\Glkron \coloneqq \{A \otimes B: A \in \text{GL}_{p_1}(\mathbb{R}), B \in \text{GL}_{p_2}(\mathbb{R})\}$. The group $\Glkron$ acts transitively on $\Skron$ according to the map:
\begin{align*}
    \Sigma_1 \otimes \Sigma_2 \mapsto (A \otimes B)(\Sigma_1 \otimes \Sigma_2)(A \otimes B)^\top =  (A\Sigma_1 A^\top) \otimes (B\Sigma_2 B^\top) \in \Skron.
\end{align*}
Each such congruence transformation preserves the $\mathcal{P}_0$ FIM  by \eqref{eqn:AIMetricIsom}. A related subgroup of $\Glkron$ consisting of tensor products of orthogonal matrices is $\Okron \coloneqq \{U \otimes V: U \in \text{O}(p_1), \; V \in \text{O}(p_2)\}$. The significance of $\Okron$ is that any Kronecker covariance can be diagonalized by orthogonal matrices in this group.

\section{Estimators for Kronecker Covariances}
\label{Sec:KronEstimators}

\subsection{The Partial Trace Estimator}
\label{Subsec:PTEstimatorDef}

Moving to the problem of estimating $\Sigma_1 \otimes \Sigma_2$, given an observation $nS^{(n)} \sim \text{Wishart}_p(\Sigma_1 \otimes \Sigma_2,n)$, 
Linton and Tang introduce an estimator for $\Sigma_1 \otimes \Sigma_2$ in \cite{LintonQuadratic}, that they call the quadratic form estimator. The estimator can be motivated by the observation that \eqref{eqn:SampleCovExpectation} implies
\begin{align}
\label{eqn:PTMeanEqn}
    E\left(\sum_{j = 1}^{p_2} S_{i_1 j,i_2 j}^{(n)} \right) = \sum_{j = 1}^{p_2}[\Sigma_1]_{i_1i_2}[\Sigma_2]_{jj} = \tr(\Sigma_2) [\Sigma_1]_{i_1i_2}.
\end{align}
Consequently, 
by summing over certain blocks of the sample covariance matrix it is possible to construct an estimator with expectation that is proportional to $\Sigma_1$. Assuming that $S^{(n)}$ is the sample covariance of the i.i.d. Gaussian matrices $Y^{(1)},\ldots, Y^{(n)}$, the summation in \eqref{eqn:PTMeanEqn} can be viewed as summing the sample covariance matrix of each row of the $Y$ matrices. In \cite{EfronMircoarray} Efron examined drawbacks of using estimates for the correlation matrices of the Kronecker factors $\Sigma_1,\Sigma_2$ that are effectively based on partial traces. 

Generalizing this block-summation, the partial trace operators $\tr_1:\Sp \rightarrow \mathcal{S}^{p_1}_{++}$ and $\tr_2:\Sp \rightarrow \mathcal{S}^{p_2}_{++}$ are defined by
\begin{align}
\label{eqn:PartialTraceDef}
    [\tr_1(\Sigma)]_{i_1i_2} \coloneqq \sum_{j = 1}^{p_2} \Sigma_{i_1 j, i_2 j}, \;\; [\tr_2(\Sigma)]_{j_1j_2} = \sum_{i = 1}^{p_1} \Sigma_{ij_1, ij_2}.
\end{align}
The notation we have chosen for the partial trace estimator emphasizes that the image of $\tr_i$ is in $\mathcal{S}^{p_i}_{++}$. Alternative notation lets the index $i$ in $\tr_i$ denote the mode(s) over which the summation is performed.

Partial trace operators can be defined similarly in the tensor case by where $p = \prod_{i = 1}^k p_i$ and $\Sigma \in \Sp$ by $\tr_j:\Sp \rightarrow \mathcal{S}^{p_i}_{++}$
\begin{align}
\label{eqn:PTTensorDefinition}
       [\tr_j(\Sigma)]_{ii'} = \sum_{i_1,\ldots i_{j-1},i_{j+1},\ldots,i_k} \Sigma_{i_1\ldots i\ldots i_k,i_1\ldots i'\ldots i_k}.
\end{align}
Mathematically, if each mode corresponds to the $p_i$-dimensional vector space $V_i$, the matrix can be viewed as a linear operator $\Sigma \in (\otimes_{i = 1}^k V_i) \otimes (\otimes_{i = 1}^k V_i)^* \cong \otimes_{i = 1}^k (V_i \otimes V_i^*)$, where $V_i^*$ is the dual space of $V_i$. The partial trace operator $\tr_j$ is equal to contracting, or equivalently taking the trace, of every pair $V_i \otimes V_i^*$ except for the $j$th pair. See the Appendix for additional details.

The quadratic form estimator of $\Sigma_1 \otimes \Sigma_2$, which we refer to as the partial trace (PT) estimator, is defined in terms of these partial traces as
\begin{align}
    P(S^{(n)}) \coloneqq \frac{1}{\tr(S^{(n)})} \tr_1(S^{(n)}) \otimes \tr_2(S^{(n)}).
\end{align}
Each Kronecker factor $\tr_i(S^{(n)})$ in $P(S^{(n)})$ can be viewed as estimating $\Sigma_i$ up to scale, while the multiplicative factor $\tr(S^{(n)})^{-1}$ rescales $P(S^{(n)})$ to have the same trace as $S^{(n)}$. The partial trace estimator not only preserves the trace of $S^{(n)}$, but it has the stronger property of being the unique function into $\Skron$ that preserves the partial traces of $S^{(n)}$. Thus, $\tr_i(P(S^{(n)})) = \tr_i(S^{(n)}), \; i = 1,2$. By the linearity of the partial trace operator, $P(S^{(n)})$ is partial trace unbiased, in the sense that $E\big(\tr_i(P(S^{(n)}))\big) = \tr_i(\Sigma_1 \otimes \Sigma_2)$.

The PT estimator has intuitive appeal as it can be understood as aggregating the sample covariance matrix for every mode of an array and then taking the Kronecker product of the resulting sample covariance aggregates. It is also trivial to compute and it exists with probability one for any sample size $n$. However, it will be shown in subsequent sections that the PT estimator is asymptotically inefficient and performs poorly near the boundary of the cone $\Skron$.


\subsection{The Kronecker MLE}
The Kronecker maximum likelihood estimator (MLE) $M(S^{(n)}) = \hat{\Sigma}_1 \otimes \hat{\Sigma}_2$ is the solution to the following likelihood equations
\begin{align}
\label{eqn:MLElikelihood1}
    & \tfrac{1}{p_2} \tr_1\big((\hat{\Sigma}_1 \otimes \hat{\Sigma}_2)^{-1/2}S^{(n)}(\hat{\Sigma}_1 \otimes \hat{\Sigma}_2)^{-1/2}\big) = I_{p_1},
    \\
    \label{eqn:MLElikelihood2}
    & \tfrac{1}{p_1}\tr_2\big((\hat{\Sigma}_1 \otimes \hat{\Sigma}_2)^{-1/2}S^{(n)}(\hat{\Sigma}_1 \otimes \hat{\Sigma}_2)^{-1/2}\big) =  I_{p_2}.
\end{align}
These equations state the average partial traces of the decorrelated sample covariance matrix $M(S^{(n)})^{-1/2}S^{(n)}M(S^{(n)})^{-1/2}$ equal the identity matrix. Using the equivariance and cyclic permutation properties of the partial trace operator provided in Lemmas \ref{lem:PTpropertiesCyclicPerm} and \ref{lem:PTEquivarianceProp} in the Appendix, equivalent likelihood equations are
\begin{align}
\label{eqn:ComputationalMLELikelihood1}
    & \tfrac{1}{p_2} \tr_1\big((I_{p_1} \otimes \hat{\Sigma}_2)^{-1}S^{(n)}\big) = \hat{\Sigma}_1,
    \\
\label{eqn:ComputationalMLELikelihood2}
    & \tfrac{1}{p_1}\tr_2\big((\hat{\Sigma}_1 \otimes I_{p_2})^{-1}S^{(n)}\big) =  \hat{\Sigma}_2.
\end{align}
While there is no general closed form expression for $M(S^{(n)})$, a block-coordinate descent algorithm that initializes $\hat{\Sigma}_1$ at an arbitrary matrix, such as $I_{p_1}$, and iterates between the equations \eqref{eqn:ComputationalMLELikelihood1} and \eqref{eqn:ComputationalMLELikelihood2} can be used to find the MLE \cite{dutilleul1999mle}. Conditions on the sample size $n$ for which the MLE  exists, is non-singular, and is unique with probability one are explored in \cite{HoffExistenceMLE,DerksenExistence}, with a simple sufficient condition being $\tfrac{p_1}{p_2} + \tfrac{p_2}{p_1} < n$ \cite{DerksenExistence}. By taking the Kronecker structure into account, this sample size requirement is much less stringent than the requirement $p_1p_2 \leq n$ needed for the sample covariance matrix $S^{(n)}$ to be non-singular with probability one.

\section{Inefficiency of the Partial Trace Estimator}
\label{Sec:AsymptoticEfficiency}

While the partial trace estimator is simple to compute and analyze, it is shown in this section that it is inefficient, having a larger asymptotic variance than the MLE.
The key result that simplifies this analysis is that asymptotic efficiency, or lack thereof, can be appraised by an orthogonality property between subspaces with respect to the Fisher information metric. We begin by providing explicit expressions for these subspaces for both the MLE and PT estimators in Lemma \ref{lem:LemmaEffOrthog}. In Theorem \ref{Thm:AngleAsymptoticVar} the asymptotic variance of the partial trace estimator is related to a principle angle between subspaces, where this angle is shown to be related to the eigenstructure of $\Sigma_1 \otimes \Sigma_2$ in Lemma \ref{lem:PTasymptoticVarRatio}.  

\begin{figure}[h]
     \centering
     \begin{subfigure}
         \centering
    \includegraphics[width=0.45\textwidth]{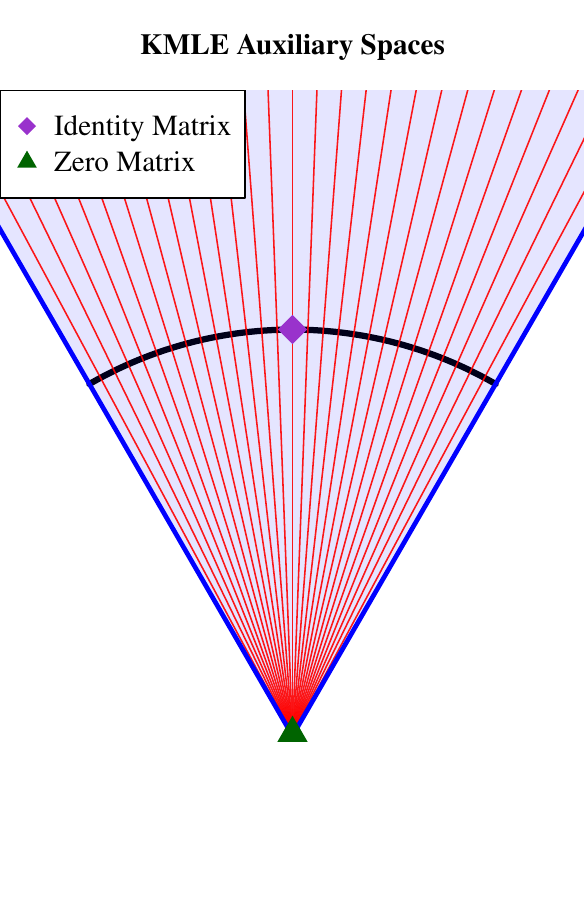}
     \end{subfigure}
     \begin{subfigure}
         \centering     \includegraphics[width=0.45\textwidth]{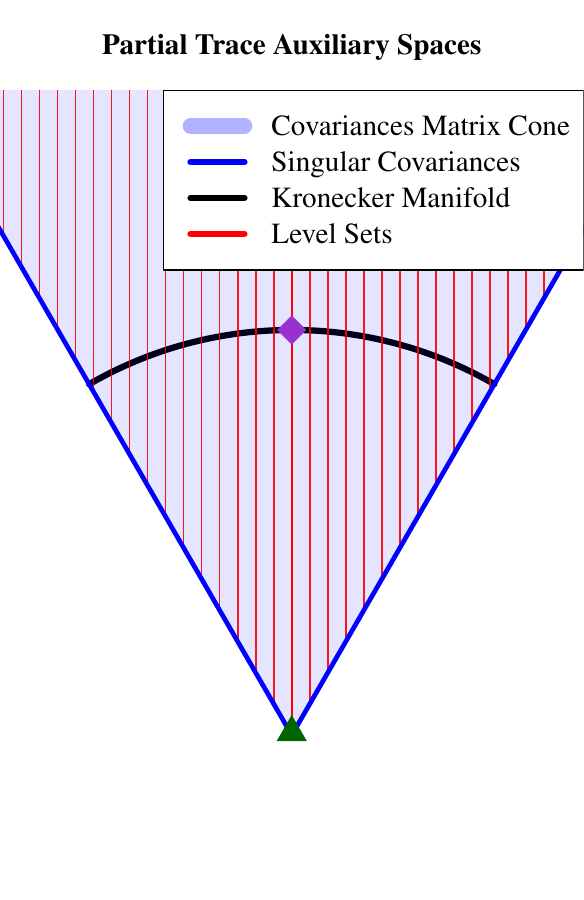}
     \end{subfigure}
        \caption[Auxiliary Spaces of the Kronecker MLE and Partial Trace Estimator.]{Sketch of the auxiliary spaces of $M(S)$ and $P(S)$ in the mean space cone $\Sp$.}
        \label{Fig:AuxSpaces}
\end{figure}

The first relevant tangent subspace, the tangent space to the auxiliary space, is one that is associated with a given estimator, and it determines whether or not an estimator is efficient. Any estimator  $\delta:\Sp \rightarrow \Skron$ of a Kronecker covariance maps a sample covariance matrix $S$ to the Kronecker covariance $\delta(S)$. Conversely, for each $\Sigma_1 \otimes \Sigma_2 \in \Skron$ the set of sample covariance matrices that get mapped to $\Sigma_1 \otimes \Sigma_2$ is defined as the the auxiliary space of $\delta$ at $\Sigma_1 \otimes \Sigma_2$ \cite[Ch 4]{AmariMethods}. That is, this auxiliary space $\mathcal{A}_{\Omega_1 \otimes \Omega_2}^\delta$ is the level set
\begin{align*}
\mathcal{A}_{\Sigma_1 \otimes \Sigma_2}^\delta = \delta^{-1}(\Sigma_1 \otimes \Sigma_2) = \{S: \delta(S) = \Sigma_1 \otimes \Sigma_2\}.
\end{align*}
The auxiliary space at $\Sigma_1 \otimes \Sigma_2$ is a possibly curved surface that passes through $\Sigma_1 \otimes \Sigma_2$. If $\delta$ is sufficiently smooth, the tangent space to $\mathcal{A}_{\Sigma_1 \otimes \Sigma_2}^\delta$ at the point $\Sigma_1 \otimes \Sigma_2$ is a subspace contained in the tangent space $T_{\Sigma_1 \otimes \Sigma_2}\Sp$. The tangent space to the Kronecker submanifold $T_{\Sigma_1 \otimes \Sigma_2}\Skron$ is also a subspace of $T_{\Sigma_1 \otimes \Sigma_2}\Sp$. The orthogonality result states that an estimator $\delta$ is asymptotically efficient if and only if $T_{\Sigma_1 \otimes \Sigma_2}\mathcal{A}_{\Sigma_1 \otimes \Sigma_2}^\delta$ is orthogonal to $T_{\Sigma_1 \otimes \Sigma_2}\Skron$ with respect to the FIM \cite[Thm 4.3]{AmariMethods}. A more precise formulation is given in Theorem \ref{Thm:AngleAsymptoticVar} below. First, we summarize the simple, affine structure of the auxiliary spaces of $M(S)$ and $P(S)$. See Figure \ref{Fig:AuxSpaces} for an illustration of these auxiliary spaces. We remark that this low-dimensional representation of $\Skron$ does not capture its cone structure, where $a\Sigma \in \Skron$ for any $a > 0$ and $\Sigma \in \Skron$.`

\begin{lemma} 
\label{lem:LemmaEffOrthog}
    Both $\mathcal{A}_{\Sigma_1 \otimes \Sigma_2}^M$ and
 $\mathcal{A}_{\Sigma_1 \otimes \Sigma_2}^P$ are the intersection of affine subspaces in $\Sym$ with the mean parameter space $\Sp$ of the exponential family $\mathcal{P}$. These auxiliary spaces are defined by the equations
 \begin{align*}
     \mathcal{A}_{\Sigma_1 \otimes \Sigma_2}^M & = \{S: \tr_i( (\Sigma_1 \otimes \Sigma_2)^{-1} (S - \Sigma_1 \otimes \Sigma_2)) = 0,\; i = 1,2\}
     \\
     \mathcal{A}_{\Sigma_1 \otimes \Sigma_2}^P & = \{S: \tr_i(S - \Sigma_1 \otimes \Sigma_2) = 0, \; i = 1,2\},
 \end{align*}
 with associated tangent subspaces
\begin{align*}
          T_{\Sigma_1 \otimes \Sigma_2}\mathcal{A}_{\Sigma_1 \otimes \Sigma_2}^M & = \{S: \tr_i( (\Sigma_1 \otimes \Sigma_2)^{-1}S) = 0,\; i = 1,2\}
     \\
T_{\Sigma_1 \otimes \Sigma_2} \mathcal{A}_{\Sigma_1 \otimes \Sigma_2}^P &  = \{S: \tr_i(S) = 0, \; i = 1,2\}.
\end{align*}
 At any positive multiple $c > 0$ of the identity $\mathcal{A}_{cI_p}^M = \mathcal{A}_{cI_p}^P$. The orthogonality condition necessary for asymptotic efficiency holds at all $\Sigma_1 \otimes \Sigma_2$ for the MLE, but only holds at the matrices that are proportional to the identity matrix for the partial trace estimator.
\end{lemma}
The auxiliary tangent spaces of $M$ and $P$ are both defined in terms of partial traces, with the distinction between these spaces being the decorrelation operation that appears in the the MLE auxiliary tangent space. As the subspaces  $T_{\Sigma_1 \otimes \Sigma_2} \mathcal{A}_{\Sigma_1 \otimes \Sigma_2}^P$  are identical as $\Sigma_1 \otimes \Sigma_2$ varies, the auxiliary spaces $\mathcal{A}_{\Sigma_1 \otimes \Sigma_2}^P$ are all parallel in $\Sp$.

Lemma \ref{lem:LemmaEffOrthog} shows that $P(S)$ is asymptotically inefficient, but it does not provide any indication of where $P$ performs poorly. Intuitively, the degree to which orthogonality between the auxiliary and Kronecker tangent spaces is violated corresponds to the degradation of the asymptotic performance of $P(S)$. Figure \ref{Fig:AuxSpaces} suggests that since the parallel, PT auxiliary spaces do not account for the bending of the Kronecker submanifold, the asymptotic performance of $P$ is the poorest near the boundary of $\Skron$. The following theorem generalizes Theorem 2.6.8 in \cite{KassVos} and relates the relative, worst-case asymptotic variance of the partial trace estimator to the degree that orthogonality is violated, as summarized by a principle angle:

\begin{theorem}
\label{Thm:AngleAsymptoticVar}
Let $\theta$ be the smallest principle angle between $T_{\Sigma_1 \otimes \Sigma_2}\mathcal{A}^T_{\Sigma_1 \otimes \Sigma_2}$ and $T_{\Sigma_1 \otimes \Sigma_2}\Skron$ with respect to the mean space inner product \eqref{eqn:AIMetricMean}. Then 
    \begin{align}
    \label{eqn:VarRatioAngleThm}
        \sup_{A \in \Sym} \frac{\text{AVar}_n(\langle A, P(S^{(n)}) \rangle_\Sigma )}{\text{AVar}_n(\langle A, M(S^{(n)}) \rangle_\Sigma)} = \sin(\theta)^{-2}
    \end{align}
    where $nS^{(n)} \sim \text{Wishart}_p(\Sigma_1 \otimes \Sigma_2,n)$, and $\text{AVar}_n$ is the asymptotic variance of the associated limiting normal distribution. 
\end{theorem}
This theorem expresses the asymptotic variance ratio of the PT and the MLE estimators, projected onto the one-dimensional line determined by $A$, in terms of a principle angle. As $M(S^{(n)})$ is asymptotically efficient, regardless of $A$, the projected asymptotic variance ratio will always be bounded below by one. Further insight into the variance ratio \eqref{eqn:VarRatioAngleThm} can be obtained by bounding the principle angle $\theta$ above. Equivariance considerations show that this principle angle only depends on the eigenvalues of $\Sigma_1 \otimes \Sigma_2$: 


\begin{lemma}
\label{lem:PrincipleAngleEquivariance}
The functions, $M$ and $P$ are $\Glkron$ and $\Okron$ equivariant respectively, implying that
the auxiliary spaces $\mathcal{A}_{\Sigma_1 \otimes \Sigma_2}^M$ are $\Glkron$ equivariant while the spaces $\mathcal{A}_{\Sigma_1 \otimes \Sigma_2}^P$ are $\Okron$ equivariant: 
 \begin{align*}
     (A \otimes B)\mathcal{A}_{\Sigma_1 \otimes \Sigma_2}^M(A \otimes B)^\top & = \mathcal{A}_{A\Sigma_1A^\top \otimes B\Sigma_2 B^\top}^M,\; \text{for all}\; A \otimes B \in \Glkron
     \\
     (U \otimes V)\mathcal{A}_{\Sigma_1 \otimes \Sigma_2}^P M(U \otimes V)^\top & = \mathcal{A}_{U\Sigma_1 U^\top \otimes V\Sigma_2 V^\top}^P, \;\text{for all}\; U \otimes V \in \Okron. 
 \end{align*}
 In particular, if $\Sigma_1 \otimes \Sigma_2$ has the eigendecomposition $U\Lambda U^\top \otimes V\Gamma V^\top$, all of the principle angles between $T_{\Sigma_1 \otimes \Sigma_2} \mathcal{A}^P_{\Sigma_1 \otimes \Sigma_2}$ and $T_{\Sigma_1 \otimes \Sigma_2}\Skron$ with respect to $\langle \cdot ,\cdot \rangle_{\Sigma_1 \otimes \Sigma_2}$ are identical to the principle angles between $T_{\Lambda \otimes \Gamma}\mathcal{A}^P_{\Lambda \otimes \Gamma}$ and $T_{\Lambda \otimes \Gamma}\Skron$ with respect to $\langle \cdot, \cdot \rangle_{\Lambda \otimes \Gamma}$.
\end{lemma}

Maximizing the inner product between vectors in the auxiliary and Kronecker tangent spaces gives the following bound on the asymptotic variance ratio:

\begin{lemma}
\label{lem:PTasymptoticVarRatio}
Denote the eigenvalues of $\Sigma_1$ and $\Sigma_2$ by $\lambda = (\lambda_1,\ldots,\lambda_{p_1})$ and $\gamma = (\gamma_1, \ldots, \gamma_{p_2})$ respectively. The maximal asymptotic variance ratio of the partial trace estimator is bounded below, as follow
    \begin{align}
    \label{eqn:VarRatioBoundEigenval1}
           \sup_{A \in \Sym} \frac{\text{AVar}_n(\langle A, P(S^{(n)}) \rangle_{\Sigma_1 \otimes \Sigma_2} )}{\text{AVar}_n(\langle A, M(S^{(n)}) \rangle_{\Sigma_1 \otimes \Sigma_2} )} \geq \max\bigg(  \frac{\Vert \lambda \Vert^2}{(\lambda^\top 1_{p_1}/\sqrt{p_1})^2},  \frac{\Vert \gamma \Vert^2}{(\gamma^\top 1_{p_2}/\sqrt{p_2})^2}  \bigg).
    \end{align}
    If $\alpha_1$ and  $\alpha_2$  are the angles between $\lambda, 1_{p_1}$ and $\gamma, 1_{p_2}$ then the lower bound \eqref{eqn:VarRatioBoundEigenval1} can be written as
    \begin{align}
        \label{eqn:VarRatioBoundEigenval2}
          \sup_{A \in \Sym} \frac{\text{AVar}_n(\langle A, P(S^{(n)}) \rangle_{\Sigma_1 \otimes \Sigma_2} )}{\text{AVar}_n(\langle A, M(S^{(n)}) \rangle_{\Sigma_1 \otimes \Sigma_2} )} \geq \max\big(\cos(\alpha_1)^{-2}, \cos(\alpha_2)^{-2}\big).
    \end{align}
\end{lemma}

The bound \eqref{eqn:VarRatioBoundEigenval2} confirms the geometric intuition garnered from Figure \ref{Fig:AuxSpaces} that the PT estimator performs poorly near the boundary of $\Skron$. Specifically, the angles $\alpha_1$ and $\alpha_2$ are large when the eigenvalues of $\Sigma_1$ and $\Sigma_2$ are widely dispersed, and are small when these eigenvalues are all nearly identical. When $\Sigma_1 \otimes \Sigma_2$ is proportional to $I_p$, the lower bound in \eqref{eqn:VarRatioBoundEigenval1} is $1$ as $P(S^{(n)})$ is efficient in this case. This qualitative behaviour is consistent with Theorem 2 of \cite{EfronMircoarray}, where it is shown that the effective sample size of estimators of Kronecker correlation matrix factors is reduced when there is significant across-row or across-column correlation.

The largest possible limiting angle between a vector of positive eigenvalues and the vector $1_{p_i}$ occurs when the vector of eigenvalues approaches any standard basis vector. Consequently, the lower bound \eqref{eqn:VarRatioBoundEigenval1} is maximized in the limit when $\Sigma_1 = \text{diag}(1,\epsilon,\ldots,\epsilon)$, $\Sigma_2 = \text{diag}(1,\epsilon,\ldots,\epsilon)$, and $\epsilon \rightarrow 0$. The limiting value of the lower bound, as $\epsilon$ converges to zero, is
\begin{align}
\label{eqn:AsympVarBoundOnBoundary}
    \lim_{\epsilon \rightarrow 0} \;\max\bigg( \frac{1 + (p_1 - 1)\epsilon^2 }{((1 + (p_1 - 1)\epsilon)/\sqrt{p_1})^2 },  \frac{1 + (p_2 - 1)\epsilon^2 }{((1 + (p_2 - 1)\epsilon)/\sqrt{p_2})^2 } \bigg) = \max(p_1,p_2).
\end{align}
The implication of \eqref{eqn:AsympVarBoundOnBoundary} is that in high-dimensional settings the PT estimator can perform much worse than the MLE, possibly having an asymptotic variance along a particular direction that is $\max(p_1,p_2)$ times larger. In the next section we provide an interpretation for why this is the case and how the PT estimator can be modified to obtain an asymptotically efficient estimator.

\section{An Efficient Rescaled Partial Trace Estimator}
\label{Sec:RescalecPT}
Consider the setting where $\Sigma_2 = \text{diag}(1,\epsilon,\ldots,\epsilon)$ and $n$ is large. The partial trace operator $\tr_1(S^{(n)}) = \sum_{j = 1}^{p_2}S^{(n)}_{\cdot j, \cdot j}$ is a sum of $p_2$ matrix-blocks that are roughly equal to 
\begin{align*}
    S_{\cdot 1,\cdot 1}^{(n)} & \approx \Sigma_1
    \\
    S_{\cdot j, \cdot j}^{(n)} & \approx \epsilon \Sigma_1  \;\; j > 1.
\end{align*}
For small values of $\epsilon$, the value of the partial trace is dominated by the first term in the sum with
\begin{align*}
    \tr_1(S^{(n)}) \approx S_{\cdot 1,\cdot 1} + (p_2-1)O(\epsilon)
\end{align*}
Due to the discrepancy of the magnitudes of the terms in the partial trace sum, the PT trace estimator effectively discards the $p_2 - 1$ terms $S_{\cdot j,\cdot j}, j > 1$ that have a negligible impact on the sum. This explains why the worst-case performance of the PT estimator compared to the MLE deteriorates as $\max(p_1,p_2)$ grows, as shown by equation \eqref{eqn:AsympVarBoundOnBoundary}. As evidenced by the efficiency of $P(S^{(n)})$ at scalar multiples of the identity, this scaling issue disappears when $\Sigma_2 \propto I_{p_2}$.  To correct for this scaling imbalance, we introduce a rescaled version of the partial trace operator.

\begin{definition}
The determinant-rescaled partial trace operators $\tr_i^{\det}: \Sp \rightarrow \mathcal{S}^{p_i}_{++}$ are defined as
\begin{align}
    \tr_1^{\det}(\Sigma) \coloneqq \sum_{j = 1}^{p_2} \frac{S_{\cdot j, \cdot j}}{\vert S_{\cdot j, \cdot j}\vert^{1/p_1}}, \;\;     \tr_2^{\det}(\Sigma) \coloneqq \sum_{j = 1}^{p_1} \frac{S_{j\cdot, j\cdot}}{\vert S_{j\cdot , j\cdot }\vert^{1/p_2}}.
\end{align}
    The rescaled partial trace estimator $(RPT)$ $R:\Sp \rightarrow \Skron$ is defined in terms of the determinant-rescaled partial trace operators by
       \begin{align}
       \label{eqn:RPTDefinition}
           R(S^{(n)}) & \coloneqq \bigg(\frac{\vert S^{(n)} \vert}{\vert \tr_1^{\det}(S^{(n)}) \otimes \tr_2^{\det}(S^{(n)}) \vert }\bigg)^{1/p}  \bigg(\tr_1^{\det}(S^{(n)}) \otimes \tr_2^{\det}(S^{(n)}) \bigg).
       \end{align}
\end{definition}
Each operator $\tr_i^{\det}$ is a sum of matrices with unit determinant. The reason for rescaling by determinants rather than a different scale-equivariant function, such as the trace, is that the determinant is a homomorphism from $\text{GL}_{p_i}$ to the multiplicative group of real numbers. This equivariance ensures that the rescaled partial trace estimator is $\Glkron$ equivariant. 

\begin{theorem}
\label{thm:RescaledPTEfficiency}
       The function $R$ is $\Glkron$ equivariant. If $nS^{(n)} \sim \text{Wishart}_p(\Sigma_1 \otimes \Sigma_2)$, the estimator $R(S^{(n)})$ is asymptotically efficient.
\end{theorem}

Similar to the PT estimator matching partial traces, the multiplicative factor appearing in \eqref{eqn:RPTDefinition} is chosen so that $R(S^{(n)})$ has the same determinant as $S^{(n)}$. Compared to the MLE, the RPT estimator has the benefit of having a closed form expression. 

It is helpful to rewrite equation \eqref{eqn:RPTDefinition} as
\begin{align}
\label{eqn:RPTDef2}
 R(S^{(n)}) = \vert S^{(n)} \vert^{1/p}  \bigg( \frac{\tr_1^{\det}(S^{(n)})}{\vert \tr_1^{\det}(S^{(n)}) \vert^{1/p_1} } \otimes \frac{\tr_2^{\det}(S^{(n)})}{\vert \tr_2^{\det}(S^{(n)}) \vert^{1/p_2}} \bigg).
\end{align}
This formulation makes clear that $R(S^{(n)})$ is composed of the scaling factor $\vert S^{(n)} \vert^{1/p}$, and the Kronecker factors $\vert \tr_i^{\det}(S^{(n)}) \vert^{-1/p_i}\tr_i^{\det}(S^{(n)})$ that each have determinant one. The RPT estimator is non-singular with probability one only when  $\vert S^{(n)} \vert \neq 0$ and thus $n \geq p_1p_2$. While this is a restrictive sample size requirement, the determinant-one Kronecker factors, which are of primary interest, exist as long as $n \geq \max(p_1,p_2)$.

\section{Asymptotic Distribution Using an Orthogonal Parameterization}
\label{Sec:OrthogonalParam}

 In this section we explore an orthogonal parameterization of $\Skron$, motivated by the factors in \eqref{eqn:RPTDef2}. The orthogonality of this parameterization implies the asymptotic independence of certain estimates of the Kronecker factors, modulo scaling. Moreover, the asymptotic distribution of efficient estimates of the Kronecker factors is shown to behave like a Wishart distribution. As an application, we demonstrate how these properties simplify hypothesis tests for structured Kronecker factors.

The individual factors in $\Sigma_1$ and $\Sigma_2$ in the Kronecker product   $\Sigma_1 \otimes \Sigma_2$ are only identifiable up to scale as $\Sigma_1 \otimes \Sigma_2 = c \Sigma_1 \otimes c^{-1}\Sigma_2$ for $c > 0$. It is possible to decompose $\Sigma_1 \otimes \Sigma_2$ into an overall scale factor, along with the scale-invariant remnants of the Kronecker factors. Using the determinant as a measure of scale, we obtain a decomposition
\begin{align}
\label{eqn:PSDDetDecomposition}
    \Sigma_1 \otimes \Sigma_2 = \vert \Sigma_1 \otimes \Sigma_2 \vert^{1/p} \bigg( \frac{\Sigma_1}{\vert \Sigma_1 \vert^{1/p_1}} \otimes \frac{\Sigma_2}{\vert \Sigma_2 \vert^{1/p_2}} \bigg) \coloneqq e^c\big( \Tilde{\Sigma}_1 \otimes \Tilde{\Sigma}_2 \big),
\end{align}
that was previously used in \cite{hoffhigherorderLQ,hoffequivariantminimax}.
Letting $\mathbb{P}(\mathcal{S}^{p_i}_{++})$ denote the collection of $p_i \times p_i$ positive definite matrices with determinant one, which can be viewed as a projective space over the cone $\Sp$, the decomposition \eqref{eqn:PSDDetDecomposition} induces a bijective map $\varphi: \Skron \rightarrow \mathbb{R} \times \mathbb{P}(\mathcal{S}^{p_1}_{++}) \times \mathbb{P}(\mathcal{S}^{p_2}_{++})$.  

Alternative decompositions, such as dividing each $\Sigma_i$ by $\tr(\Sigma_i)$ are possible; however, \eqref{eqn:PSDDetDecomposition} has the useful property of providing an orthogonal parameterization \cite{reidcoxOrthog} of $\Skron$. Moreover, this parameterization inherits an equivariance property, detailed in Lemma \ref{lem:OrthogParamEquivariance} in the Appendix, with respect to the product group $\mathbb{R} \times \text{SL}_{p_1} \times \text{SL}_{p_2}$ of real numbers under addition and the special linear group under matrix multiplication. 

\begin{theorem}
\label{thm:OrthogFIMDecomposition}
    The map $\varphi: \Skron \rightarrow \mathbb{R} \times \mathbb{P}(\mathcal{S}^{p_1}_{++}) \times \mathbb{P}(\mathcal{S}^{p_2}_{++})$, with $\varphi^{-1}$ defined as $\varphi^{-1}(c,\Tilde \Sigma_1, \Tilde \Sigma_2) = e^c(\Tilde\Sigma_1 \otimes \Tilde\Sigma_2)$ is a bijection that is an orthogonal parameterization. The Riemannian metric on the Cartesian product $\mathbb{R} \times \mathbb{P}(\mathcal{S}^{p_1}_{++}) \times \mathbb{P}(\mathcal{S}^{p_2}_{++})$ induced from the FIM on $\Skron$ is a product metric, where $\mathbb{R}$ has the standard Euclidean metric scaled by a factor of $p_1p_2$ and each $\mathbb{P}(\mathcal{S}^{p_i}_{++}), \; i = 1,2$ is equipped with the restriction of the FIM \eqref{eqn:AIMetricMean} on $\mathcal{S}^{p_i}_{++}$ scaled by a factor of $p_j$, $j \neq i$.     
\end{theorem}

Restricting our attention to the first two components $\mathbb{R} \times \mathbb{P}(\mathcal{S}^{p_1}_{++})$ in the orthogonal parameterization, the Kronecker FIM is seen to be the same as the FIM from a related Wishart model:

\begin{corollary}
The map $\varphi_1: \mathcal{S}^{p_1}_{++} \rightarrow \mathbb{R} \times \mathbb{P}(\mathcal{S}^{p_1}_{++})$, with $\varphi_1^{-1}$  defined as $\varphi_1^{-1}(c,\Tilde{\Sigma}_1) = e^c \Tilde{\Sigma}_1$ is an orthogonal parameterization of the Wishart model $\mathcal{P}_1 = \{\text{Wishart}_{p_1}(\Sigma_1,p_2):\Sigma_1 \in \mathcal{S}^{p_1}_{++}\}$. The FIM induced on $\mathbb{R} \times \mathbb{P}(\mathcal{S}^{p_1}_{++})$ by $\varphi_1$ and $\mathcal{P}_1$ equals the FIM induced on $\mathbb{R} \times \mathbb{P}(\mathcal{S}^{p_1}_{++})$ by the Kronecker model $\mathcal{P} = \{ \text{Wishart}_{p}(\Sigma_1 \otimes \Tilde{\Sigma}_2,1): \vert \Tilde{\Sigma}_2 \vert = 1\}$. Equivalently, the block of the Fisher information matrix corresponding to $\Sigma_1$ is the same in both $\mathcal{P}_1$ and $\mathcal{P}$. An analogous statement holds for the $\mathbb{R} \times \mathbb{P}(\mathcal{S}^{p_2}_{++})$ component of $\Skron$. 
\end{corollary}



An implication of this result is that the asymptotic variance of an efficient estimator of $\Sigma_1$ in the Kronecker model where $\Sigma_2$ is parameterized to have determinant one, is the same as that in a Wishart model with $p_2$ more observations. The dimension of the second mode in the Kronecker model acts as a sample size multiplier when estimating $\Sigma_1$. This sample size amplification is explicit in the expression for the partial trace estimator, as $\tr_1$ sums over $p_2$ matrix blocks.  In the next section we explore this phenomena further in the context of large-$p$ asymptotics.    

In terms of the Fisher information matrix, parameter orthogonality amounts to the matrix being block diagonal. Although we do not investigate this here, orthogonality can be useful for Newton-Raphson or Fisher-scoring algorithms where the information matrix only has to be inverted block-wise \cite{pourahmadiOrthogonality,reidcoxOrthog}. As the inverse of the Fisher information matrix is the asymptotic variance matrix for an asymptotically efficient estimator, orthogonality also implies an asymptotic independence property:

\begin{corollary}
\label{Cor:AsymptIndependenceofRPT}
As $n \rightarrow \infty$, the components of $\sqrt{n}\big(\varphi(R(S^{(n)})) - \varphi(\Sigma)\big)$, namely
\begin{align*}
    \sqrt{n}\bigg(\tfrac{1}{p}\log( \vert S^{(n)} \vert) - c,\frac{\tr_1^{\det}(S^{(n)})}{\vert \tr_1^{\det}(S^{(n)}) \vert^{1/p_1} } - \Tilde\Sigma_1, \frac{\tr_2^{\det}(S^{(n)})}{\vert \tr_2^{\det}(S^{(n)}) \vert^{1/p_2}} - \Tilde\Sigma_2 \bigg),
\end{align*}
are asymptotically independent.
Similarly, the components of $\sqrt{n}\big( \varphi(M(S^{(n)})) - \varphi(\Sigma)\big)$ are asymptotically independent.
\end{corollary}

An application of this asymptotic independence is to hypothesis testing problems concerning the Kronecker factors. In many instances it is hypothesized that the Kronecker factors $\Sigma_1$ and $\Sigma_2$ are structured covariance matrices, such an autoregressive covariance, a graphical model, an exchangeable covariance, or an isotropic covariance \cite{herospatiotemporalKron,roykhatreecompoundsymmetry,sristavacompoundsymm,herographicalLasso}. Specifically, let  $\mathcal{H}_i \subset \mathcal{S}^{p_i}_{++}$ and define the hypotheses $H_i: \Sigma_i \in \mathcal{H}_i$, $i = 1,2$. Further, assume that the subset $\mathcal{H}_i$ is cone, implying that if $A \in \mathcal{H}_i$ then so is $cA$ for $c > 0$. The majority of relevant hypotheses, including all of those listed above and hypotheses that can be specified using only correlation matrices, satisfy this scale-invariance assumption.
The asymptotic independence of efficient estimates can be used to test the intersection hypothesis $H_1 \cap H_2$, and asymptotically control the Type I error rate:
\begin{lemma}
\label{lem:TestingIndep}
Assume that $\psi_i(S_i^{(n)})$ is a smooth, scale-invariant test statistic for the hypothesis $H_i:\Sigma_i \in \mathcal{H}_i$ in the model $np_j S_i^{(n)} \sim \text{Wishart}_{p_i}(\Sigma_i,np_j)$, $i,j = 1,2$, $i \neq j$, with $S_1^{(n)}$ and $S_2^{(n)}$ independent. Further assume that $c_i^{(n)}(\psi_i(S^{(n)}) - \psi(\Sigma_i))$ converges weakly to a non-degenerate distribution that does not depend on $\Sigma_i$ under $\mathcal{H}_i$. Let $H_1 \cap H_2: \{\Sigma_1 \otimes \Sigma_2: \Sigma_i \in \mathcal{H}_i, i = 1,2\}$ be the intersection hypothesis with respect to the model $nS^{(n)} \sim \text{Wishart}_{p}(\Sigma_1 \otimes \Sigma_2,n)$. The asymptotic distribution of the test statistics 
\begin{align}
\label{eqn:TestStatIndepWish}
\left(  c_1^{(n)}(\psi_1(S_1^{(n)}) - \psi_1(\Sigma_1)),c_2^{(n)}(\psi_2(S_2^{(n)}) - \psi_2(\Sigma_2))\right)
\end{align}
equal the asymptotic distribution of
\begin{align}
\label{eqn:TestStatRSPT}
\left(  c_1^{(n)}(\psi_1(\tr_1^{\det}(S^{(n)})) - \psi_1(\Sigma_1)),c_2^{(n)}(\psi_2(\tr_2^{\det}(S^{(n)})) - \psi_2(\Sigma_2))\right).
\end{align}
In particular, if $q_{i\alpha}$ is an asymptotic level-$\alpha$ quantile of $c_i^{(n)}(\psi_i(S_i^{(n)}) - \psi_i(\Sigma_i))$, the test that rejects $H_1 \cap H_2$ when $c_i^{(n)}(\psi_i(\tr_i^{\det}(S^{(n)})) - \psi_i(\Sigma_i)) \geq q_{i\alpha}$ for either $i = 1$ or $i = 2$ has asymptotic level $1 - (1-\alpha)^2 = \alpha(2-\alpha)$. Similar statements apply, by replacing all occurrences of $\tr_i^{\det}(S^{(n)})$ above with the Kronecker MLE factors $\hat{\Sigma}_i$. 
\end{lemma}


To illustrate Lemma \ref{lem:TestingIndep}, we test the hypotheses $H_1$ that $\Sigma_1$ is spherical and $H_2$ that $\Sigma_2$ follows a compound symmetry model given an observation $nS \sim \text{Wishart}_p(\Sigma_1 \otimes \Sigma_2,n)$. We use the test statistics $\psi_1(\hat{\Sigma}_1),\psi_2(\hat{\Sigma}_2)$ where the $\hat{\Sigma}_i$ are the respective Kronecker factors of the MLE, $\psi_1$ is a pivotal test statistic based on the affine-invariant distance, and $\psi_2$ is the likelihood ratio test statistic. Additional details regarding these tests are provided in Appendix \ref{App:SimulationDetails}. The contingency Table \ref{tab:TestingIndependenceTable} records the observed  joint proportion of rejections of $H_1$ and $H_2$, computed from $5000$ Monte Carlo samples. Also provided in Table \ref{tab:TestingIndependenceTable} are the fitted expected proportions of rejection if $\hat{\Sigma}_1$ and $\hat{\Sigma}_2$ were independent. As dictated by Lemma \ref{lem:TestingIndep}, it is seen that for large sample sizes the hypotheses $H_1$ and $H_2$ are approximately rejected independently. Even when $n$ is only $10$ independence approximately holds.

\begin{table}[ht]
\label{tab:TestingIndependenceTable}
\centering
\begin{tabular}{|r|r|r|rr|rr|}
\hline
&&& \multicolumn{2}{|c|}{Null} &  \multicolumn{2}{|c|}{Alternative} 
\\
\cline{4-7}
&&& \multicolumn{4}{|c|}{$H_2$} 
  \\
  \cline{4-7}
& & & No Reject & Reject & No Reject & Reject\\ 
  \cline{2-7}
\multirow{ 6}{*}{$H_1$}& \multirow{ 2}{*}{$n = 10$}  & No Reject & 0.818 (0.814) & 0.106 (0.109) & 0.789 (0.786) & 0.132 (0.135) \\ 
& & Reject & 0.064 (0.067)  & 0.012 (0.009) & 0.064 (0.068) & 0.015 (0.012) \\ 
  \cline{2-7}
 & \multirow{ 2}{*}{$n = 50$} & No Reject & 0.887 (0.888) & 0.059 (0.058) & 0.657 (0.652)  & 0.209 (0.215) \\ 
& & Reject & 0.051 (0.051) & 0.003 (0.003) & 0.095 (0.100) & 0.039 (0.033) \\ 
  \cline{2-7}
 & \multirow{ 2}{*}{$n = 200$} & No Reject & 0.896 (0.896) & 0.049 (0.050) & 0.119 (0.116) & 0.399 (0.402) \\ 
& &  Reject & 0.051 (0.051) & 0.003 (0.003) & 0.105 (0.108) & 0.377 (0.374) \\ 
   \hline
\end{tabular}
\caption{Proportion of hypotheses rejected at a level of $0.95$ under the null hypothesis and under an alternative hypothesis for various sample sizes. The expected proportion of rejections under the assumption that the tests of $H_1$ and $H_2$ are independent are provided in parentheses.}
\end{table}

Operationally, Lemma \ref{lem:TestingIndep} can be used to simplify testing procedures for the hypothesis $H_1 \cap H_2$. A likelihood ratio test of $H_1 \cap H_2$ requires that the MLE under $H_1 \cap H_2$ be computed. This can be computationally expensive as it involves iterating a block-coordinate descent algorithm, where at each step a constrained estimate of $\Sigma_i$ under $H_i$ is found. Lemma \ref{lem:TestingIndep} shows that when $n$ is sufficiently large, it is enough to compute $M(S)$, or any other asymptotically efficient estimator of $\Sigma_1 \otimes \Sigma_2$, and apply tests developed for the Wishart distribution to each estimated Kronecker factor.

\section{Large-$p$ Asymptotics: A Blessing of Dimensionality}
\label{Sec:LargePAsymptotics}
The $p_2$ columns of the random matrix $Y \sim N(0,\Sigma_1\otimes \Sigma_2)$ each follow the normal distribution $Y_{\cdot i} \sim N(0, [\Sigma_2]_{ii}\Sigma_1)$. As $p_2$ increases, it is conceivable that up to scaling, the variability of estimates of $\Sigma_1$ decreases. In the special case that $\Sigma_2 = I_{p_2}$ this intuition exactly holds since $Y_{\cdot i} \overset{i.i.d.}{\sim} N(0,\Sigma_1)$, $i = 1,\ldots,p_2$. While $\Sigma_1$ may be easier to estimate as $p_2$ grows, there is a tradeoff in the Kronecker model as concomitantly the dimension of $\Sigma_2$ increases. In this section we provide convergence rates of the partial trace estimator as the dimensions $p_i$ grow. For matrix data, the effects of simultaneously increasing $p_1$ and $p_2$ cancel out with each other, and the Kronecker estimation problem in high-dimensions is no more difficult than in low-dimensions. Even more surprising is that for tensor-valued data, it becomes easier to estimate a Kronecker covariance in high-dimensional settings, with the estimator converging in relative Frobenius norm to the underlying Kronecker covariance as the $p_i$s jointly increase.

In this section it is assumed that $n_mS^{(m)} \sim \text{Wishart}_{p_m}(\otimes_{i = 1}^k \Sigma_i^{(m)},n_m)$, $\Sigma_i^{(m)}$ is $p_{im}$-dimensional, and $p_m = \prod_{i = 1}^k p_{im}$ is the dimension of $S^{(m)}$. The partial trace estimator is extended to the tensor setting by defining
\begin{align}
    P(S^{(m)}) \coloneqq \frac{1}{\tr(S^{(m)})^{k-1}} \bigotimes_{i = 1}^k \tr_i(S^{(m)}),
\end{align}
with the partial trace operators $\tr_i$ given by \eqref{eqn:PTTensorDefinition}. The scale factor $\tr(S^{(m)})^{-(k-1)}$ ensures that $\tr_i(P(S^{(m)})) = \tr_i(S^{(m)})$. 

Our central asymptotic result is the following:

\begin{theorem}
\label{thm:MainLargePConvergence}
    Let $n_mS^{(m)} \sim \text{Wishart}_{p_m}(\otimes_{i = 1}^k \Sigma_i^{(m)},n_m)$ with $\Sigma_i^{(m)} \in \mathcal{S}^{p_{im}}_{++}$ and $p_m = \prod_{i = 1}^k p_{im}$. Denote $\otimes_{i = 1}^k \Sigma_i^{(m)}$ by $\Sigma^{(m)}$, take $\lambda_i^{(m)}$ to be the length $p_{im}$ vector of eigenvalues of $\Sigma_i^{(m)}$, and let 
 $\lambda^{(m)}$ to be the length $p_{m}$ vector of eigenvalues of $\Sigma^{(m)}$. If
    \begin{align}
    \label{eqn:PTConvThmAssumptions}
        \frac{ \sup_{j \in \{1,\ldots,k\}}p_{jm}\cos^2\big(\lambda_j^{(m)},1_{p_{jm}}\big)}{\sqrt{n_m} \sqrt{p_m} \cos(\angle \lambda^{(m)}, 1_{p_{m}})}  = O(1), \; \frac{1}{\sqrt{n_m} \sqrt{p_m} \cos(\angle \lambda^{(m)},1_{p_{m}})} = o(1)
    \end{align}
    then 
\begin{align}
\label{eqn:PTConThmFrobeniusScale}
    \frac{\big\Vert P(S^{(m)}) - \Sigma^{(m)} \big\Vert_F}{\big\Vert \Sigma^{(m)} \big\Vert_F}  = O_P\bigg(       \frac{ \sup_{j \in \{1,\ldots,k\}}p_{jm} \cos^2\big(\angle \lambda_j^{(m)},1_{p_{jm}}\big) }{\sqrt{n_m} \sqrt{p_m} \cos(\angle \lambda^{(m)}, 1_{p_{m}})}\bigg)
\end{align}
\end{theorem}
Related convergence results in the matrix case for block-coordinate descent iterates of the MLE  can be found in Theorem 2 of \cite{herographicalLasso}. As compared to the theorem in \cite{herographicalLasso} Theorem \ref{thm:MainLargePConvergence},  although applied to a more tractable estimator, has advantages in that it gives a faster convergence rate by a factor of $\log(\max(\{p_{1m},\ldots,p_{km},n\})$, does not require any constraints on the sample size, and does not require a uniform bound on the spectrum of $\Sigma^{(m)}$. We note that the proof of Theorem \ref{thm:MainLargePConvergence} can be extended beyond the normal-Wishart family, to sample covariance matrices obtained from elliptically symmetric distributions with finite fourth-order moments.

The angle between $\lambda^{(m)}$ and the vector of all ones appears in \eqref{eqn:PTConThmFrobeniusScale}, where related angles were encountered previously in Lemma \ref{lem:PTasymptoticVarRatio}. This angle effectively replaces uniform spectral bounds. It is worth examining what kind of eigenvalue behaviour results in  $\cos^{-1}(\angle \lambda^{(m)}, 1_{p_{m}})$ being bounded asymptotically. The following lemma shows that $\sup_m\cos^{-1}(\angle \lambda^{(m)}, 1_{p_{m}})$ is finite under the mild condition that the eigenvalues grow at a polynomial rate.   

\begin{lemma}
  Let $C(x) = \sum_{q = 0}^{Q} c_q x^{q}$ be a non-zero polynomial with non-negative coefficients $c_q \geq 0$.
    If $\lambda^{(m)} = (C(1),C(2),\ldots,C(m))$, then $\cos^{-2}(\angle \lambda^{(m)},1_m) = O(1)$ as $m \rightarrow \infty$. If  $\lambda^{(m)} = (a^1,a^2,\ldots,a^m)$, $a > 1$ grows exponentially then $\cos^{-2}(\angle \lambda^{(m)},1_m) = O(m)$ as $m \rightarrow \infty$. Lastly, if $\lambda^{(m)} = (a_m + b_m,\ldots,a_m + b_m,b_m,\ldots,b_m)$ has a spiked eigenvalue structure with $q_m > 0$ entries that equal $a_m + b_m$, $m - q_m$ entries that equal $b_m$, and $a_m,b_m \geq 0$, then
    \begin{align*}
        \cos^{-2}(\angle \lambda^{(m)},1_m) = \frac{r_ma_m^2 + 2r_m a_mb_m + b_m^2}{r_m^2 a_m^2 + 2r_m a_mb_m + b_m^2}
    \end{align*}
where $r_m = q_m/m$. This expression is $O(1)$ if $a_m/b_m = O(1)$ or $m/q_m = O(1)$.
\end{lemma}
Under the assumption that $\sup_m\cos^{-1}(\angle \lambda^{(m)}, 1_{p_{m}}) < \infty$, and $p_{im} = \Theta(m)$ for all $i$, Theorem \ref{thm:MainLargePConvergence} shows that in the matrix case $P(S^{(m)})$ converges to $\Sigma^{(m)}$ in relative Frobenius norm at a rate of
\begin{align*}
        \frac{\big\Vert P(S^{(m)}) - \Sigma^{(m)} \big\Vert_F}{\big\Vert \Sigma^{(m)} \big\Vert_F} = O_P\big( n_m^{-1/2} m^{-(k-2)/2}\big).
\end{align*}
Consequently, a sufficient condition for relative Frobenius norm consistency in the matrix setting under these conditions is $n_m \rightarrow \infty$. For higher-order tensors, we  obtain the following consistency result that holds even when the sample size is only one:
\begin{corollary}
    If $k \geq 3$, $\sup_m\cos^{-1}(\angle \lambda^{(m)}, 1_{p_{m}}) < \infty$, and $p_{im} = \Theta(m)$, then under the conditions of Theorem \ref{thm:MainLargePConvergence},
    \begin{align*}
            \frac{\big\Vert P(S^{(m)}) - \Sigma^{(m)} \big\Vert_F}{\big\Vert \Sigma^{(m)} \big\Vert_F} = o_P(1),
    \end{align*}
    regardless of the sample size $n_m$.
\end{corollary}
A consequence of this corollary is that, even though the dimension of the parameter space is growing, it actually becomes easier to estimate $\Sigma^{(m)}$ in high-dimensions. In the $k = 3$ setting, there are $p_2p_3$ ``pseudo''-sample vectors $Y_{\cdot jk}$ that can be used to estimate $\Sigma_1^{(m)}$. The intuition is that as long as $p_2p_3/p_1$ is large enough, $\Sigma_1^{(m)}$ can be estimated consistently in relative Frobenius norm. This consistency result is possible because the Kronecker covariance assumption is a strong assumption, especially in high-dimensions where the dimension of the parameter space is reduced from $O(\prod_{i = 1}^k p_i^2)$ to $O(\sum_{i = 1}^k p_i^2)$.

\section{Simulation Study and Data Example}
\label{Sec:SimulationandData}

In this last section we provide an empirical assessment of the performance of the PT estimator, and compare it to the performance of the MLE in different settings. We also apply the PT estimator to a high-dimensional data set that only contains a single observation ($n = 1$). Despite this sample size limitation, the PT estimator is seen to produce plausible estimates for an underlying Kronecker covariance structure. 

\subsection{International Trade Data}
The UN Comtrade database contains information about the amount, in US dollars, of pairwise imports and exports between 229 countries across 9 different commodity types \cite{hoff2017amen}. We demonstrate in this section how the partial trace estimator can be applied to a subset of the high-dimensional UN Comtrade data. These data are high-dimensional in the sense that the sample size $(n = 1)$ is significantly smaller than the dimension of the observed tensor. Let $Y_{ijt}$ be the order-3 tensor with entries that equal the sum of the total quantity of trade, aggregated over all commodity types, exported from country $i$ to country $j$ in year $t$. Both countries $i$ and $j$ report the quantity of exports and imports respectively. These possibly differing reports are averaged when constructing $Y_{ijt}$. As some country pairings report no trade, a subset of 45 of the largest countries by trade, that all have non-zero pairwise imports and exports, is selected. The tensor $Y_{ijt}$ is $45 \times 45 \times 11$ dimensional, with missing values along the country $\times$ country diagonal $Y_{iit}$. Although a measure, such as the GDP for country $i$ in year $t$, could be imputed for $Y_{iit}$, the other entries of $Y$ are sufficiently distinct quantities that we opt to analyze the array without this imputation. 

Letting $L_{ijt} = \log(Y_{ijt})$,  it is postulated that the tensor $L$ is distributed according to the following autoregressive model on the log-scale:
\begin{align}
\label{eqn:TradeModel}
    L_{ijt} & =  L_{ij(t-1)} + \beta_t + \epsilon_{ijt} 
    \\
    \nonumber
    \epsilon & \sim N_{45 \times 45 \times 10}(0,\Sigma_1 \otimes \Sigma_2 \otimes \Sigma_3).
\end{align}
According to this model, the log-trade between countries $i$ and $j$ at time $t$ is equal to the trade in the previous year plus a noise term $\epsilon_{ijt}$ and a constant $\beta_t$, that captures the state of the global economy in year $t$. The noise $\epsilon_{ijt}$ is assumed to have a Kronecker covariance structure that reflects the correlation between the import and export behaviour across countries. For example, the import and export behaviour of countries that are signatories on a multilateral trade agreement may be similar. 
The baseline trade values $Y_{ij1}$ at $t = 1$, the magnitude of which differs greatly across country pairings, will not be analyzed here.

\begin{figure}
    \centering
    \includegraphics[width = .9\textwidth,height = 0.4\textwidth]{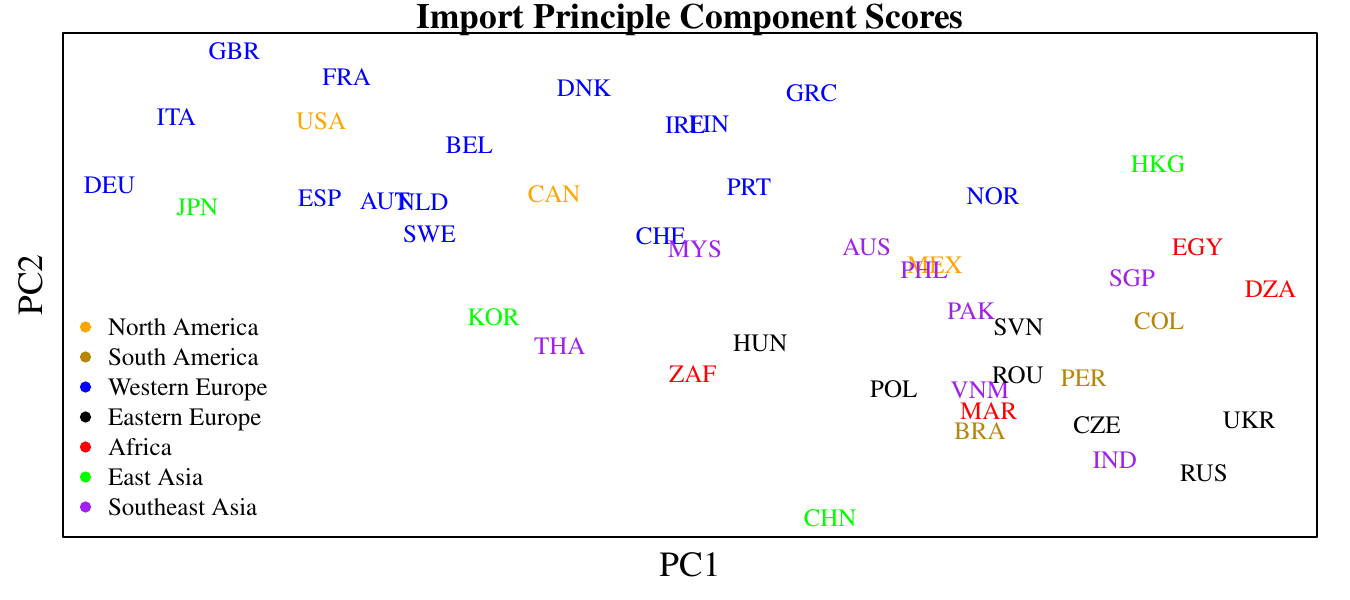}
    \caption{The first two principle component scores of the correlation matrix associated with $\hat{\Sigma}_2$, along with approximate country locations.}
    \label{fig:PCImportScores}
\end{figure}

We fit the model \eqref{eqn:TradeModel} by first defining $L_{ijt}^* = L_{ijt} - L_{ij(t-1)} \sim N(\beta_t,[\otimes_{a = 1}^3\Sigma_a]_{ijt})$. The  averages $\hat{\beta}_t = \tfrac{1}{45^2}\sum_{i,j = 1}^{45} L_{ijt}^*$ are used as estimates for the $\beta_t$'s. Estimates of the Kronecker covariance of the error term are found using $\Tilde{L}_{ijt} \coloneqq L_{ijt}^* - \hat{\beta}_t$. There are benefits to using partial trace estimator in this problem as it is simple to compute, even though the entries $\Tilde{L}_{iit}$ are not defined. The MLE is more challenging to compute as the iterates of the block-descent algorithm \eqref{eqn:ComputationalMLELikelihood1} and \eqref{eqn:ComputationalMLELikelihood2} are no longer well-defined. We use the approximate partial trace estimates
\begin{align*}
    [\hat{\Sigma}_1]_{ii'} & \propto \sum_{j: j \neq i,i'}\sum_t \Tilde{L}_{ijt}\Tilde{L}_{i'jt}
    \\
    [\hat{\Sigma}_2]_{jj'} & \propto \sum_{i: i \neq j,j'}\sum_t \Tilde{L}_{ijt}\Tilde{L}_{i'jt} 
    \\
        [\hat{\Sigma}_3]_{tt'} & \propto \sum_{i,j:i \neq j}\Tilde{L}_{ijt}\Tilde{L}_{i'jt'},
\end{align*}
and examine the correlation matrices associated with these estimates. The estimated correlation matrices along with a collection of other plots are provided in Appendix \ref{App:SimulationDetails}. It is seen that the estimated time correlation matrix  is close to the identity matrix, showing that the model \eqref{eqn:TradeModel} has largely accounted for time-effects through the autoregressive structure and $\beta_t$. The estimated import and export correlation matrices primarily have positive entries, with the import matrix generally having larger correlations. As might be expected, the import and export behaviour of high-income countries is positively correlated. To visualize the relationships between the import behaviour of countries, the first two principle component scores obtained from $\hat{\Sigma}_2$ for each country are plotted in Figure \ref{fig:PCImportScores}. It is clear that there is some relationship between these principle components and country income.

\begin{figure}
    \centering
    \includegraphics[width = .8\textwidth,height = .55\textwidth]{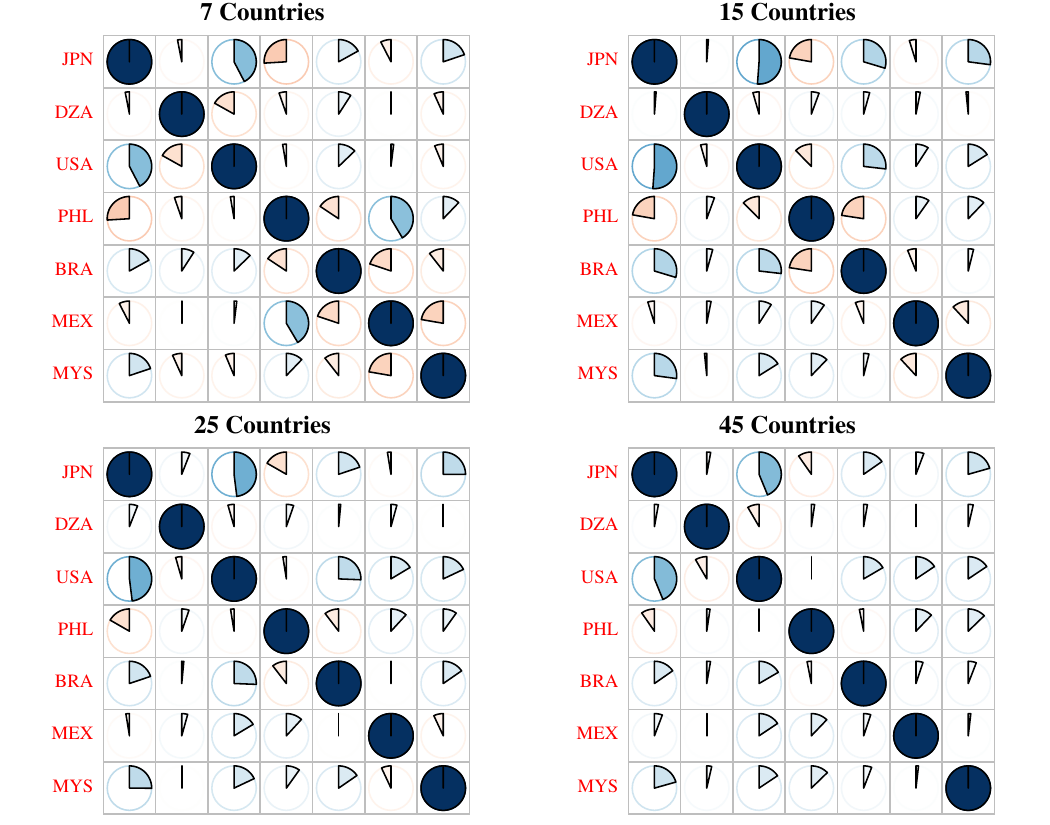}
    \caption{Estimated correlation matrices associated with $\hat{\Sigma}_1$. Each estimate is computed from a randomly selected subset of the entire $45$ countries.}
    \label{fig:CorplotHoldOut}
\end{figure}

It was argued in Section \ref{Sec:LargePAsymptotics} that the performance of partial trace estimator can improve as the dimension of the underlying tensor increases. To illustrate this, Figure \ref{fig:CorplotHoldOut} displays the estimated correlation matrix associated with $\Sigma_1$ that is estimated after holding out data. Specifically, each correlation matrix is computed from a sub-tensor of $\Tilde{L}_{ijt}$, where $i$ and $j$ range over a subset of the countries. As more countries are included, the correlation estimates appear to stabilize. In the next section we provide further empirical evidence of this phenomena in simulations.

\subsection{Simulation Study}
We present three different simulation regimes in this section. The first regime (Table \ref{tab:LargeNSimRisk}) has fixed values of $p_1,p_2$ and a large value of $n$, illustrating the asymptotic efficiency properties of the PT, RPT and MLE estimators. The second regime (Table \ref{tab:ModeratenSimRisk}) compares the PT and MLE estimators for sample sizes the are of a similar order to the array dimensions. In the last setting (Table \ref{tab:LargePKSimRisk}) we illustrate the performance of the PT estimator for high-dimensional tensor-valued data with a sample size of one. The squared, relative Frobenius norm, loss function $L(\delta,\Sigma) = n\Vert \Sigma \Vert^{-2}\Vert \delta - \Sigma \Vert^2$ is used to gauge the performance of the aforementioned estimators. As this loss function is orthogonally invariant and all of the estimators are orthogonally equivariant, the various risk functions will only depend on the eigenvalues of the Kronecker model from which observations are drawn. It is thus sufficient to examine the risk function for Kronecker models that have diagonal covariance matrices.

\begin{table}[h]
\centering
\begin{tabular}{|cc|cc|cc|cc|}
  \hline
 & & \multicolumn{2}{|c|}{$\lambda = (1,\ldots,1)$} &  \multicolumn{2}{|c|}{$\lambda = (1,\ldots,p_i)$} & \multicolumn{2}{|c|}{$\lambda = (e^1,\ldots,e^{p_i})$}
  \\
  \hline
 &  $p_1 = p_2 = $  & 3 & 7 & 3 & 7 & 3 & 7 \\ 
 \hline
& $n$ & & & & & &
\\
  \hline
\multirow{ 3}{*}{PT} &50 & 2.46 & 2.26 & 10.04 & 9.30 & 128.69 & 130.61 \\ 
 & 200 & 0.61 & 0.56 & 2.52 & 2.30 & 31.92 & 31.88 \\ 
  & 2500 & 0.05 & 0.04 & 0.20 & 0.18 & 2.52 & 2.54 \\
  \hline
\multirow{ 3}{*}{MLE} &  50 & 2.52 & 2.31 & 8.93 & 7.63 & 92.46 & 43.43 \\ 
 &  200 & 0.62 & 0.56 & 2.19 & 1.86 & 22.22 & 10.63 \\ 
  & 2500 & 0.05 & 0.05 & 0.17 & 0.15 & 1.74 & 0.87 \\
  \hline
\multirow{ 3}{*}{RPT} &  50 & 2.43 & 18.79 & 8.84 & 74.93 & 92.88 & 928.63 \\ 
 &  200 & 0.61 & 1.19 & 2.17 & 4.43 & 22.31 & 45.94 \\ 
  & 2500 & 0.05 & 0.05 & 0.17 & 0.16 & 1.74 & 1.08 \\ 
   \hline
\end{tabular}
\caption{Large $n$, squared, relative Frobenius norm, risk for the PT, MLE and RPT estimators. The underlying covariance matrix has the form $\Sigma_i = \text{diag}(\lambda)$, $i = 1,2$. Monte Carlo standard errors in all tables are less then $5\%$ of the magnitude of the cell entries.}
\label{tab:LargeNSimRisk}
\end{table}

\begin{table}[h]
\centering
\begin{tabular}{|cc|ccc|ccc|ccc|}
  \hline
 & & \multicolumn{3}{|c|}{$\lambda = (1,\ldots,1)$} &  \multicolumn{3}{|c|}{$\lambda = (1,\ldots,p_i)$} & \multicolumn{3}{|c|}{$\lambda = (e^1,\ldots,e^{p_i})$}
  \\
  \hline
 &  $p_1 = p_2 = $  & 3 & 7 & 20 & 3 & 7& 20 & 3 & 7 & 20\\ 
 \hline
& $n$ & & & & & & &&&
\\
  \hline
\multirow{ 3}{*}{PT} & 5 & 2.69 & 2.50 & 2.31 & 8.33 & 7.77 & 7.02 & 15.92 & 17.07 & 16.58 \\ 
 & 15 & 0.85 & 0.78 & 0.72 & 2.61 & 2.38 & 2.19 & 5.14 & 5.15 & 5.37 \\ 
 & 30 & 0.41 & 0.38 & 0.36 & 1.28 & 1.17 & 1.08 & 2.56 & 2.65 & 2.56 \\
  \hline
\multirow{ 3}{*}{MLE} &   5 & 4.72 & 4.12 & 3.67 & 12.21 & 9.85 & 8.27 & 18.39 & 8.40 & 9.03 \\ 
&  15 & 0.92 & 0.85 & 0.79 & 2.49 & 2.09 & 1.83 & 3.91 & 1.91 & 1.37 \\ 
&  30 & 0.43 & 0.40 & 0.37 & 1.16 & 0.98 & 0.86 & 1.84 & 0.90 & 0.49 \\ 
  \hline
\end{tabular}
\caption{Squared, relative Frobenius norm, risks in the moderate $n$ regime.}
\label{tab:ModeratenSimRisk}
\end{table}

\begin{table}[h]
\centering
\begin{tabular}{|c|ccc|ccc|}
  \hline
& & & & $\lambda = (1,\ldots,1) $ & $\lambda = (1,\ldots,p_i) $ & $\lambda = (e^1,\ldots,e^{p_i})$ \\ 
  \hline
\multirow{6}{*}{PT} &\multirow{3}{*}{$k = 2$} & \multirow{3}{*}{$p_1 = p_2 = $} & 5 & 0.91 & 1.02 & 1.30 \\ 
& & & 50 & 0.25 & 0.29 & 1.27 \\ 
& & & 200 & 0.12 & 0.14 & 1.29 \\ 
  \cline{2-7}
&\multirow{3}{*}{$p_i = 2$} &  \multirow{3}{*}{$k = $} & 4 & 1.12 & 1.20 & 1.27 \\ 
&& & 8 & 0.37 & 0.50 & 0.63 \\ 
&& & 16 & 0.03 & 0.07 & 0.12 \\ 
   \hline
\end{tabular}
\caption{Squared, Frobenius norm risk of the PT estimator for large $p_i$ and $k$ respectively, with a sample size of $n = 1$.}
\label{tab:LargePKSimRisk}
\end{table}

Summarizing the results in the tables, first note  that the RPT estimator has a noticeably large risk function in the large $n$ simulations. A large sample size of $n = 2500$ is needed for risk of the RPT estimator to be comparable to the MLE. Consequently, we do not recommend using the RPT estimator in practice, as the determinant rescaling introduces too much additional variability in the estimator. From both the large and moderate $n$ simulations we see that the PT estimator has performance that is comparable to the MLE in the cases where the eigenvalues are identical or linearly increasing. Interestingly, the first column of Table \ref{tab:ModeratenSimRisk} shows that the PT estimator has better small sample performance than the MLE when $\Sigma_i = I_{p_i}$. 
As predicted by our efficiency results for the PT estimator, the MLE appreciably outperforms the PT estimator when the eigenvalues have the form $\lambda = (e^1,\ldots,e^{p_i})$. These simulations suggest that if the eigenvalues of the true covariance matrix are not expected to be extremely dispersed then the PT estimator is a reasonable alternative to the MLE, especially in small sample settings or high-dimensional settings where it is computationally expensive to find the MLE.

Table \ref{tab:LargePKSimRisk} corroborates the consistency results in Section \ref{Sec:LargePAsymptotics}. When $k = 3$ the relative Frobenius norm risk of the PT estimator decreases as $p_i$ increases for eigenvalue structures that increase at a polynomial rate. However, when the eigenvalues increase at an exponential rate the risk of the PT estimator remains approximately constant as $p_i$ increases. Moreover, when the number of modes $k$ increases but $p_i$ is held fixed at $2$, the risk decreases at an exponential rate. The implication of these results is that the Kronecker covariance assumption, assuming that it is approximately true, is especially useful for modelling the covariance matrices of high-dimensional data, since the partial trace estimator is consistent with respect to $p_i$ under mild eigenvalues assumptions.

\section{Conclusion}
We have examined the performance of the partial trace estimator in various asymptotic settings in this article. A useful orthogonal parameterization of the Kronecker submodel was introduced in the course of our examinations. To conclude, we remark on two extensions that are relevant to the present work.  

All of our results are predicated on the covariance structure of the data under consideration having a Kronecker covariance matrix. Shrinkage estimators were that allow for the Kronecker covariance assumption to be violated are described in \cite{bersson2023bayesian,hoff2022core}. An alternative is to consider the larger class of models where the covariance matrix is assumed to be a sum of separable covariance matrices, meaning that it has the form $\sum_{r = 1}^R \otimes_{i = 1}^k \Sigma_i^{(r)}$ \cite{masakseparable,masak2022random,panaretosRanktest,greenewald2013kronecker}. We postulate that it is possible to obtain fixed $n$ and increasing $p_i$ consistency results for the Kronecker sum model with a fixed rank $R$ and $k > 2$. For this model, $\text{Cov}(Y_{\cdot j_2\ldots j_k},Y_{\cdot j_2\ldots j_k}) = \sum_{r = 1}^R \prod_{ i =2}^k [\Sigma_i^{(r)}]_{j_ij_i} \Sigma_1^{(r)}$. Letting $j_2,\ldots,j_k$ vary, it ostensibly is possible to recover each of $\Sigma_1^{(r)},\ldots,\Sigma_1^{(R)}$ up to scaling, and likewise for the Kronecker factors corresponding to the other modes of the array.

A second extension of the results presented here is to other curved exponential families that have a tensor product structure on the mean parameter space. An example of this is the following log-linear Poisson model for contingency table data $Y \in \mathbb{N}_0^{p_1} \times \cdots \mathbb{N}_0^{p_k}$:
\begin{align}
\label{eqn:SeparablePoisson}
    Y \sim \text{Poisson}( \otimes_{i = 1}^k \lambda_i), \;\;\lambda_i \in \mathbb{R}^{p_i},
\end{align}
where this notation is equivalent to $Y_{j_1\ldots j_k} \sim \text{Poisson}(\prod_{i = 1}^k [\lambda_{i}]_{j_i})$ independently. Each mode of the tensor $Y$ can be viewed as representing a factor or variable, with the $i$th factor having $p_i$ levels. The analogue of the $i$th partial trace operator is the sum over every mode but the $i$th mode. Reminiscent of \eqref{eqn:PTMeanEqn}, the expectation of the analogue of the first partial operator is
\begin{align*}
    E\left( \sum_{j_2,\ldots,j_k}Y_{\cdot j_2 \ldots j_k}\right) =\left(\prod_{i > 1} 1_{p_i}^\top \lambda_i\right) \lambda_1 \propto \lambda_1.
\end{align*}
It is expected that fixed-$n$ consistency results with a growing number of factors or levels within each factor are possible for this separable Poisson model. Additionally, the orthogonal parameterization introduced Section \ref{Sec:OrthogonalParam} has an analogue for the model \eqref{eqn:SeparablePoisson}. If $k = 2$, 
\begin{align*}
    \lambda_1 \otimes \lambda_2 \mapsto \left(1_{p_1}^\top \lambda_1 1^\top_{p_2}\lambda_2, \frac{\lambda_1}{1_{p_1}^\top \lambda_1},\frac{\lambda_2}{1_{p_2}^\top \lambda_2}\right) 
\end{align*}
is an orthogonal parameterization. In future work it would be interesting to characterize classes of exponential families for which the results in this article can be generalized.

\bibliographystyle{abbrv}
\bibliography{biblio}

\newpage

\begin{appendix}
\section{Proofs}

\begin{customlem}{0.1}
\label{Lem:CovInnerProdExpFam}
    Let $\mathcal{P} = \{P_\theta:\theta \in \Theta\}$ be a regular exponential family with densities that have the form $p(y|\theta) \propto \exp(\theta^\top y - \kappa(\theta))$. If $Y \sim P_\theta$  and $\mu(\theta) = E_\theta(Y)$ is the mean parameter of the exponential family then
    \begin{align*}
        \text{Cov}_\theta(\langle v,Y\rangle_\mu, \langle w,Y\rangle_\mu)  = \langle v,w \rangle_\mu,
    \end{align*}
    where $\langle v,w\rangle = v^\top \mathcal{I}(\mu) w$.
\end{customlem}
\begin{proof}
In a regular exponential family $\kappa(\theta)$ is the cumulant generating function of $Y$ with $\text{Var}_\theta(Y) = \nabla^2\kappa(\theta) = \mathcal{I}(\theta)$. Moreover, the information matrix in the mean parameterization is $\mathcal{I}(\mu(\theta)) = \mathcal{I}(\theta)^{-1}$. We compute
\begin{align*}
     \text{Cov}(\langle v,Y\rangle_\mu, \langle w,Y\rangle_\mu) & =    \text{Cov}_\theta\big( (\mathcal{I}(\mu)v)^\top Y , (\mathcal{I}(\mu)w)^\top Y \big) 
     \\
     & =  \big(\mathcal{I}(\mu) v\big)^\top \mathcal{I}(\theta) \big(\mathcal{I}(\mu) w\big)
     \\
     & = \langle v,w\rangle_\mu.
\end{align*}
Equation \eqref{eqn:CovMetricFormula} follows by using the sufficient statistic $-S/2$ and the corresponding canonical parameter $\Sigma^{-1}$ in the Wishart exponential family. 
\end{proof}

As noted in Section \ref{Subsec:PTEstimatorDef} the partial trace operator is formally defined as performing contractions on upper and lower indices of a tensor in $\bigotimes_{i = 1}^k (V_i^* \otimes V_i)$. Let $T \in \bigotimes_{i = 1}^k (V_i^* \otimes V_i)$ have coordinates $T^{i_1\ldots i_k}_{j_1\ldots j_k}$. Using the Einstein summation notation, where there is an implicit summation between repeated upper and lower indices, the $\ell$th partial trace operator is defined as 
\begin{align}
\label{eqn:PTContractionDef}
    \tr_\ell(T)^{i_\ell}_{j_\ell} = \delta_{j_1}^{i_1} \cdots \delta_{j_{\ell-1}}^{i_{\ell-1}} \delta_{j_{\ell+1}}^{i_{\ell+1}} \cdots \delta_{j_k}^{i_k} T^{i_1\ldots i_k}_{j_1\ldots j_k} = T^{i_1\ldots i_k}_{i_1 \ldots i_{\ell-1} j_\ell i_{\ell+1} \ldots i_k}
\end{align}
where $\delta_j^i$ is the Kronecker-delta tensor with coordinates $\delta_j^i = 1$ if $i = j$, $\delta_j^i = 0$ if $i \neq j$ \cite{McculaghTensor}. To be precise we have explicitly represented the tensor $T$ as having covariant (lower) and contravariant (upper) indices here. Throughout this article it is assumed that each vector space $V_i$ is endowed with an inner product so that there is a canonical identification $V_i \cong V_i^*$, and there is no distinction between upper and lower indices. Thus, the tensor $\Sigma_{i_1i_2,j_1j_2}$ is equal to $\Sigma^{i_1i_2}_{j_1j_2}$ with lowered indices. This tensor notation simplifies the verification of properties of the partial trace operator. For example, we have
\begin{align*}
    \tr(\tr_\ell(T)) = \tr_\ell(T)^{i_\ell}_{i_\ell} = T^{i_1\ldots i_k}_{i_1\ldots i_k} = \tr(T),
\end{align*}
and 
\begin{align*}
    \tr_\ell\left( \bigotimes_{\alpha = 1}^k \tensor[_\alpha]{A}{} \right)^{i_\ell}_{j_\ell} = \left( \bigotimes_{\alpha = 1}^k \tensor[_\alpha]{A}{} \right)^{i_1\ldots i_k}_{i_1 \ldots i_{\ell-1} j_\ell i_{\ell+1} \ldots i_k} = \tensor[_\ell]{A}{}^{i_\ell}_{j_\ell} \prod_{\alpha \neq \ell} \tensor[_\alpha]{A}{}^{i_\alpha}_{i_\alpha} =  \tensor[_\ell]{A}{}^{i_\ell}_{j_\ell} \prod_{\alpha \neq \ell} \tr\left(\tensor[_\alpha]{A}{} \right).
\end{align*}
The above equation generalizes \eqref{eqn:PTMeanEqn}.

\begin{customlem}{0.2}
    \label{lem:PTpropertiesCyclicPerm}
If $A \in \mathbb{R}^{p \times q}$ and $B \in \mathbb{R}^{nq \times np}$ 
then the partial trace operator satisfies the cyclic permutation property
\begin{align}
\label{eqn:CyclicPerm}
    \tr_1\left((I_n \otimes A)B\right) =     \tr_1\left(B(I_n \otimes A)\right). 
\end{align}
\end{customlem}
\begin{proof}
By the formula \eqref{eqn:PTContractionDef}
\begin{align*}
  \tr_1\left( (I_n \otimes A)B\right)^{i_1}_{j_1} & = \left((I_n \otimes A)B \right)^{i_1i_2}_{j_1i_2}
  \\
  & = (I_n \otimes A)^{i_1i_2}_{k_1k_2} B^{k_1k_2}_{j_1i_2}
  \\
  & = \delta^{i_1}_{k_1}A^{i_2}_{k_2} B^{k_1k_2}_{j_1i_2}
  \\
  & =  B^{i_1 k_2}_{j_1i_2} A^{i_2}_{k_2},
\end{align*}
while 
\begin{align*}
    \tr_1\left(B(I_n \otimes A)\right)^{i_1}_{j_1} & = \left(B(I_n \otimes A)\right)^{i_1i_2}_{j_1 i_2}
    \\
    & = B^{i_1i_2}_{k_1k_2} (I_n \otimes A)^{k_1k_2}_{j_1i_2}
    \\
     & =  B^{i_1i_2}_{k_1k_2} \delta^{k_1}_{j_1} A^{k_2}_{i_2}
     \\
     & = B^{i_1i_2}_{j_1 k_2}A^{k_2}_{i_2},
\end{align*}
giving the desired equality after relabeling indices.
\end{proof}

\begin{customlem}{0.3}
    \label{lem:PTEquivarianceProp}
    If $A \in \mathbb{R}^{p \times q}$, $B \in \mathbb{R}^{qn \times pn}$, and $C \in \mathbb{R}^{p \times r}$ 
then the partial trace operator satisfies the equivariance property
\begin{align}
\label{eqn:PTEquivariance}
    \tr_1\left( (A \otimes I_n)B(C \otimes I_n)\right) = A\tr_1(B)C
\end{align}
\end{customlem}
\begin{proof}
    Using formula \eqref{eqn:PTContractionDef} gives
    \begin{align*}
        \tr_1\left( (A \otimes I_n)B(C \otimes I_n)\right)^{i_1}_{j_1} & = \left((A \otimes I_n)B(C \otimes I_n)\right)^{i_1i_2}_{j_1i_2}
        \\
        & = (A \otimes I_n)^{i_1i_2}_{k_1k_2} B^{k_1k_2}_{l_1l_2}  (C \otimes I_n)^{l_1l_2}_{j_1i_2}
        \\
        & = A^{i_1}_{k_1} \delta^{i_2}_{k_2}  B^{k_1k_2}_{l_1l_2} C^{l_1}_{j_1} \delta^{l_2}_{i_2} 
        \\
        & = A^{i_1}_{k_1}  B^{k_1i_2}_{l_1i_2} C^{l_1}_{j_1}
        \\
        & = A^{i_1}_{k_1}\tr_1(B)^{k_1}_{l_1} C^{l_1}_{j_1}
        \\
        & = (A \tr_1(B) C)^{i_1}_{j_1},
    \end{align*}
    as needed.
\end{proof}

Both Lemma \ref{lem:PTpropertiesCyclicPerm} and Lemma \ref{lem:PTEquivarianceProp} similarly hold for $\tr_2$ and generalize to higher-order tensors.

\begin{customlem}{1}
    Both $\mathcal{A}_{\Sigma_1 \otimes \Sigma_2}^M$ and
 $\mathcal{A}_{\Sigma_1 \otimes \Sigma_2}^P$ are the intersection of affine subspaces in $\Sym$ with the mean parameter space $\Sp$ of the exponential family $\mathcal{P}$. These auxiliary spaces are defined by the equations
 \begin{align*}
     \mathcal{A}_{\Sigma_1 \otimes \Sigma_2}^M & = \{S: \tr_i( (\Sigma_1 \otimes \Sigma_2)^{-1} (S - \Sigma_1 \otimes \Sigma_2)) = 0,\; i = 1,2\}
     \\
     \mathcal{A}_{\Sigma_1 \otimes \Sigma_2}^P & = \{S: \tr_i(S - \Sigma_1 \otimes \Sigma_2) = 0, \; i = 1,2\},
 \end{align*}
 with associated tangent subspaces
\begin{align*}
          T_{\Sigma_1 \otimes \Sigma_2}\mathcal{A}_{\Sigma_1 \otimes \Sigma_2}^M & = \{S: \tr_i( (\Sigma_1 \otimes \Sigma_2)^{-1}S) = 0,\; i = 1,2\}
     \\
T_{\Sigma_1 \otimes \Sigma_2} \mathcal{A}_{\Sigma_1 \otimes \Sigma_2}^P &  = \{S: \tr_i(S) = 0, \; i = 1,2\}.
\end{align*}
 At any positive multiple $c > 0$ of the identity $\mathcal{A}_{cI_p}^M = \mathcal{A}_{cI_p}^P$. The orthogonality condition necessary for asymptotic efficiency holds at all $\Sigma_1 \otimes \Sigma_2$ for the MLE, but only holds at the matrices that are proportional to the identity matrix for the partial trace estimator.
\end{customlem}
\begin{proof}
    The form of $\mathcal{A}^{M}_{\Sigma_1 \otimes\Sigma_2}$ follows immediately as the condition $S \in \{S: \tr_i( (\Sigma_1 \otimes \Sigma_2)^{-1} (S - \Sigma_1 \otimes \Sigma_2)) = 0,\; i = 1,2\}$ is equivalent to $\Sigma_1 \otimes \Sigma_2$ solving the likelihood equations \eqref{eqn:MLElikelihood1} and \eqref{eqn:MLElikelihood2} with respect to an observation $S$. 

    If $S \in \{S: \tr_i(S - \Sigma_1 \otimes \Sigma_2) = 0, \; i = 1,2\}$ then $\tr_1(S) = \tr(\Sigma_2)\Sigma_1$, $\tr_2(S) = \tr(\Sigma_1)\Sigma_2$, and $\tr(S) = \tr(\Sigma_1)\tr(\Sigma_2)$. By the definition of $P$, $P(S) = \Sigma_1 \otimes \Sigma_2$, so $\{S: \tr_i(S - \Sigma_1 \otimes \Sigma_2) = 0, \; i = 1,2\} \subset \mathcal{A}^P_{\Sigma_1 \otimes \Sigma_2}$. Conversely, if $S \in \mathcal{A}^P_{\Sigma_1 \otimes \Sigma_2}$ then 
    \begin{align*}
      \tr_1(S) =  \tr_1(P(S)) = \tr_1(\Sigma_1 \otimes \Sigma_2) = \tr(\Sigma_2)\Sigma_1,
    \end{align*}
    with a symmetric equation holding for $\tr_2$. Thus, $S \in   \{S: \tr_i(S - \Sigma_1 \otimes \Sigma_2) = 0, \; i = 1,2\}$, proving the desired equality $\mathcal{A}^P_{\Sigma_1 \otimes \Sigma_2} =  \{S: \tr_i(S - \Sigma_1 \otimes \Sigma_2) = 0, \; i = 1,2\}$.
As $\mathcal{A}^M_{\Sigma_1 \otimes \Sigma_2}$ and $\mathcal{A}^P_{\Sigma_1 \otimes \Sigma_2}$ are affine, $T_{\Sigma_1 \otimes \Sigma_2}\mathcal{A}_{\Sigma_1 \otimes \Sigma_2}^M$ and $  T_{\Sigma_1 \otimes \Sigma_2}\mathcal{A}_{\Sigma_1 \otimes \Sigma_2}^P$ equal the translations of $\mathcal{A}^M_{\Sigma_1 \otimes \Sigma_2}$ and $\mathcal{A}^P_{\Sigma_1 \otimes \Sigma_2}$ to the origin.

Lastly, the orthogonality $ T_{\Sigma_1 \otimes \Sigma_2}\mathcal{A}_{\Sigma_1 \otimes \Sigma_2}^M \perp T_{\Sigma_1\otimes \Sigma_2}\Skron$ is easily checked from the definition of the Kronecker tangent space \eqref{eqn:KronTangentSpace} and the likelihood equations \eqref{eqn:MLElikelihood1} and \eqref{eqn:MLElikelihood2}. As there is a unique orthogonal complement to $T_{\Sigma_1\otimes \Sigma_2}\Skron$, $ T_{\Sigma_1 \otimes \Sigma_2}\mathcal{A}_{\Sigma_1 \otimes \Sigma_2}^P$ is orthogonal to $T_{\Sigma_1\otimes \Sigma_2}\Skron$ if and only if $ T_{\Sigma_1 \otimes \Sigma_2}\mathcal{A}_{\Sigma_1 \otimes \Sigma_2}^P =  T_{\Sigma_1 \otimes \Sigma_2}\mathcal{A}_{\Sigma_1 \otimes \Sigma_2}^M$. This occurs only for $(\Sigma_1 \otimes \Sigma_2)^{-1} \propto I_{p_1} \otimes I_{p_2}$. 
\end{proof}

\begin{customthm}{1}
\label{thm:AppendixAvarRatio}
Let $\theta$ be the smallest principle angle between $T_{\Sigma_1 \otimes \Sigma_2}\mathcal{A}^P_{\Sigma_1 \otimes \Sigma_2}$ and $T_{\Sigma_1 \otimes \Sigma_2}\Skron$ with respect to the mean space inner product \eqref{eqn:AIMetricMean}. Then 
    \begin{align}
        \sup_{A \in \Sym} \frac{\text{AVar}_n(\langle A, P(S^{(n)}) \rangle_\Sigma )}{\text{AVar}_n(\langle A, M(S^{(n)}) \rangle_\Sigma)} = \sin(\theta)^{-2}
    \end{align}
    where $nS^{(n)} \sim \text{Wishart}_p(\Sigma_1 \otimes \Sigma_2,n)$, and $\text{AVar}_n$ is the asymptotic variance of the associated limiting normal distribution. 
\end{customthm}
\begin{proof}
Let $D^M_\Sigma$ be the differential of $M$ at the point $\Sigma$. This differential can be viewed as a linear map from $T_\Sigma \Sp$ to $T_\Sigma \Skron$. The subspace $T_\Sigma \mathcal{A}_\Sigma^M$ is contained in the kernel of $D^M_\Sigma$ since $\tfrac{d}{dt}\big\vert_{t = 0} M(\Sigma(t)) = D^M_\Sigma(\Sigma'(0)) = 0$ for any curve $\Sigma(t)$ with image contained in $\mathcal{A}^M_\Sigma$ that passes through $\Sigma$ at $t = 0$. Any tangent vector in $T_\Sigma \Skron$ is not in the kernel of $D^M_\Sigma$ as $\tfrac{d}{dt}\big\vert_{t = 0} M((\Sigma_1 \otimes \Sigma_2)(t)) = D^M_\Sigma( (\Sigma_1 \otimes \Sigma_2)'(0)) =  (\Sigma_1 \otimes \Sigma_2)'(0)$ for any curve $(\Sigma_1 \otimes \Sigma_2)(t)$ that passes through $\Sigma$ at $t = 0$. This shows that the differential is idempotent. By counting dimensions it follows that $D^M_\Sigma$ has kernel $T_\Sigma \mathcal{A}_\Sigma^M$ and range $T_\Sigma \Skron$. The orthogonality established in Lemma \ref{lem:LemmaEffOrthog} implies that $D^M_\Sigma$ is the orthogonal projection onto $T_\Sigma \Skron$ with respect to the FIM.

An identical argument shows that the differential of $P$ at the point $\Sigma$ is idempotent, has kernel $T_\Sigma \mathcal{A}^P_\Sigma$ and range $T_\Sigma \Skron$. However, $D^P_\Sigma$ is in general only a projection, not an orthogonal projection.  

Using the delta method, equation \eqref{eqn:CovMetricFormula}, and the fact that $D_\Sigma^M$ is self-adjoint we compute
\begin{align*}
   \text{AVar}_n(\sqrt{n}\langle A, M(S^{(n)}) \rangle_\Sigma) & =      \text{AVar}_n(\langle A, D^M_\Sigma(\sqrt{n}S^{(n)}) \rangle_\Sigma) 
    \\
& =  \text{AVar}_n(\langle D_\Sigma^M(A), \sqrt{n}(S^{(n)} - \Sigma) \rangle_\Sigma) 
\\
& = 4 \Vert D_\Sigma^M(A) \Vert_\Sigma^2.
\end{align*}
Similarly,
\begin{align*}
       \text{AVar}_n(\sqrt{n}\langle A, P(S^{(n)}) \rangle_\Sigma) & =      \text{AVar}_n(\langle A, D^P_\Sigma(\sqrt{n}S^{(n)}) \rangle_\Sigma) 
    \\
& =  \text{AVar}_n(\langle D^M_\Sigma(A), D^P_\Sigma(\sqrt{n}S^{(n)}) \rangle_\Sigma) 
\\
& =  \text{AVar}_n(\langle (D^P_\Sigma)^*D^M_\Sigma(A), \sqrt{n}(S^{(n)} - \Sigma) \rangle_\Sigma) 
\\
& = 4\Vert (D^P_\Sigma)^*D^M_\Sigma(A) \Vert^2_\Sigma,
\end{align*}
where the second equality is a result of $D^M_\Sigma$ being the orthogonal projection onto $T_\Sigma\Skron$ and $D^T_\Sigma$ having range $T_\Sigma\Skron$. Defining $B = D^M_\Sigma(A)$, the maximum asymptotic variance ratio is
\begin{align*}
    \sup_{A \in \Sym} \frac{\text{AVar}_n(\langle A, P(S^{(n)}) \rangle_\Sigma )}{\text{AVar}_n(\langle A, M(S^{(n)}) \rangle_\Sigma)} & =  \sup_{B \in T_\Sigma \Skron} \frac{\Vert (D^P_\Sigma)^*(B)\Vert_\Sigma^2}{\Vert B \Vert^2_\Sigma} 
    \\
    & = \sup_{B \in T_\Sigma \Skron: \Vert B \Vert^2_\Sigma = 1} \Vert (D^P_\Sigma)^*(B)\Vert_\Sigma^2.
\\
& = \sup_{B \in T_\Sigma \Skron, C \in \Sym: \Vert B \Vert^2_\Sigma = 
 \Vert C \Vert^2_\Sigma = 1} \langle (D_\Sigma^P)^*(B),C \rangle_\Sigma^2 
    \\
    & = \sup_{B \in T_\Sigma \Skron, C \in \Sym: \Vert B \Vert^2_\Sigma = \Vert C \Vert^2_\Sigma = 1} \langle B ,D^P_\Sigma(C) \rangle_\Sigma^2 
    \\
    & =  \sup_{C \in \Sym: \Vert C \Vert^2_\Sigma = 1} \Vert D^P_\Sigma(C)\Vert_\Sigma^2 
\end{align*}
For any $C \in \Sym$, decompose $C = C_1 + C_2$ with $C_1 \in T_\Sigma \Skron$ and $C_2 \in T_\Sigma\mathcal{A}^P_\Sigma$ so that $D_\Sigma^P(C) = C_1$, giving
\begin{align*}
      \sup_{C \in \Sym: \Vert C \Vert^2_\Sigma = 1} \Vert D^P_\Sigma(C)\Vert_\Sigma^2 & = \sup_{C_1 \in T_\Sigma \Skron, C_2 \in T_\Sigma\mathcal{A}^T_\Sigma} \left\Vert D^P_\Sigma\left(\frac{C_1 + C_2}{\Vert C_1  + C_2 \Vert^2_\Sigma}  \right) 
      \right\Vert_\Sigma^2 
      \\
      & =  \sup_{C_1 \in T_\Sigma \Skron, C_2 \in T_\Sigma\mathcal{A}^T_\Sigma} \frac{\Vert C_1 \Vert^2_\Sigma}{\Vert C_1 + C_2 \Vert^2_\Sigma}
\end{align*}
We now argue that if $C_1$ and $C_2$ are fixed, up to a choice of scaling, then the maximum of $ \Vert C_1 \Vert^2_\Sigma\Vert C_1 + C_2 \Vert^{-2}_\Sigma$ is equal to $\sin(\Tilde{\theta}_{\ell_1,\ell_2})^{-2}$ where $\Tilde{\theta}_{\ell_1,\ell_2}$ is the angle under the FIM between the lines $\ell_1 = \text{span}\{C_1\}$ and $\ell_2 = \text{span}\{C_2\}$. The scaling of $C_1$ does not change the value of  $\Vert C_1 \Vert^2_\Sigma\Vert C_1 + C_2 \Vert^{-2}_\Sigma$ so without loss of generality we fix $C_1$ and minimize $\Vert C_1 + C_2 \Vert^2_\Sigma$ over $C_2$. Any point $C_1 + C_2$ necessarily lies on the line $\ell_2 + C_1$. The minimizing value of $C_1 + C_2$ is the FIM orthogonal projection of $0$ onto the affine subspace $\ell_2 + C_1$. The points $0,C_1,C_2,C_1 + C_2$ form a parallelogram with $\Tilde{\theta}_{\ell_1,\ell_2}$ equalling the angle between $\ell_1$ and $\ell_2$ at $0$, which is equal to the angle at $C_1$ between $\ell_1$ and $C_1 + \ell_2$. The triangle $\Delta 0(C_1+C_2)C_1$ is a right triangle with right angle at $C_1 + C_2$ as $C_1 + C_2$ is the orthogonal projection onto $\ell_2 + C_1$. The side opposite to this right angle has length $\Vert C_1 \Vert_\Sigma$, while the side opposite the angle $\Tilde{\theta}_{\ell_1,\ell_2}$ has length $\Vert C_1 + C_2 \Vert_\Sigma$, from which it follows that $\sin(\Tilde\theta) = \Vert C_1 + C_2 \Vert_\Sigma \Vert C_1 \Vert_\Sigma^{-1}$. This gives the desired result as
\begin{align*}
    \sup_{C_1 \in T_\Sigma \Skron, C_2 \in T_\Sigma\mathcal{A}^T_\Sigma} \frac{\Vert C_1 \Vert^2_\Sigma}{\Vert C_1 + C_2 \Vert^2_\Sigma} = \sup_{{C_1 \in T_\Sigma \Skron, C_2 \in T_\Sigma\mathcal{A}^T_\Sigma}} \sin(\Tilde{\theta}_{\ell_1,\ell_2})^{-2},
\end{align*}
is maximized when $\Tilde{\theta}_{\ell_1,\ell_2}$ is the smallest angle between lines $\ell_1$ and $\ell_2$ of $T_\Sigma \Skron$ and $T_\Sigma \mathcal{A}_\Sigma^P$, namely the principle angle between these two subspaces.
\end{proof}

 Reasoning similar to that above provides the following useful characterization of the asymptotic variance of efficient estimators in curved exponential families.

\begin{customlem}{0.4}
\label{lem:AsympVarianceCurvedFamilies}
Let $\Theta$ be an open set in the vector space $\mathbb{R}^k$ and take $\Theta_0$ to be a smooth submanifold of $\Theta$. Additionally, let $\mathcal{P}_0 = \{P_\theta:\theta \in \Theta_0\}$ be a curved exponential family contained in the regular exponential family $\mathcal{P} = \{P_\theta:\theta \in \Theta\}$ with Fisher information metric $\langle \cdot,\cdot\rangle_\theta$. If $T_n$ is the sufficient statistic based on $n$ independent observations from $P_{\theta_0}$, and $\varphi(T_n)$ is an smooth, efficient estimator of $\theta_0 \in \Theta$ in the submodel $\mathcal{P}_0$, the variance of the asymptotic normal distribution of $\langle v, \sqrt{n}(\varphi(T_n) - \varphi(\theta_0)\rangle_{\theta_0}$ is 
    $\Vert P_{\theta_0}v\Vert^2_{\theta_0}$, where $P_{\theta_0}$ is the orthogonal projection onto $T_{\theta_0}\Theta$.
\end{customlem}
\begin{proof}
We equip $\mathbb{R}^k$ with an arbitrary inner product $\langle\cdot,\cdot\rangle$, where the self-adjoint Fisher information matrix $\mathcal{I}(\theta)$ associated with the Fisher information metric satisfies
\begin{align*}
    \langle v,w\rangle_\theta = \langle v,\mathcal{I}(\theta)w\rangle =  \langle \mathcal{I}(\theta) v,w\rangle.
\end{align*}

As was argued in the proof of Theorem \ref{thm:AppendixAvarRatio}, if $\varphi$ is efficient and $\varphi(\theta_0) = \theta_0$ for all $\theta_0 \in \Theta_0$,  the differential $D^\varphi_{\theta_0}$ of $\varphi$ at $\theta_0$ is equal to $P_{\theta_0}$. 
The distribution of $\sqrt{n}(T_n - \theta_0)$ converges weakly to $N_k(0,\mathcal{I}(\theta_0)^{-1})$, where if $Y \sim N_k(0,\mathcal{I}(\theta_0)^{-1})$ then $\mathcal{I}(\theta_0)^{-1}$ is the matrix that satisfies \cite[Ch1]{eatonmultivariate}
\begin{align}
\label{eqn:VarianceDefinition}
    \text{Var}(\langle v,Y\rangle)  = \langle v, \mathcal{I}(\theta_0)^{-1}v\rangle
\end{align}
By the delta method and the fact that $P_{\theta_0}$ is self-adjoint with respect to the FIM, the asymptotic distribution of $\langle v, \sqrt{n}(\varphi(T_n) - \varphi(\theta_0))\rangle_{\theta_0}$ equals the asymptotic distribution of 
\begin{align*}
    \langle v, \sqrt{n}P_{\theta_0}(T_n - \theta_0)\rangle_{\theta_0} & =   \langle P_{\theta_0}v, \sqrt{n}(T_n - \theta_0)\rangle_{\theta_0} 
    \\
    & = \langle \mathcal{I}(\theta_0)P_{\theta_0}v,\sqrt{n}(T_n - \theta_0)\rangle.
\end{align*}
Applying \eqref{eqn:VarianceDefinition}, the variance of the asymptotic distribution of the above quantity is
\begin{align*}
    \langle  \mathcal{I}(\theta_0)P_{\theta_0}v, \mathcal{I}(\theta_0)^{-1} \mathcal{I}(\theta_0)P_{\theta_0}v\rangle = \Vert P_{\theta_0}v\Vert^2_{\theta_0}.
\end{align*}
\end{proof}

\begin{customlem}{0.5}
\label{lem:ReparamVariance}
Define $\sqrt{n}(\varphi(T_n) - \varphi(\theta))$ as in Lemma \ref{lem:AsympVarianceCurvedFamilies}. 
    Let $\phi:\Theta_0 \rightarrow \Gamma$ be a smooth reparameterization of $\Theta_0$, with the associated estimator sequence $\sqrt{n}(\phi(\varphi(T_n)) - \phi(\varphi(\theta)))$. If the manifold $\Gamma$ is equipped with the push-forward Riemannian metric defined by
    \begin{align*}
        \langle D^\phi_{\theta_0}(v),D^\phi_{\theta_0}(w) \rangle_{\phi(\theta_0)}^* \coloneqq \langle v, w\rangle_{\theta_0}, \;\; v,w \in T_{\theta_0}\Theta_0
    \end{align*}
    then the asymptotic variance of $\langle D^\phi_{\theta_0}(v),\sqrt{n}(\phi(\varphi(T_n)) - \phi(\varphi(\theta)))\rangle_{\phi(\theta_0)}^*$ is $\Vert v \Vert_{\theta_0}^2$.
\end{customlem}
\begin{proof}
    The delta method shows that the asymptotic variance of $\langle D^\phi_{\theta_0}(v),\sqrt{n}(\phi(\varphi(T_n)) - \phi(\varphi(\theta)))\rangle_{\varphi(\theta_0)}^*$ equals that of
    \begin{align*}
        \langle D^\phi_{\theta_0}(v),\sqrt{n}D^\phi_{\theta_0}(\varphi(T_n) - \varphi(\theta) )\rangle_{\phi(\theta_0)}^* = \langle v, \sqrt{n}(\varphi(T_n) - \varphi(\theta)) \rangle_{\theta_0} = \Vert v \Vert^2_{\theta_0},
    \end{align*}
    where the last equality follows from Lemma \ref{lem:AsympVarianceCurvedFamilies}. We note that as $\phi$ is a map with domain $\Theta_0$, the differential $D^\phi_{\theta_0}$ is only defined on vectors that lie in the tangent space $T_{\theta_0}\Theta_0$. We extend $D^\phi_{\theta_0}$ so that it takes the value $0$ on $(T_{\theta_0}\Theta_0)^\perp$. This ensures that the quantity $\sqrt{n}D^\phi_{\theta_0}(\varphi(T_n) - \varphi(\theta))$ is well-defined.
\end{proof}

\begin{customlem}{2}
The functions, $M$ and $P$ are $\Glkron$ and $\Okron$ equivariant respectively, implying that
the auxiliary spaces $\mathcal{A}_{\Sigma_1 \otimes \Sigma_2}^M$ are $\Glkron$ equivariant while the spaces $\mathcal{A}_{\Sigma_1 \otimes \Sigma_2}^P$ are $\Okron$ equivariant: 
 \begin{align*}
     (A \otimes B)\mathcal{A}_{\Sigma_1 \otimes \Sigma_2}^M(A \otimes B)^\top & = \mathcal{A}_{A\Sigma_1A^\top \otimes B\Sigma_2 B^\top}^M,\; \text{for all}\; A \otimes B \in \Glkron
     \\
     (U \otimes V)\mathcal{A}_{\Sigma_1 \otimes \Sigma_2}^P M(U \otimes V)^\top & = \mathcal{A}_{U\Sigma_1 U^\top \otimes V\Sigma_2 V^\top}^P, \;\text{for all}\; U \otimes V \in \Okron. 
 \end{align*}
 In particular, if $\Sigma_1 \otimes \Sigma_2$ has the eigendecomposition $U\Lambda U^\top \otimes V\Gamma V^\top$, all of the principle angles between $T_{\Sigma_1 \otimes \Sigma_2} \mathcal{A}^P_{\Sigma_1 \otimes \Sigma_2}$ and $T_{\Sigma_1 \otimes \Sigma_2}\Skron$ with respect to $\langle \cdot ,\cdot \rangle_{\Sigma_1 \otimes \Sigma_2}$ are identical to the principle angles between $T_{\Lambda \otimes \Gamma}\mathcal{A}^P_{\Lambda \otimes \Gamma}$ and $T_{\Lambda \otimes \Gamma}\Skron$ with respect to $\langle \cdot, \cdot \rangle_{\Lambda \otimes \Gamma}$.
\end{customlem}
\begin{proof}
 The $\Glkron$ equivariance of $M$ follows from the likelihood equations \eqref{eqn:ComputationalMLELikelihood1} and \eqref{eqn:ComputationalMLELikelihood2}. If $\hat{\Sigma}_1 \otimes \hat{\Sigma}_2$ satisfies
    \begin{align*}
        \tr_1\big((\hat{\Sigma}_1 \otimes \hat{\Sigma}_2)^{-1}S\big) = p_2 I_{p_1}, \;\;   \tr_2\big((\hat{\Sigma}_1 \otimes \hat{\Sigma}_2)^{-1}S\big) = p_1 I_{p_2},
    \end{align*}
    then using \eqref{eqn:CyclicPerm} and \eqref{eqn:PTEquivariance}  if $A \otimes B \in \Glkron$, 
   \begin{align*}
        & \tr_1\big((A\hat{\Sigma}_1A^\top \otimes B\hat{\Sigma}_2B^\top)^{-1}(A \otimes B)S(A \otimes B)^\top\big) = p_2 I_{p_1},
        \\
        &  \tr_2\big((A\hat{\Sigma}_1A^\top \otimes B\hat{\Sigma}_2B^\top)^{-1}(A \otimes B)S(A \otimes B)^\top\big) = p_1 I_{p_2}.
   \end{align*} 
   The equivariance claim $P\left((U \otimes V)S(U \otimes V)^\top\right) = (U \otimes V)P(S)(U \otimes V)^\top$ also  follows from \eqref{eqn:CyclicPerm} and \eqref{eqn:PTEquivariance} since for $U \otimes V \in \Okron$
   \begin{align*}
       \tr_1\big( (U \otimes V)S(U \otimes V)^\top\big) = U \tr_1\big( (I_{p_1} \otimes V^\top)(I_{p_1} \otimes V)S \big)U^\top = U \tr_1(S) U^\top.
   \end{align*}
A similar equation holds for $\tr_2$, yielding
\begin{align*}
    P\big( (U \otimes V)S(U \otimes V)^\top \big) & = \frac{\tr_1\left((U \otimes V)S(U \otimes V)^\top\right) \otimes \tr_2\left((U \otimes V)S(U \otimes V)^\top\right)}{\tr\big((U \otimes V)S(U \otimes V)^\top\big)} 
    \\
    & = \frac{ U\tr_1(S)U^\top \otimes V\tr_2(S)V^\top }{\tr\left(S\right)} 
    \\
    & = (U \otimes V) P(S) (U \otimes V)^\top.
\end{align*}  
The corresponding equivariance of $\mathcal{A}^P_{\Sigma_1 \otimes \Sigma_2}$ and $\mathcal{A}^M_{\Sigma_1 \otimes \Sigma_2}$ holds by the definition of an auxiliary space. 

   To prove the last statement regarding principle angles, consider the map on that sends $H \in T_{\Sigma_1 \otimes \Sigma_2}\Sp$ to $(U \otimes V)^\top H (U \otimes V) \in T_{\Lambda \otimes \Gamma}\Sp$. By the $\Okron$ equivariance of the linear space $T_{\Sigma_1 \otimes \Sigma_2}\mathcal{A}^P_{\Sigma_1 \otimes \Sigma_2}$, the image of this subpace under this map equals $T_{\Lambda \otimes \Gamma}\mathcal{A}^P_{\Lambda \otimes \Gamma}$. Likewise, the image of $T_{\Sigma_1 \otimes \Sigma_2}\Skron$ is $T_{\Lambda \otimes \Gamma}\Skron$. This map is an inner product space isometry between the inner product spaces $(\Sym,\langle \cdot ,\cdot \rangle_{\Sigma_1 \otimes \Sigma_2})$ and $(\Sym,\langle \cdot ,\cdot \rangle_{\Lambda \otimes \Gamma})$ by \eqref{eqn:AIMetricIsom}. Principle angles, being defined in terms of a sequential maximization of inner products, are preserved by vector space isometries, completing the proof.
\end{proof}

\begin{customlem}{3}
Denote the eigenvalues of $\Sigma_1$ and $\Sigma_2$ by $\lambda = (\lambda_1,\ldots,\lambda_{p_1})$ and $\gamma = (\gamma_1, \ldots, \gamma_{p_2})$ respectively. The maximal asymptotic variance ratio of the partial trace estimator is bounded below, as follows
    \begin{align}
           \sup_{A \in \Sym} \frac{\text{AVar}_n(\langle A, P(S^{(n)}) \rangle_{\Sigma_1 \otimes \Sigma_2} )}{\text{AVar}_n(\langle A, M(S^{(n)}) \rangle_{\Sigma_1 \otimes \Sigma_2} )} \geq \max\bigg(  \frac{\Vert \lambda \Vert^2}{(\lambda^\top 1_{p_1}/\sqrt{p_1})^2},  \frac{\Vert \gamma \Vert^2}{(\gamma^\top 1_{p_2}/\sqrt{p_2})^2}  \bigg).
    \end{align}
    If $\alpha_1$ and  $\alpha_2$  are the angles between $\lambda, 1_{p_1}$ and $\gamma, 1_{p_2}$ then the lower bound \eqref{eqn:VarRatioBoundEigenval1} can be written as
    \begin{align}
          \sup_{A \in \Sym} \frac{\text{AVar}_n(\langle A, P(S^{(n)}) \rangle_{\Sigma_1 \otimes \Sigma_2} )}{\text{AVar}_n(\langle A, M(S^{(n)}) \rangle_{\Sigma_1 \otimes \Sigma_2} )} \geq \max\big(\cos(\alpha_1)^{-2}, \cos(\alpha_2)^{-2}\big).
    \end{align}
\end{customlem}
\begin{proof}
Lemma \ref{lem:PrincipleAngleEquivariance} shows that without loss of generality $\Sigma_1 \otimes \Sigma_2$ can be taken to be the diagonal matrix $\Lambda \otimes \Gamma$. That is, the maximum asymptotic variance ratio only depends on the eigenvalues of $\Sigma_1 \otimes \Sigma_2$. To ease notation, let $\Sigma = \Lambda \otimes \Gamma$. We also let  $\circ$ represent the Hadamard product $a \circ b = (a_1b_1,\ldots,a_kb_k)$ and define the element-wise division of vectors $a/b = (a_1/b_1,\ldots,a_k/b_k)$ for $a,b \in \mathbb{R}^k$.

Take $x \in \mathbb{R}^{p_1}$, $y \in \mathbb{R}^{p_2}$ with $1^\top x = 1^\top y = 0$. We define 
\begin{align*}
    V = \text{diag}(x) \otimes \text{diag}(y),
\end{align*}
where by construction $\tr_1(V) = 0$ and $\tr_2(V) = 0$ so that $V \in T_\Sigma\mathcal{A}^T_{\Sigma}$. The norm of $V$ can be represented as
\begin{align*}
    \Vert V \Vert_\Sigma^2 = \tr\big( (\text{diag}(x/\lambda) \otimes \text{diag}(y/\gamma))^2\big) = \Vert x/\lambda \Vert^2 \Vert y/\gamma \Vert^2,
\end{align*}
where the norm on the right hand side is the usual Euclidean norm of vectors.

We want to find a vector $H \in T_\Sigma \Skron$ that has a large inner product with $V$. Assume that $H$ has the form $H = \Lambda \otimes H_2$. Then
\begin{align*}
    \Vert H \Vert_\Sigma^2 = \tr\big((I_{p_1} \otimes \Gamma^{-1/2}H_2\Gamma^{-1/2})^2\big) = p_1 \tr\big((\Gamma^{-1/2}H_2\Gamma^{-1/2})^2\big).
\end{align*}
If $\Vert H \Vert^2_\Sigma = 1$ then $\tr\big((\Gamma^{-1/2}H_2\Gamma^{-1/2})^2\big)^{1/2} = p_1^{-1/2}$. We compute
\begin{align*}
    \langle V,H \rangle_\Sigma & = \tr\big(  \text{diag}(x/\lambda) \otimes \text{diag}(y/\gamma) (I_{p_1} \otimes \Gamma^{-1/2}H_2\Gamma^{-1/2})\big)
    \\
    & = (1^\top (x/\lambda)) \tr(\text{diag}(y/\gamma) \Gamma^{-1/2}H_2\Gamma^{-1/2}).
\end{align*}
The optimal choice of $H_2$ maximizing the above inner product, subject to the scale constraint, is the matrix 
\begin{align*}
    H_2 = \frac{1}{\Vert y/\gamma \Vert p_1^{1/2}}\Gamma^{1/2} \text{diag}(y/\gamma)\Gamma^{1/2}.
\end{align*}
With this choice of $H_2$, $\Vert H \Vert^2_\Sigma = 1$ and 
\begin{align*}
  \langle V,H\rangle_\Sigma =  ( 1^\top (x/\lambda)) \Vert y/\gamma \Vert p_1^{-1/2},
\end{align*}
so that
\begin{align}
\label{eqn:AngleOptimizationEq1}
    \left\langle \frac{V}{\Vert V\Vert},\frac{H}{\Vert H\Vert} \right\rangle_\Sigma = \frac{(1^\top (x/\lambda))}{\Vert x/\lambda \Vert p_1^{1/2}}. 
\end{align}
Now we optimize this expression with respect to $x$. As this expression is scale invariant we assume that $\Vert x/\lambda \Vert^2 = 1$. Make a change of variables $z = x/\lambda$, or equivalently $z \circ \lambda = x$. The partial trace constraint $1^\top x = 0$ can be rewritten as $\lambda^\top z = 1$. The optimization problem becomes
\begin{align*}
    \max_{z} \;1^\top z, \text{ such that } z^\top z = 1 \text{ and } \lambda^\top z = 0.
\end{align*}
Using Lagrange multipliers the optimal $z$ equals $a(1 - \tfrac{\lambda^\top 1}{\Vert \lambda\Vert^2} \lambda)$ where $a > 0$ ensures that this expression has unit norm. The optimal $x$ is therefore $x = a(\lambda - \frac{\lambda^\top 1 }{\Vert \lambda\Vert^2}\lambda \circ \lambda)$. Plugging this optimal $x$ back into \eqref{eqn:AngleOptimizationEq1} yields
\begin{align*}
      \left\langle \frac{V}{\Vert V\Vert},\frac{H}{\Vert H\Vert} \right\rangle_\Sigma = \frac{1^\top (1 - \tfrac{\lambda^\top 1}{\Vert \lambda\Vert^2} \lambda)}{\Vert(1 - \tfrac{\lambda^\top 1}{\Vert \lambda\Vert^2} \lambda)\Vert p_1^{1/2}} = \left(1 -  \frac{(\lambda^\top 1)^2}{p_1\Vert \lambda \Vert^2}\right)^{1/2}.
\end{align*}
If we take $\Tilde{\theta}$ to be the angle between $V$ and $H$ and $\theta$ to be the principle angle between $T_\Sigma\mathcal{A}^T_\Sigma$ and $T_\Sigma \Skron$ then 
\begin{align*}
\sin(\theta)^{-2} \geq  \sin(\Tilde{\theta})^{-2} = \bigg(1 -   \big\langle \frac{V}{\Vert V\Vert},\frac{H}{\Vert H\Vert} \big\rangle_\Sigma^2\bigg)^{-1} =  \frac{p_1 \Vert \lambda \Vert^2}{(\lambda^\top 1)^2} = \frac{\Vert \lambda \Vert^2}{(\lambda^\top 1/\sqrt{p_1})^2}
\end{align*}
Note that $\tfrac{(\lambda^\top 1/\sqrt{p_1})^2 }{ \Vert \lambda \Vert^2}$ is equal to $\cos(\alpha_1)^2$.

Symmetric equations hold with respect to $\gamma$, where in this case $H$ has the form $H_1 \otimes \Gamma$ and the subsequent maximization is performed over $y$.
 \end{proof}

\begin{customthm}{2}
       The function $R$ is $\Glkron$ equivariant. If $nS^{(n)} \sim \text{Wishart}_p(\Sigma_1 \otimes \Sigma_2)$, the estimator $R(S^{(n)})$ is asymptotically efficient.
\end{customthm}
\begin{proof}
In showing the equivariance of $R(S)$, note that
\begin{align*}
        \big((A \otimes B)S(A\otimes B)^\top\big)_{\cdot j,\cdot j} = (B_{j,j})^2 (AS_{\cdot j,\cdot j}A^\top)
\end{align*}
    The equivariance of $R(S)$ follows from the computation
    \begin{align*}
        \tr_1^{\det}((A\otimes B)S(A\otimes B)^\top) & =  \sum_{j = 1}^{p_2} \frac{(A\otimes B)S(A\otimes B)^\top)_{\cdot j, \cdot j}}{\vert (A\otimes B)S(A\otimes B)^\top)_{\cdot j,\cdot j}\vert^{1/p_1}}
        \\
        & =  \sum_{j = 1}^{p_2} \frac{(B_{i,i})^2 (AS_{\cdot j,\cdot j}A^\top)}{\vert (B_{i,i})^2 (AS_{\cdot j,\cdot j}A^\top) \vert^{1/p_1}}
        \\
        & = \frac{A}{\vert A\vert^{1/p_1}}\bigg(\sum_{j = 1}^{p_2} \frac{ S_{\cdot j,\cdot j}}{\vert S_{\cdot j,\cdot j} \vert^{1/p_1}}\bigg)\frac{A^\top}{\vert A \vert^{1/p_1}}.
    \end{align*}
    Symmetrically, $\tr_2^{\det}((A\otimes B)S(A\otimes B)^\top) = \vert B\vert^{-2/p_2} B \tr_2^{\det}(S)B^\top$. In particular, 
    \begin{align*}
        \vert \tr_i^{\det}((A\otimes B)S(A\otimes B)^\top) \vert = \vert \tr_i^{\det}(S) \vert.
    \end{align*}
    This proves that $R$ is $\Glkron$-equivariant:
    \begin{align*}
        R((A\otimes B)S(A\otimes B)^\top)) = & \bigg(\frac{\vert ((A\otimes B)S(A\otimes B)^\top) \vert^{1/p} }{ \vert \tr_1^{\det}((A\otimes B)S(A\otimes B)^\top)) \otimes \tr_2^{\det}((A\otimes B)S(A\otimes B)^\top)) \vert^{1/p} }\bigg) \times
        \\
        & \tr_1^{\det}((A\otimes B)S(A\otimes B)^\top) \otimes \tr_2^{\det}((A\otimes B)S(A\otimes B)^\top)
        \\
        = & \bigg(\frac{\vert S \vert \vert A \vert^{2p_2} \vert B \vert^{2p_1} }{ \vert \tr_1^{\det}(S) \otimes \tr_2^{\det}(S) \vert}\bigg)^{1/p} \frac{A\tr_1^{\det}(S)A^\top  \otimes B\tr_2^{\det}(S)B^\top }{\vert A \vert^{2/p_1}\vert B \vert^{2/p_2}}
        \\
        & = (A \otimes B) R(S) (A \otimes B)^\top.
    \end{align*}
The map $R:\Sp \rightarrow \Skron$ is a smooth, surjective map and by equivariance it has constant rank \cite[Thm 4.14 a)]{lee2013smooth}.

To prove that $R(S)$ is asymptotically efficient, by equivariance it is enough to show that the auxiliary space of $R(S)$ at the identity $\mathcal{A}^{R}_{I_p}$ is orthogonal to the tangent space $T_{I_p}\Skron$, meaning that $\mathcal{A}^R_{I_p} = (T_{I_p}\Skron)^\perp$. As $R$ is a submersion, the dimension of $\mathcal{A}^R_{I_p}$ equals the dimension of $(T_{I_p}\Skron)^\perp$. Consequently, if it is shown that  $(T_{I_p}\Skron)^\perp \subset \mathcal{A}^R_{I_p}$ then this containment will imply the desired orthogonality property. 

Let $S(t)$ be a curve in $\Sp$ with $S(0) = I_p$ and $S'(0) = \Delta \in (T_{I_p}\Skron)^\perp$. For an arbitrary $\Delta$, we want to show that $\tfrac{d}{dt}R(S(t))\big\vert_{t = 0} = 0$, which implies the containment  $(T_{I_p}\Skron)^\perp \subset \mathcal{A}^R_{I_p}$. As 
\begin{align*}
    (T_{I_p}\Skron)^\perp = \mathcal{A}^M_{I_p} = \{\Delta:\tr_i(\Delta) = 0, \; i = 1,2\},
\end{align*}
assume that $\tr_1(\Delta) = 0$ and $\tr_2(\Delta) = 0$. 

To ease notation, form the matrices
\begin{align*}
    D_1 & =  \text{diag}\big( \vert S_{1 \cdot ,1\cdot} \vert^{-1/p_2}, \ldots, \vert S_{p_1 \cdot, p_1 \cdot} \vert^{-1/p_2} \big)
    \\
        D_2 & = \text{diag}\big( \vert S_{\cdot 1,\cdot 1} \vert^{-1/p_1}, \ldots, \vert S_{\cdot p_2,\cdot p_2} \vert^{-1/p_1} \big),
\end{align*}
so that
\begin{align*}
    \tr_1^{\det}(S) &= \tr_1( (I_{p_1} \otimes D_2) S),
    \\
    \tr_2^{\det}(S) &= \tr_2\big( (D_1 \otimes I_{p_2})S\big).
\end{align*}
By the product rule we have
\begin{align*}
    \frac{d}{dt} \tr_1^{\det}\big(S(t)\big) \big\vert_{t = 0} & = \tr_1\big( (I_{p_1} \otimes D_2'(0))S(0)\big) + \tr_1\big((I_{p_1} \otimes D_2(0))S'(0)\big)
    \\
    & = \tr_1\big(I_{p_1} \otimes D_2'(0)\big) + \tr_1\big(S'(0)\big)
    \\
    & = p_1 D_2'(0) + \tr_1(\Delta) 
    \\
    & = p_1 D'_2(0).
\end{align*}
The derivatives of the diagonal entries of $D_2$ are
\begin{align*}
    \frac{d}{dt} \vert S_{\cdot j,\cdot j}(t)\vert^{-1/p_1} \bigg\vert_{t = 0} & = \bigg(-\tfrac{1}{p_1}  \vert S_{\cdot j,\cdot j}(0)\vert^{-1/p_1 - 1}\bigg) \frac{d}{dt}  \vert S_{\cdot j,\cdot j}(t)\vert \bigg\vert_{t = 0}
    \\
    & =  \frac{d}{dt} \vert S_{\cdot j,\cdot j}(t)\vert \bigg\vert_{t = 0}
    \\
    & =   \vert S_{\cdot j,\cdot j}(0)\vert \tr\big(  (S_{\cdot j,\cdot j})^{-1}(0) (S_{\cdot j,\cdot j})'(0)\big)
    \\
    & = \tr( \Delta_{\cdot j,\cdot j})
    \\
    & = 0,
\end{align*}
where the second equality follows from Jacobi's formula and the last equality follows from $\tr_2(\Delta) = 0$. We conclude that $D_2'(0) = 0$ and the derivative of $\tr_1^{\det}\big(S(t)\big)$ is $0$ at $t = 0$. Similarly, $\tfrac{d}{dt}\tr_2^{\det}\big(S(t)\big)\vert_{t = 0 } = 0$. It remains to find the derivative of the multiplicative constant 
\begin{align*}
\bigg(\frac{\vert S(t) \vert }{\vert \tr_1^{\det}\big(S(t)\big) \otimes \tr_2^{\det}\big(S(t)\big) \vert}\bigg)^{1/p} = \vert S(t) \vert^{1/p} \vert \tr_1^{\det}\big(S(t)\big) \otimes \tr_2^{\det}\big(S(t)\big) \vert^{-1/p}.
\end{align*}
From what we have just shown with respect to $\tr_1^{\det}$ and $\tr_2^{\det}$, the derivative of the right-most factor above will be $0$. Thus, by the product rule
\begin{align*}
    \frac{d}{dt}\bigg(\frac{\vert S(t) \vert }{\vert \tr_1^{\det}\big(S(t)\big) \otimes \tr_2^{\det}\big(S(t)\big) \vert}\bigg)^{1/p} \bigg\vert_{t = 0} & = \frac{d}{dt} \vert S(t) \vert^{1/p} \bigg\vert_{t = 0} \vert \tr_1^{\det}\big(S(0)\big) \otimes \tr_2^{\det}\big(S(0)\big) \vert^{-1/p}
    \\
    & = \frac{d}{dt} \vert S(t) \vert^{1/p} \bigg\vert_{t = 0}
    \\
    & = \tfrac{1}{p} \vert S(0) \vert^{1/p} \tr\big( S^{-1}(0) S'(0)\big)
    \\
    & = \tr(\Delta)
    \\
    & = \tr\big(\tr_1(\Delta)\big) = 0.
\end{align*}
Recalling that for Kronecker products, $\tfrac{d}{dt} \big(A(t) \otimes B(t)\big) \big\vert_{t = 0} = A'(0) \otimes B(0) + A(0) \otimes B'(0)$, 
one last application of the product rule shows that 
\begin{align*}
    \frac{d}{dt} R\big(S(t)\big) \bigg\vert_{t = 0} = 0,
\end{align*}
completing the proof.
\end{proof}

\begin{customthm}{3}
    The map $\varphi: \Skron \rightarrow \mathbb{R} \times \mathbb{P}(\mathcal{S}^{p_1}_{++}) \times \mathbb{P}(\mathcal{S}^{p_2}_{++})$, with $\varphi^{-1}$ defined as $\varphi^{-1}(c,\Tilde \Sigma_1, \Tilde \Sigma_2) = e^c(\Tilde\Sigma_1 \otimes \Tilde\Sigma_2)$ is a bijection that is an orthogonal parameterization. The Riemannian metric on the Cartesian product $\mathbb{R} \times \mathbb{P}(\mathcal{S}^{p_1}_{++}) \times \mathbb{P}(\mathcal{S}^{p_2}_{++})$ induced from the FIM on $\Skron$ is a product metric, where $\mathbb{R}$ has the standard Euclidean metric scaled by a factor of $p_1p_2$ and each $\mathbb{P}(\mathcal{S}^{p_i}_{++}), \; i = 1,2$ is equipped with the restriction of the FIM on $\mathcal{S}^{p_i}_{++}$ scaled by a factor of $p_j$, $j \neq i$. 
\end{customthm}
\begin{proof}
    Let $\tfrac{d}{dt}\big\vert_{t = 0}(c(t), \Tilde{\Sigma}_1(t),\Tilde{\Sigma}_2(t)) =(v,H_1,H_2)$ be a tangent vector of $\mathbb{R} \times \mathbb{P}(\mathcal{S}^{p_1}_{++}) \times \mathbb{P}(\mathcal{S}^{p_2}_{++})$ at the point $(c(0), \Tilde{\Sigma}_1(0),\Tilde{\Sigma}_2(0)) = (c,\Tilde{\Sigma}_1,\Tilde{\Sigma}_2)$. We define $\Sigma = \Sigma_1 \otimes \Sigma_2 =  \varphi^{-1}(c,\Tilde{\Sigma}_1,\Tilde{\Sigma}_2)$. The determinant one constraint on the spaces $\mathbb{P}(\mathcal{S}^{p_i}_{++})$ implies that each $H_i$ has trace $0$. The tangent vector $D\varphi^{-1}(v,H_1,H_2)$ in $T_\Sigma \Skron$ is
    \begin{align*}
        \frac{d}{dt}\bigg\vert_{t = 0} \varphi^{-1}(c(t), \Tilde{\Sigma}_1(t),\Tilde{\Sigma}_2(t)) & =       \frac{d}{dt}\bigg\vert_{t = 0} e^{c(t)} \Tilde{\Sigma}_1(t) \otimes \Tilde{\Sigma}_2(t) 
        \\
        & = v(\Sigma_1 \otimes \Sigma_2) +  e^c(H_1 \otimes \Tilde{\Sigma}_2) + e^c(\Tilde{\Sigma}_1 \otimes H_2).
    \end{align*}
    The induced Riemannian metric on $T_{(c,\Tilde{\Sigma}_1,\Tilde{\Sigma}_2)}\mathbb{R} \times \mathbb{P}(\mathcal{S}^{p_1}_{++}) \times \mathbb{P}(\mathcal{S}^{p_2}_{++})$ is 
    \begin{align*}
        \Vert (v,H_1,H_2)\Vert_{(c,\Tilde{\Sigma}_1,\Tilde{\Sigma}_2)}^2  \coloneqq & \Vert  v(\Sigma_1 \otimes \Sigma_2) +  e^c(H_1 \otimes \Tilde{\Sigma}_2) + e^c(\Tilde{\Sigma}_1 \otimes H_2) \Vert_{\Sigma}^2
        \\
         = & \tr\bigg(  \left(v(\Sigma_1 \otimes \Sigma_2) +  e^c(H_1 \otimes \Tilde{\Sigma}_2) + e^c(\Tilde{\Sigma}_1 \otimes H_2) \right)\Sigma^{-1}
        \\
        &\left( v(\Sigma_1 \otimes \Sigma_2) +  e^c(H_1 \otimes \Tilde{\Sigma}_2) + e^c(\Tilde{\Sigma}_1 \otimes H_2) \right) \Sigma^{-1}\bigg)
        \\
        = & v^2 \tr(I_{p_1} \otimes I_{p_2}) + \tr( \Tilde{\Sigma}_1^{-1}H_1\Tilde{\Sigma}_1^{-1}H_1 \otimes I_{p_2}) + \tr( I_{p_1} \otimes \Tilde{\Sigma}_2^{-1}H_2\Tilde{\Sigma}_2^{-1}H_2)
        \\
         = & p_1p_2 v^2 + p_2 \Vert H_1 \Vert^2_{\Tilde{\Sigma}_1} + p_1\Vert H_2 \Vert^2_{\Tilde{\Sigma}_2},
    \end{align*}
    where the third equality uses $\tr(H_i) = 0$ to show that the three terms are orthogonal to each other. By polarization, the above formula shows that the induced inner product on  $T_{(c,\Tilde{\Sigma}_1,\Tilde{\Sigma}_2)}\mathbb{R} \times \mathbb{P}(\mathcal{S}^{p_1}_{++}) \times \mathbb{P}(\mathcal{S}^{p_2}_{++})$ is exactly the inner product on the stated product Riemannian manifold $\mathbb{R} \times \mathbb{P}(\mathcal{S}^{p_1}_{++}) \times \mathbb{P}(\mathcal{S}^{p_2}_{++})$.
\end{proof}

\begin{customlem}{0.6}
    \label{lem:OrthogParamEquivariance}
    Define a $\Glkron$-action on $\mathbb{R} \times \mathbb{P}(\mathcal{S}^{p_1}_{++}) \times \mathbb{P}(\mathcal{S}^{p_2}_{++})$ by 
    \begin{align*}
       (A \otimes B)\cdot(c,\Tilde\Sigma_1,\Tilde\Sigma_2) \coloneqq  \bigg(c + \frac{2}{p}\log(\vert A \otimes B \vert) ,  \frac{A\Tilde\Sigma_1 A^\top}{\vert A \vert^{2/p_1}}, \frac{ B\Tilde\Sigma_2 B^\top}{\vert B \vert^{2/p_2}}\bigg). 
    \end{align*}
    The parameterization map $\varphi:\mathbb{R} \times \mathbb{P}(\mathcal{S}^{p_1}_{++}) \times \mathbb{P}(\mathcal{S}^{p_2}_{++}) \rightarrow \Skron$, or equivalently $\varphi^{-1}$, is $\Glkron$-equivariant with respect to the above action and the conjugation action on $\Skron$. The map $\phi: \Glkron \rightarrow \mathbb{R} \times \text{SL}_{p_1} \times \text{SL}_{p_2}$ defined by 
    \begin{align*}
        \phi( A \otimes B) = \bigg( \tfrac{1}{p}\log(\vert A \otimes B \vert), \frac{A}{\vert A \vert^{1/p_1}}, \frac{B}{\vert B \vert^{1/p_2}}\bigg)
    \end{align*}
    is a group isomorphism that extends the parameterization map $\varphi$. 
\end{customlem}
\begin{proof}
    The equivariance of $\varphi^{-1}$ follows from  
    \begin{align*}
        \varphi^{-1} \bigg(c + \frac{2}{p}\log(\vert A \otimes B \vert) ,  \frac{A\Tilde\Sigma_1 A^\top}{\vert A \vert^{2/p_1}}, \frac{ B\Tilde\Sigma_2 B^\top}{\vert B \vert^{2/p_2}}\bigg) & = e^c \vert A \otimes B \vert^{2/p}\bigg(  \frac{A\Tilde\Sigma_1 A^\top}{\vert A \vert^{2/p_1}} \otimes \frac{ B\Tilde\Sigma_2 B^\top}{\vert B \vert^{2/p_2}}\bigg)
        \\
        & = (A \otimes B) \big(e^c \Tilde\Sigma_1 \otimes \Tilde\Sigma_2)(A \otimes B)^\top
        \\
        & = (A \otimes B) \cdot \varphi^{-1}(c,  \Tilde\Sigma_1, \Tilde\Sigma_2).
    \end{align*}
    Any $A \otimes B \in \Glkron$ can be uniquely decomposed as 
\begin{align*}
    e^{\log(\vert A \otimes B\vert^{1/p})} \left( \frac{A}{\vert A \vert^{1/p_1} } \otimes \frac{B}{\vert B \vert^{1/p_2}} \right),
\end{align*}
with $A \vert A\vert^{-1/p_1} \in \text{SL}_{p_1}$ and $B \vert B\vert^{-1/p_2} \in \text{SL}_{p_2}$, showing that $\phi$ is bijective.

 That $\phi$ is a homomorphism is seen from
    \begin{align*}
        \phi((A \otimes B)(C \otimes D)) & =  \bigg( \tfrac{1}{p}\log(\vert AC \otimes BD \vert), \frac{AC}{\vert AC \vert^{1/p_1}}, \frac{BD}{\vert BD \vert^{1/p_2}}\bigg)
        \\
        & =   \bigg( \tfrac{1}{p}\log(\vert A \otimes B \vert) + \tfrac{1}{p}\log(\vert C \otimes D \vert), \frac{A}{\vert A \vert^{1/p_1}}\frac{C}{\vert C \vert^{1/p_1}} , \frac{B}{\vert B \vert^{1/p_2}} \frac{D}{\vert D \vert^{1/p_2}}\bigg).  
    \end{align*}
\end{proof}

\begin{customcor}{1}
The map $\varphi_1: \mathcal{S}^{p_1}_{++} \rightarrow \mathbb{R} \times \mathbb{P}(\mathcal{S}^{p_1}_{++})$, with $\varphi_1^{-1}$  defined as $\varphi_1^{-1}(c,\Tilde{\Sigma}_1) = e^c \Tilde{\Sigma}_1$ is an orthogonal parameterization of the Wishart model $\mathcal{P}_1 = \{\text{Wishart}_{p_1}(\Sigma_1,p_2):\Sigma_1 \in \mathcal{S}^{p_1}_{++}\}$. The FIM induced on $\mathbb{R} \times \mathbb{P}(\mathcal{S}^{p_1}_{++})$ by $\varphi_1$ and $\mathcal{P}_1$ equals the FIM induced on $\mathbb{R} \times \mathbb{P}(\mathcal{S}^{p_1}_{++})$ by the Kronecker model $\mathcal{P} = \{ \text{Wishart}_{p}(\Sigma_1 \otimes \Tilde{\Sigma}_2,1): \vert \Tilde{\Sigma}_2 \vert = 1\}$. Equivalently, the block of the Fisher information matrix corresponding to $\Sigma_1$ is the same in both $\mathcal{P}_1$ and $\mathcal{P}$. An analogous statement holds for the $\mathbb{R} \times \mathbb{P}(\mathcal{S}^{p_2}_{++})$ component of $\Skron$. 
\end{customcor}
\begin{proof}
The Wishart model $\mathcal{P}_1$ has the Riemannian metric on $T_{\Sigma_1} \mathcal{S}^{p_1}_{++}$
\begin{align*}
    \langle A_1,B_1 \rangle_{\Sigma_1} = p_2\tr(\Sigma_1^{-1}A_1\Sigma_1^{-1}B_1)
\end{align*}
The model $\mathcal{P}$ has the Fisher information metric
\begin{align}
\nonumber
     \langle (A_1 \otimes \Tilde{\Sigma}_2) + & (\Sigma_1 \otimes A_2), (B_1 \otimes \Tilde{\Sigma}_2) + (\Sigma_1 \otimes B_2) \rangle_{\Sigma_1 \otimes \Tilde{\Sigma}_2} 
    \\
    \nonumber
    & =  \tr( (\Sigma_1^{-1}A_1\Sigma_1^{-1}B_1 \otimes I_{p_2}) + \tr(I_{p_1} \otimes \Tilde{\Sigma}_2^{-1}A_2\Tilde{\Sigma}_2^{-1}B_2) 
    \\
    \label{eqn:FIMEquality1}
    & = p_2\tr((\Sigma_1^{-1}A_1\Sigma_1^{-1}B_1) + p_1\tr(\Tilde{\Sigma}_2^{-1}A_2\Tilde{\Sigma}_2^{-1}B_2) 
\end{align}
on $\mathcal{S}^{p_1}_{++}\times \mathbb{P}(\mathcal{S}^{p_2}_{++})$, where $A_2,B_2 \in T_{\Tilde{\Sigma}_2}\mathbb{P}(\mathcal{S}^{p_2}_{++})$ have the property that 
\begin{align*}
    \tr(\Tilde{\Sigma}_2^{-1/2}A_2 \Tilde{\Sigma}_2^{-1/2}) = \tr(\Tilde{\Sigma}_2^{-1/2}B_2 \Tilde{\Sigma}_2^{-1/2}) = 0.
\end{align*} 
The first term in the expression \eqref{eqn:FIMEquality1} involving $A_1,B_1 \in T_{\Sigma_1}\mathcal{S}^{p_1}_{++}$ is the same as the FIM in the model $\mathcal{P}_1$, proving that the restricted FIM on the first two components of $\mathcal{P}$ is identical to the FIM on $\mathcal{P}_1$. The orthogonality of the $\varphi_1$ parameterization of $\mathcal{P}_1$ follows from a similar computation to that in the proof of Theorem \ref{thm:OrthogFIMDecomposition}. 
\end{proof}

\begin{customcor}{2}
As $n \rightarrow \infty$, the components of $\sqrt{n}\big(\varphi(R(S^{(n)})) - \varphi(\Sigma)\big)$, namely
\begin{align*}
    \sqrt{n}\bigg(\tfrac{1}{p}\log( \vert S^{(n)} \vert) - c,\frac{\tr_1^{\det}(S^{(n)})}{\vert \tr_1^{\det}(S^{(n)}) \vert^{1/p_1} } - \Tilde\Sigma_1, \frac{\tr_2^{\det}(S^{(n)})}{\vert \tr_2^{\det}(S^{(n)}) \vert^{1/p_2}} - \Tilde\Sigma_2 \bigg),
\end{align*}
are asymptotically independent.
Similarly, the components of $\sqrt{n}\big( \varphi(M(S^{(n)})) - \varphi(\Sigma)\big)$ are asymptotically independent.
\end{customcor}
\begin{proof}
    The delta method shows that $\sqrt{n}\big(\varphi(R(S^{(n)})) - \varphi(\Sigma)\big)$ is asymptotically normal with mean zero. It remains compute the asymptotic covariance of this random vector. Define $\langle \cdot,\cdot\rangle_{\varphi(\Sigma)}^*$ to be the induced Riemannian metric on $\mathbb{R} \times \mathbb{P}(\mathcal{S}^{p_1}_{++}) \times \mathbb{P}(\mathcal{S}^{p_2}_{++})$. By the reparameterization Lemma \ref{lem:ReparamVariance}  
    \begin{align*}
      \text{AVar}_n\left( \langle D^\varphi_{\varphi(\Sigma)}(H), \sqrt{n}\big(\varphi(R(S^{(n)})) - \varphi(\Sigma)\big) \rangle_{\varphi(\Sigma)}^* \right) =  \text{AVar}_n( \langle H, \sqrt{n}(R(S^{(n)}) - \Sigma)\rangle_\Sigma),
    \end{align*}
    for $H \in T_{\Sigma}\Skron$.
    As $R(S^{(n)})$ is asymptotically efficient this variance equals $\Vert H\Vert^2_{\Sigma} = (\Vert D^\varphi_{\varphi(\Sigma)}(H)\Vert^*_{\varphi
(\Sigma)})^2$ by Lemma \ref{lem:AsympVarianceCurvedFamilies}. By polarization, the asymptotic variance matrix determines the asymptotic covariance matrix, where we conclude that 
\begin{align*}
    \text{ACov}_n\left( \langle D^\varphi_{\varphi(\Sigma)}(H_1), \sqrt{n}\big(\varphi(R(S^{(n)})) - \varphi(\Sigma)\big) \rangle_{\varphi(\Sigma)}^*, \langle D^\varphi_{\varphi(\Sigma)}(H_2), \sqrt{n}\big(\varphi(R(S^{(n)})) - \varphi(\Sigma)\big) \rangle_{\varphi(\Sigma)}^*\right) 
    \\
     = \langle D^\varphi_{\varphi(\Sigma)}(H_1), D^\varphi_{\varphi(\Sigma)}(H_2) \rangle_{\varphi(\Sigma)}.
\end{align*}
The orthogonality of the $\varphi$-parameterization shown in Theorem \ref{thm:OrthogFIMDecomposition} implies for instance that if $D^{\varphi}_{\varphi(\Sigma)}(H_1) = (v,0,0)$ and $D^{\varphi}_{\varphi(\Sigma)}(H_1) = (0,A,0)$ then $\langle D^\varphi_{\varphi(\Sigma)}(H_1), D^\varphi_{\varphi(\Sigma)}(H_2) \rangle_{\varphi(\Sigma)} = 0$. This implies that the components of $R(\varphi(S^{(n)}))$, namely $\sqrt{n}(\tfrac{1}{p}\log(\vert S^{(n)}\vert) - c)$, and $\sqrt{n}\big( \frac{\tr_i^{\det}(S^{(n)})}{\vert \tr_i^{\det}(S^{(n)}) \vert^{1/p_i} } - \Tilde\Sigma_i\big)$ $i = 1,2$, are asymptotically pairwise-uncorrelated. As these components are jointly asymptotically normal, the desired independence statement follows. 
\end{proof}

\begin{customlem}{4}
\label{lem:TestingIndepAppendix}
Assume that $\psi_i(S_i^{(n)})$ is a smooth, scale-invariant test statistic for the hypothesis $H_i:\Sigma_i \in \mathcal{H}_i$ in the model $np_j S_i^{(n)} \sim \text{Wishart}_{p_i}(\Sigma_i,np_j)$, $i,j = 1,2$, $i \neq j$, with $S_1^{(n)}$ and $S_2^{(n)}$ independent. Further assume that $c_i^{(n)}(\psi_i(S^{(n)}) - \psi(\Sigma_i))$ converges weakly to a non-degenerate distribution that does not depend on $\Sigma_i$ under $\mathcal{H}_i$. Let $H_1 \cap H_2: \{\Sigma_1 \otimes \Sigma_2: \Sigma_i \in \mathcal{H}_i, i = 1,2\}$ be the intersection hypothesis with respect to the model $nS^{(n)} \sim \text{Wishart}_{p}(\Sigma_1 \otimes \Sigma_2,n)$. The asymptotic distribution of the test statistics 
\begin{align}
\left(  c_1^{(n)}(\psi_1(S_1^{(n)}) - \psi_1(\Sigma_1)),c_2^{(n)}(\psi_2(S_2^{(n)}) - \psi_2(\Sigma_2))\right)
\end{align}
equal the asymptotic distribution of
\begin{align}
\left(  c_1^{(n)}(\psi_1(\tr_1^{\det}(S^{(n)})) - \psi_1(\Sigma_1)),c_2^{(n)}(\psi_2(\tr_2^{\det}(S^{(n)})) - \psi_2(\Sigma_2))\right).
\end{align}
In particular, if $q_{i\alpha}$ is an asymptotic level-$\alpha$ quantile of $c_i^{(n)}(\psi_i(S_i^{(n)}) - \psi_i(\Sigma_i))$, the test that rejects $H_1 \cap H_2$ when $c_i^{(n)}(\psi_i(\tr_i^{\det}(S^{(n)})) - \psi_i(\Sigma_i)) \geq q_{i\alpha}$ for either $i = 1$ or $i = 2$ has asymptotic level $1 - (1-\alpha)^2 = \alpha(2-\alpha)$. Similar statements apply, by replacing all occurrences of $\tr_i^{\det}(S^{(n)})$ above with the Kronecker MLE factors $\hat{\Sigma}_i$. 
\end{customlem}
\begin{proof}
Implicitly we are assuming that the asymptotic distribution of \eqref{eqn:TestStatIndepWish} is obtained by performing a Taylor expansion on the $\psi_i$s and applying the delta method. Consequently, using the scale invariance of the test statistics, it is enough to show that 
\begin{align}
\label{eqn:Indepwishstats}
    \left( \sqrt{n}\left( \frac{S_1^{(n)}}{\vert S_1^{(n)}\vert^{1/p_1}} - \Tilde{\Sigma}_1\right),\sqrt{n}\left( \frac{S_2^{(n)}}{\vert S_2^{(n)}\vert^{1/p_2}} - \Tilde{\Sigma}_2\right)  \right)
\end{align}
has the same asymptotic normal distribution as
\begin{align}
\label{eqn:trdetstats}
        \left( \sqrt{n}\left( \frac{\tr_1^{\det}(S^{(n)})}{\vert \tr_1^{\det}(S^{(n)})\vert^{1/p_1}} - \Tilde{\Sigma}_1\right),\sqrt{n}\left( \frac{\tr_2^{\det}(S^{(n)})}{\vert \tr_2^{\det}(S^{(n)})\vert^{1/p_2}} - \Tilde{\Sigma}_2\right)  \right)
\end{align}

Both of the components in \eqref{eqn:Indepwishstats} and \eqref{eqn:trdetstats} are asymptotically independent, the former by construction, and the latter by Corollary \ref{Cor:AsymptIndependenceofRPT}. It remains to establish that $S_i^{(n)} \vert S_i^{(n)}\vert^{-1/p_i}$ has the same asymptotic variance as $tr_i^{\det}(S^{(n)}) \vert tr_i^{\det}(S^{(n)})  \vert^{-1/p_i}$. The proof of Corollary \ref{Cor:AsymptIndependenceofRPT} shows that
\begin{align*}
    \text{AVar}_n\left( p_j \left\langle A, \sqrt{n}\left(\frac{\tr_i^{\det}(S^{(n)})}{\vert \tr_i^{\det}(S^{(n)}) \vert^{1/p_i} } - \Tilde{\Sigma}_i \right) \right\rangle_{\Tilde{\Sigma}_i} \right) = p_j \Vert A \Vert_{\Tilde{\Sigma}_i}^2, \;\; j \neq i
\end{align*}
where $\Vert \cdot \Vert_{\Tilde{\Sigma}_i}$ is the FIM \eqref{eqn:AIMetricMean} and $A$ is in $T_{\Tilde{\Sigma}_i}\mathbb{P}(\mathcal{S}^{p_i}_{++})$. 

The same reasoning as in Theorem \ref{thm:OrthogFIMDecomposition} can be applied to show that the map $\Sigma \mapsto (\tfrac{1}{p}\log(\vert \Sigma_i \vert),\Sigma_i \vert \Sigma_i \vert^{-/p_i})$ gives an orthogonal parameterization of the family $\{\text{Wishart}_{p_i}(\Sigma_i,1): \Sigma_i \in \mathcal{S}^{p_i}_{++}\}$ onto the space $\mathbb{R} \times \mathbb{P}(\mathcal{S}^{p_i}_{++})$, where the factor $\mathbb{P}(\mathcal{S}^{p_i}_{++
})$ has FIM \eqref{eqn:AIMetricMean} restricted to $T_{\Sigma_i \vert \Sigma_i \vert^{-/p_i}}\mathbb{P}(\mathcal{S}^{p_i}_{++})$. As $S_i^{(n)}$ is an efficient estimator of $\Sigma_i$, the reparameterization Lemma \ref{lem:ReparamVariance} implies that 
\begin{align*}
     \text{AVar}_n\left( \left\langle A, \sqrt{np_j}\left(\frac{S_i^{(n)}}{\vert S_i^{(n)}\vert^{1/p_i}} - \Tilde{\Sigma}_i \right) \right\rangle_{\Tilde{\Sigma}_i} \right) = \Vert A \Vert_{\Tilde{\Sigma}_i}^2,
\end{align*}
from which we obtain the desired equality 
\begin{align*}
        \text{AVar}_n\left( \left\langle A, \sqrt{n}\left(\frac{\tr_i^{\det}(S^{(n)})}{\vert \tr_i^{\det}(S^{(n)}) \vert^{1/p_i} } - \Tilde{\Sigma}_i \right) \right\rangle_{\Tilde{\Sigma}_i} \right) & = \tfrac{1}{p_j} \Vert A \Vert_{\Tilde{\Sigma}_i}^2
        \\
       & = \text{AVar}_n\left( \left\langle A, \sqrt{n}\left(\frac{S_i^{(n)}}{\vert S_i^{(n)}\vert^{1/p_i}} - \Tilde{\Sigma}_i \right) \right\rangle_{\Tilde{\Sigma}_i} \right).
\end{align*}
The statement regarding the asymptotic level of the joint test of $H_1 \cap H_2$ is a consequence of the asymptotic independence of the components of \eqref{eqn:trdetstats}. As the maximum likelihood estimator is asymptotically efficient, the asymptotic normal distribution of \eqref{eqn:trdetstats} will be the same as the same statistic that has the occurrences of $\tr_i^{\det}(S^{(n)})$ replaced by $\hat{\Sigma}_i$.  
\end{proof}

\begin{customlem}{0.7}
\label{lem:PTComponentConvergence}
     Let $n_mS^{(m)} \sim \text{Wishart}_{p_m}(\otimes_{i = 1}^k \Sigma_i^{(m)},n_m)$ with $\Sigma_i^{(m)} \in \mathcal{S}^{p_{im}}_{++}$ and $p_m = \prod_{i = 1}^k p_{im}$. Denote the length $p_{im}$ vector of eigenvalues of $\Sigma_i^{(m)}$ by $\lambda_i^{(m)}$. Then 
     \begin{align}
     E\bigg(    \bigg\Vert \frac{\tr_j(S^{(m)})}{\tr( \otimes_{i = 1}^k \Sigma_i^{(m)})} - \frac{\Sigma_j^{(m)}}{\tr(\Sigma_j^{(m)})} \bigg\Vert_F^2\bigg) \leq    \frac{2}{n_m \prod_{i \neq j}^k\big(p_{im} \cos^2(\angle \lambda_i^{(m)}, 1_{p_{im}})  \big)},
     \end{align}
     and 
     \begin{align*}
        \frac{1}{\Vert \Sigma_j \Vert_F} \bigg\Vert \frac{\tr_j(S^{(m)})}{\tr( \otimes_{i \neq j} \Sigma_i^{(m)})} - \Sigma_j^{(m)} \bigg\Vert_F = O_P\bigg(    \frac{ \sqrt{p_{jm}} \cos\big(\angle \lambda_j^{(m)},1_{p_{jm}}\big) }{\sqrt{n_m} \prod_{i \neq j}^k\big(\sqrt{p_{im}} \cos(\angle \lambda_i^{(m)}, 1_{p_{im}})  \big)}\bigg). 
     \end{align*}
\end{customlem}
\begin{proof}
    Without loss of generality we take $j = 1$ and by the orthogonal equivariance of the partial trace operators we assume that $\Sigma_i^{(m)} = \text{diag}(\lambda_i^{(m)})$ is diagonal for $i = 1,\ldots,k$.

 We can write $S^{(m)} \overset{d}{=} (\otimes_{i = 1}^k \Sigma_i^{(m)})^{1/2}W^{(m)} (\otimes_{i = 1}^k \Sigma_i^{(m)})^{1/2}$ with $n_m W^{(m)} \sim \text{Wishart}_{p_m}(I_{p_m},n_m)$. Using the equivariance and cyclic permutation properties of the partial trace,   
 \begin{align}
 \label{eqn:PTasymp1}
       \bigg\Vert \frac{\tr_1(S^{(m)})}{\tr( \otimes_{i = 1}^k \Sigma_i^{(m)})} - \frac{\Sigma_1^{(m)}}{\tr(\Sigma_1^{(m)})} \bigg\Vert_F^2 & \overset{d}{=}  \frac{1}{\tr(\Sigma_1^{(m)})^2} \bigg\Vert (\Sigma_1^{(m)})^{1/2} \bigg( \frac{\tr_1((I_{p_{1m}} \otimes_{i > 1} \Sigma_i^{(m)}) W^{(m)}) }{ \prod_{i > 1}\tr(\Sigma_i^{(m)})} - I_{p_{1m}}\bigg) (\Sigma_1^{(m)})^{1/2} \bigg\Vert_F^2.
 \end{align}
To ease notation define
 \begin{align*}
     M \coloneqq  \bigg( \frac{\tr_1((I_{p_{1m}} \otimes_{i > 1} \Sigma_i^{(m)}) W^{(m)}) }{ \prod_{i > 1}\tr(\Sigma_i^{(m)})} - I_{p_{1m}}\bigg), 
 \end{align*}
 and let $\Sigma_1^{(m)} = \text{diag}(\alpha_1,\ldots,\alpha_{p_{1m}}) = \text{diag}(\alpha)$.
 The right hand side of \eqref{eqn:PTasymp1}
 becomes
 \begin{align}
\nonumber   \frac{1}{\tr(\Sigma_1^{(m)})^2} \bigg\Vert (\Sigma_1^{(m)})^{1/2} M (\Sigma_1^{(m)})^{1/2} \bigg\Vert_F^2 & = \frac{1}{(1_{p_{1m}}^\top \alpha)^2} \sum_{i,j} \big(\alpha_i^{1/2} \alpha_j^{1/2} M_{ij}\big)^2
   \\
   \label{eqn:PTasympeqn2}
   & = \frac{1}{(1_{p_{1m}}^\top \alpha)^2}  \sum_{i,j} \alpha_i \alpha_j M_{ij}^2
 \end{align}
We denote the unit-trace diagonal matrix $\otimes_{i > 1} \big(\tr(\Sigma_i^{(m)})^{-1}\Sigma_i^{(m)}\big)$  by $A$ and decompose $M$ as 
 \begin{align}
 \nonumber   \frac{\tr_1((I_{p_{1m}} \otimes_{i > 1} \Sigma_i^{(m)}) W^{(m)}) }{ \prod_{i > 1}\tr(\Sigma_i^{(m)})} - I_{p_{1m}} & =\tr_1((I_{p_{1m}} \otimes A)W^{(m)}) - I_{p_{1m}} 
     \\
     \label{eqn:NormofPTDif}
     & =  \sum_{i_2,\ldots,i_k} a_{i_2\ldots i_k} \big(W_{\cdot i_2\ldots i_k,\cdot i_2\ldots i_k}^{(m)} - I_{p_{1m}}\big),
 \end{align}
 where $a_{i_2\ldots i_k}$ are the diagonal entries of $A$ and $W_{\cdot i_2\ldots i_k,\cdot i_2\ldots i_k}^{(m)}$ are the $p_{1m} \times p_{1m}$ dimensional blocks of $W^{(m)}$. All of the scaled blocks $n_m W_{\cdot i_2\ldots i_k,\cdot i_2\ldots i_k}^{(m)}$ are independent and follow the distribution $\text{Wishart}_{p_m}(I_{p_{1m}},n_m)$. The expectation of $M_{ij}^2$ equals
 \begin{align*}
    E(M_{ij}^2) & =  E\bigg(\big(\sum_{i_2,\ldots,i_k} a_{i_2\ldots i_k} (W_{i i_2\ldots i_k,j i_2\ldots i_k}^{(m)} - \delta_{ij})\big)^2\bigg)
    \\
    & =  \sum_{i_2,\ldots,i_k} a_{i_2\ldots i_k}^2 E((W_{i i_2\ldots i_k,j i_2\ldots i_k}^{(m)} - \delta_{ij})^2)
    \\
    & =  \sum_{i_2,\ldots,i_k} a_{i_2\ldots i_k}^2 \text{Var}\big(W_{i i_2\ldots i_k,j i_2\ldots i_k}^{(m)}\big)
    \\
    & \leq 2 n_m^{-1}  \sum_{i_2,\ldots,i_k} a_{i_2\ldots i_k}^2
    \\
    & = 2n_m^{-1} \Vert A\Vert_F^2
\end{align*}
The Frobenius norm of $A$ can be computed as
\begin{align*}
    \Vert A \Vert_F^2 = \prod_{i > 1} \frac{\Vert \Sigma_i^{(m)} \Vert_F^2}{\tr(\Sigma_i^{(m)})^2} =  \prod_{i > 1} \frac{ \Vert \lambda_i^{(m)} \Vert^2 }{\langle 1_{p_{im}},\lambda_i^{(m)} \rangle^2} =  \prod_{i > 1} \frac{1}{p_{im} \cos^2\big( \angle \lambda_i^{(m)},1_{p_{im}}\big)} 
\end{align*}
Plugging this back into    \eqref{eqn:PTasympeqn2} gives 
\begin{align*}
      E\bigg(    \bigg\Vert \frac{\tr_1(S^{(m)})}{\tr( \otimes_{i = 1}^k \Sigma_i^{(m)})} - \frac{\Sigma_1^{(m)}}{\tr(\Sigma_1^{(m)})} \bigg\Vert_F^2\bigg) & \leq \frac{2}{(1_{p_{1m}}^\top \alpha)^2 n_m\prod_{i > 1}^kp_{im} \cos^2\big( \angle \lambda_i^{(m)},1_{p_{im}}\big)} \big(\sum_{i = 1}^{p_{1m}} \alpha_i \big)^2
      \\
      & = \frac{2}{ n_m\prod_{i > 1}^kp_{im} \cos^2\big( \angle \lambda_i^{(m)},1_{p_{im}}\big)}.
\end{align*}
The remaining result follows from
\begin{align*}
     \bigg\Vert \frac{\tr_1(S^{(m)})}{\tr( \otimes_{i = 1}^k \Sigma_i^{(m)})} - \frac{\Sigma_1^{(m)}}{\tr(\Sigma_1^{(m)})} \bigg\Vert_F^2 =  \bigg(\frac{\Vert \Sigma_1^{(m)}\Vert_F}{\tr(\Sigma_1^{(m)})}\bigg)^2  \bigg(\frac{1}{\Vert \Sigma_1 \Vert_F}\bigg\Vert \frac{\tr_1(S^{(m)})}{\tr( \otimes_{i > 1} \Sigma_i^{(m)})} - \Sigma_1^{(m)}\bigg\Vert_F\bigg)^2,
\end{align*}
where
\begin{align*}
    \bigg(\frac{\Vert \Sigma_1^{(m)}\Vert_F}{\tr(\Sigma_1^{(m)})}\bigg)^2  = \frac{1}{ p_{1m} \cos^2(\angle \lambda_1^{(m)},1_{p_{1m}})}.
\end{align*}
 \end{proof}

\begin{customthm}{4}
    Let $n_mS^{(m)} \sim \text{Wishart}_{p_m}(\otimes_{i = 1}^k \Sigma_i^{(m)},n_m)$ with $\Sigma_i^{(m)} \in \mathcal{S}^{p_{im}}_{++}$ and $p_m = \prod_{i = 1}^k p_{im}$. Denote $\otimes_{i = 1}^k \Sigma_i^{(m)}$ by $\Sigma^{(m)}$ and take $\lambda^{(m)}$ to be the length $p_{m}$ vector of eigenvalues of $\Sigma^{(m)}$. If
    \begin{align}
        \frac{ \sup_{j \in \{1,\ldots,k\}}p_{jm}\cos^2\big(\lambda_j^{(m)},1_{p_{jm}}\big)}{\sqrt{n_m} \sqrt{p_m} \cos(\angle \lambda^{(m)}, 1_{p_{m}})}  = O(1), \; \frac{1}{\sqrt{n_m} \sqrt{p_m} \cos(\angle \lambda^{(m)},1_{p_{m}})} = o(1)
    \end{align}
    then 
\begin{align}
    \frac{\big\Vert P(S^{(m)}) - \Sigma^{(m)} \big\Vert_F}{\big\Vert \Sigma^{(m)} \big\Vert_F}  = O_P\bigg(       \frac{ \sup_{j \in \{1,\ldots,k\}}p_{jm} \cos^2\big(\angle \lambda_j^{(m)},1_{p_{jm}}\big) }{\sqrt{n_m} \sqrt{p_m} \cos(\angle \lambda^{(m)}, 1_{p_{m}})}\bigg)
\end{align}
\end{customthm} 
\begin{proof}
    First we compute the variance of $\tr(S^{(m)})/\tr(\Sigma^{(m)}) $. Without loss of generality $\Sigma^{(m)}$ can be taken to be diagonal, and after rescaling by $\tr(\Sigma^{(m)})$, it can be assumed that $n_mS^{(m)} \sim \text{Wishart}_{p_m}(\text{diag}(\alpha^{(m)}),n_m)$, with $\alpha^{(m)}$ being equal to the diagonal entries of $\Sigma^{(m)}$ divided by their sum. Taking $n_mW^{(m)} \sim \text{Wishart}_{p_m}(I_{p_m},n_m)$, $S^{(m)}$ has the same distribution as $\text{diag}(\alpha^{(m)})^{1/2}W^{(m)}\text{diag}(\alpha^{(m)})^{1/2}$. Then 
    \begin{align}
     \nonumber   \text{Var}(\tr(S^{(m)})) & = \text{Var}\big(\sum_{i = 1}^{p_m} \alpha^{(m)}_i w_i^2\big)
        \\
  \nonumber      & = \frac{2}{n_m}\sum_{i = 1}^{p_m} (\alpha_i^{(m)})^2
        \\
   \nonumber     & = \frac{2}{n_m} \frac{ \langle \otimes_{i = 1}^k \lambda_i^{(m)}, \otimes_{i = 1}^k \lambda_i^{(m)} \rangle }{ \langle 1_{p_m}, \otimes_{i = 1}^k \lambda_i^{(m)}\rangle^2}
        \\
\label{eqn:PTConvergenceProof1}
        & = \frac{2}{n_m p_m} \prod_{i = 1}^k \cos^{-2}(\angle \lambda_i^{(m)},1_{p_{im}}), 
    \end{align}
    where we have used the fact that the diagonal entries $w_i^2$ of $W^{(m)}$ have independent scaled chi-squared distributions with variance $2n_m^{-1}$. By Chebychev's inequality,
    \begin{align*}
        \frac{\tr(S^{(m)})}{\tr(\Sigma^{(m)})} - 1 = O_P\bigg( \frac{1}{\sqrt{n_m} \sqrt{p_m} \prod_{i = 1}^k \cos(\angle \lambda_i^{(m)},1_{p_{im}})}\bigg).
    \end{align*}
    Note $\lambda^{(m)} = \otimes_{i = 1}^k \lambda_i^{(m)}$ so that 
    \begin{align*}
        \prod_{i = 1}^k \cos(\angle \lambda_i^{(m)},1_{p_{im}}) = \cos( \angle \lambda^{(m)}, 1_{p_m}).
    \end{align*}

Define $z^{(m)} \coloneqq \tr(S^{(m)})/\tr(\Sigma^{(m)})$. We can decompose the norm of interest as
\begin{align*}
    \frac{ \Vert P(S^{(m)}) - \Sigma^{(m)} \big\Vert_F}{\big\Vert \Sigma^{(m)} \big\Vert_F} = \bigg\Vert (z^{(m)})^{-(k-1)} \bigotimes_{i = 1}^k \frac{\tr_i(S^{(m)})}{\Vert \Sigma_i^{(m)} \Vert_F \tr( \otimes_{j \neq i} \Sigma_j^{(m)})} - \bigotimes_{i = 1}^k \frac{\Sigma_i^{(m)}}{\Vert \Sigma_i^{(m)}\Vert_F} \bigg\Vert_F,
\end{align*}
which is bounded above by 
\begin{align*}
&\bigg\Vert \bigotimes_{i = 1}^k \frac{\tr_i(S^{(m)})}{\Vert \Sigma_i^{(m)} \Vert_F \tr( \otimes_{j \neq i} \Sigma_j^{(m)})} - \bigotimes_{i = 1}^k \frac{\Sigma_i^{(m)}}{\Vert \Sigma_i^{(m)}\Vert_F} \bigg\Vert_F + \bigg(\frac{1}{(z^{(m)})^{k-1}} - 1 \bigg)  \prod_{i = 1}^k \frac{\Vert \tr_i(S^{(m)}) \Vert_F}{\Vert \Sigma_i^{(m)} \Vert_F \tr( \otimes_{j \neq i} \Sigma_j^{(m)})}
\\
&\coloneqq (A) + (B)
\end{align*}


We start by bounding the first term $(A)$.
By Lemma \ref{lem:PTComponentConvergence} 
\begin{align*}
        \bigg\Vert \frac{\tr_1(S^{(m)})}{\Vert \Sigma_1 \Vert_F\tr( \otimes_{j \neq 1} \Sigma_j^{(m)})} - \frac{\Sigma_1^{(m)}}{\Vert \Sigma_1^{(m)}\Vert_F} \bigg\Vert_F = O_P\bigg(    \frac{ \sqrt{p_{1m}} \cos\big(\angle \lambda_1^{(m)},1_{p_{1m}}\big) }{\sqrt{n_m} \prod_{j \neq 1}^k\big(\sqrt{p_{jm}} \cos(\angle \lambda_j^{(m)}, 1_{p_{jm}})  \big)}\bigg)
\end{align*}
and 
\begin{align*}
    \frac{\Vert\tr_i(S^{(m)})\Vert_F}{\Vert \Sigma_i^{(m)} \Vert_F \tr( \otimes_{j \neq i} \Sigma_j^{(m)})} & \leq 1 +      \bigg\Vert \frac{\tr_i(S^{(m)})}{\Vert \Sigma_1 \Vert_F\tr( \otimes_{j \neq i} \Sigma_j^{(m)})} - \frac{\Sigma_i^{(m)}}{\Vert \Sigma_i^{(m)}\Vert_F} \bigg\Vert_F
    \\
    & = 1 + O_P\bigg(    \frac{ \sqrt{p_{im}} \cos\big(\angle \lambda_i^{(m)},1_{p_{im}}\big) }{\sqrt{n_m} \prod_{j \neq i}\big(\sqrt{p_{jm}} \cos(\angle \lambda_j^{(m)}, 1_{p_{jm}})  \big)}\bigg),
\end{align*}

where by the first assumption \eqref{eqn:PTConvThmAssumptions} this term is $O_P(1)$.

Decomposing $(A)$, we get
\begin{align*}
   (A)  \leq &  \bigg\Vert \bigg(\frac{\tr_1(S^{(m)})}{\Vert \Sigma_1 \Vert_F\tr( \otimes_{j \neq 1} \Sigma_j^{(m)})} - \frac{\Sigma_1^{(m)}}{\Vert \Sigma_1^{(m)}\Vert_F}\bigg) \otimes \bigg(\bigotimes_{i > 1}^k \frac{\tr_i(S^{(m)})}{\Vert \Sigma_i^{(m)} \Vert_F \tr( \otimes_{j \neq i} \Sigma_j^{(m)})} \bigg)\bigg\Vert_F  + 
    \\
    & \bigg\Vert \bigg(\frac{\Sigma_1^{(m)}}{\Vert \Sigma_1^{(m)}\Vert_F}\bigg) \otimes \bigg(\bigotimes_{i > 1}^k \frac{\tr_i(S^{(m)})}{\Vert \Sigma_i^{(m)} \Vert_F \tr( \otimes_{j \neq i} \Sigma_j^{(m)})} \bigg) -  \bigg( \bigotimes_{i = 1}^k \frac{\Sigma_i^{(m)}}{\Vert \Sigma_i^{(m)}\Vert_F}\bigg) \bigg\Vert_F
    \\
    & \coloneqq (C) + (D).
\end{align*}
Using the fact that $\Vert \otimes_{i = 1}^k H_i\Vert_F = \prod_{i = 1}^k \Vert H_i \Vert_F$
\begin{align*}
    (C) = O_P\bigg(    \frac{ \sqrt{p_{1m}} \cos\big(\angle \lambda_1^{(m)},1_{p_{1m}}\big) }{\sqrt{n_m} \prod_{j \neq 1}^k\big(\sqrt{p_{jm}} \cos(\angle \lambda_j^{(m)}, 1_{p_{jm}})  \big)}\bigg).
\end{align*}
Inductively decomposing $(D)$ in a similar manner, we obtain
\begin{align*}
    (A) = O_P\bigg(  \sup_{i \in \{1,\ldots,k\}}  \frac{ \sqrt{p_{im}} \cos\big(\angle \lambda_i^{(m)},1_{p_{im}}\big) }{\sqrt{n_m} \prod_{j \neq i}^k\big(\sqrt{p_{jm}} \cos(\angle \lambda_j^{(m)}, 1_{p_{jm}})  \big)}\bigg).
\end{align*}
Next we bound $(B)$ by
\begin{align*}
(B) & =  \bigg(\frac{1}{(z^{(m)})^{k-1}} - 1 \bigg)  \prod_{i = 1}^k \frac{\Vert \tr_i(S^{(m)}) \Vert_F}{\Vert \Sigma_i^{(m)} \Vert_F \tr( \otimes_{j \neq i} \Sigma_j^{(m)})}
\\
& = \bigg(\frac{1 - (z^{(m)})^{k-1} }{(z^{(m)})^{k-1}}\bigg) O_P(1)
    \\
 & =  \bigg(\frac{1 - (z^{(m)})^{k-1} }{1 + o_P(1)}\bigg) O_P(1).  
\end{align*}
The second assumption in \eqref{eqn:PTConvThmAssumptions} is used to get that $z^{(m)} \overset{P}{\rightarrow} 1$ in the denominator above.
By the binomial theorem
\begin{align*}
   (z^{(m)})^{k-1} - 1 = \sum_{a = 1}^{k-1}(z^{(m)} - 1)^{a} = O_P\bigg( \frac{1}{\sqrt{n_m} \sqrt{p_m} \prod_{i = 1}^k \cos(\angle \lambda_i^{(m)},1_{p_{im}})}\bigg),
\end{align*}
gives the correct convergence rate for $(B)$. 
\end{proof}

\begin{customlem}{5}
Let $C(x) = \sum_{q = 0}^{Q} c_q x^{q}$ be a non-zero polynomial with non-negative coefficients $c_q \geq 0$.
    If $\lambda^{(m)} = (C(1),C(2),\ldots,C(m))$, then $\cos^{-2}(\angle \lambda^{(m)},1_m) = O(1)$ as $m \rightarrow \infty$. If  $\lambda^{(m)} = (a^1,a^2,\ldots,a^m)$, $a > 1$ grows exponentially then $\cos^{-2}(\angle \lambda^{(m)},1_m) = O(m)$ as $m \rightarrow \infty$.  Lastly, if $\lambda^{(m)} = (a_m + b_m,\ldots,a_m + b_m,b_m,\ldots,b_m)$ has a spiked eigenvalue structure with $q_m > 0$ entries that equal $a_m + b_m$, $m - q_m$ entries that equal $b_m$, and $a_m,b_m \geq 0$, then
    \begin{align*}
        \cos^{-2}(\angle \lambda^{(m)},1_m) = \frac{r_ma_m^2 + 2r_m a_mb_m + b_m^2}{r_m^2 a_m^2 + 2r_m a_mb_m + b_m^2}
    \end{align*}
where $r_m = q_m/m$. This expression is $O(1)$ if $a_m/b_m = O(1)$ or $m/q_m = O(1)$.
\end{customlem}
\begin{proof}
    In the case $\lambda^{(m)} = (1^q,2^q,\ldots,m^q)$, $q \geq 0$  we compute
    \begin{align*}
        \cos(\angle \lambda^{(m)},1_m) & = \frac{\langle \lambda^{(m)},1_m \rangle}{ \Vert \lambda^{(m)} \Vert 1_m \Vert}
        \\
        & = \frac{\sum_{i = 1}^m i^q}{\big( \sum_{i = 1}^m i^{2q} \big)^{1/2} \sqrt{m}}
    \end{align*}
    By Faulhaber's formula $\sum_{i = 1}^m i^q \sim (q+1)^{-1}m^{q+1}$ and thus
\begin{align*}
    \frac{\sum_{i = 1}^m i^q}{\big( \sum_{i = 1}^m i^{2q} \big)^{1/2} \sqrt{m}} \sim \frac{(2q+1)^{1/2}m^{q+1}}{(q+1) m^{q + 1/2} \sqrt{m} } = O(1).
\end{align*}
The case for a general polynomial follows
similarly since
\begin{align*}
    \langle \lambda^{(m)},1_m\rangle = \sum_{i = 1}^m \sum_{q = 0}^Q c_q i^q =  \Theta\left(m^{Q+1}\right),
    \\
    \Vert \lambda^{(m)} \Vert^2 = \sum_{i = 1}^m \left( \sum_{q = 0}^Q c_q i^q \right)^2 = \Theta\left( m^{2Q+1} \right),
\end{align*}
assuming that $c_Q \neq 0$.

In the exponential growth case, $\lambda^{(m)} = (a^1,a^2,\ldots,a^m)$ and  
\begin{align*}
         \cos(\angle \lambda^{(m)},1_m) & = \frac{\sum_{i = 1}^m a^i}{\big( \sum_{i = 1}^m (a^2)^i \big)^{1/2} \sqrt{m}}
        \\
        & = \frac{(a^m - 1)(a^2 - 1)^{1/2}}{(a-1)(a^{2m} - 1)^{1/2} \sqrt{m}} \sim \frac{(a^2 - 1)^{1/2}}{(a-1) \sqrt{m}}.
\end{align*}
It follows that $\cos^{-2}(\angle \lambda^{(m)},1_m) \sim (a^2 - 1)(a-1)^{-2}m$, as needed.

Dropping the subscripts on $a_m,b_m$, for readability, the formula for the remaining case $\lambda^{(m)} = (a + b,\ldots,a + b,b,\ldots,b)$ follows from expression
\begin{align*}
     \cos^{-2}(\angle \lambda^{(m)},1_m) & = \bigg(\frac{q_m(a+b) + (m - q_m)b}{ \big(q_m (a + b)^2 + (m-q_m)b^2\big)^{1/2} \sqrt{m}}\bigg)^{-2}
     \\
     & = \frac{(q_m a^2 + 2q_m ab + mb^2)m}{q_m^2a^2 + 2mq_m  ab + m^2 b^2}
     \\
     & = \frac{\tfrac{q_m}{m} a^2 + 2\tfrac{q_m}{m} ab + b^2}{(\tfrac{q_m}{m})^2 a^2 + 2 \tfrac{q_m}{m}  ab +  b^2}
\end{align*}
If $b_m \neq 0$, $a_m/b_m \leq c$ for all $m$ then 
\begin{align*}
     \frac{r_m a_m^2 + 2r_m a_mb_m + b_m^2}{r_m^2 a_m^2 + 2 r_m  a_mb_m +  b_m^2} \leq r_m c^2 + 2r_m c + 1 = O(1).
\end{align*}
If $r_m^{-1} \leq d$ for all $m$ then 
\begin{align*}
      \frac{r_m a_m^2 + 2r_m a_mb_m + b_m^2}{r_m^2 a_m^2 + 2 r_m  a_mb_m +  b_m^2} \leq       \frac{d a_m^2 + 2d a_mb_m + r_m^{-2} b_m^2}{ a_m^2 + 2 r_m^{-1}  a_mb_m +  r_m^{-2} b_m^2} \leq d + 2d + 1 = O(1)
\end{align*}

\end{proof}

\section{Simulation Study Details}
\label{App:SimulationDetails}
\subsection{Hypothesis Testing}
The rejection probabilities presented in the simulation study in Section \ref{Sec:OrthogonalParam} are obtained from first drawing observations $nS \sim \text{Wishart}_{15}(\Sigma_1 \otimes \Sigma_2,n)$. Under the null hypothesis of sphericity and compound symmetry $\Sigma_1$ and $\Sigma_2$ are taken to be $\text{diag}(1,2,3,4,5)$ and $2I +  11^\top$ respectively. Under the simulations for the alternative hypothesis $\Sigma_1$ equals $\text{diag}(1,2,3,4,5) + 0.1 \times 11^\top$, and $\Sigma_2$ equals $2I +  11^\top - 0.3\times e_3e_3^\top$, where $e_3^\top = (0,0,1)$. 

The test statistic used to test $H_1$ is $d_{AI}(\hat{\Sigma}_1, f(\hat{\Sigma}_1))$ where $f(\hat{\Sigma}_1)$ equals the diagonal matrix that has the same diagonal entries as $\hat{\Sigma}_1$. The affine-invariant distance $d_{AI}$ is equal to the Riemannian distance on the set of positive definite matrices induced by the Fisher information metric. It is defined by $d_{AI}(A,B) = \Vert \log(A^{-1/2}BA^{-1/2})\Vert_F$, where $\log$ is the matrix logarithm. The sphericity test statistic is pivotal under the null hypothesis due to the affine invariance of $d_{AI}$:
\begin{align*}
    d_{AI}(CAC^\top,CBC^\top) = d_{AI}(A,B),\;\; A,B \in \Sp,\; C \in \Glp.
\end{align*}
In particular, this test statistic is scale-invariant. The $95\%$ quantile of $d_{AI}(\hat{\Sigma}_1, f(\hat{\Sigma}_1))$ is approximated by computing $d_{AI}(S_i,f(S_i))$ over $5000$ independent draws of $3n S_i \sim \text{Wishart}_5(I,3n)$. Invariance of the test statistic under the group of diagonal matrices, or equivalently the fact that the test statistic is pivotal, ensures that the quantile does not depend on which diagonal matrix is chosen for $E(S_i)$. Lemma \ref{lem:TestingIndep} implies that the
Wishart distribution with $p_2n = 3n$  degrees of freedom gives the correct quantile asymptotically. 

The test of compound symmetry that is used in this simulation study is the likelihood ratio test. The likelihood ratio test statistic is compared to the $95\%$ quantile of a chi-squared distribution with four degrees of freedom. The likelihood ratio statistic is scale-invariant, as needed for Lemma \ref{lem:TestingIndep} to apply. 

\subsection{Comtrade Data Plots}

\begin{figure}[h]
    \centering
    \includegraphics[width = 1\textwidth, height = 0.6\textwidth]{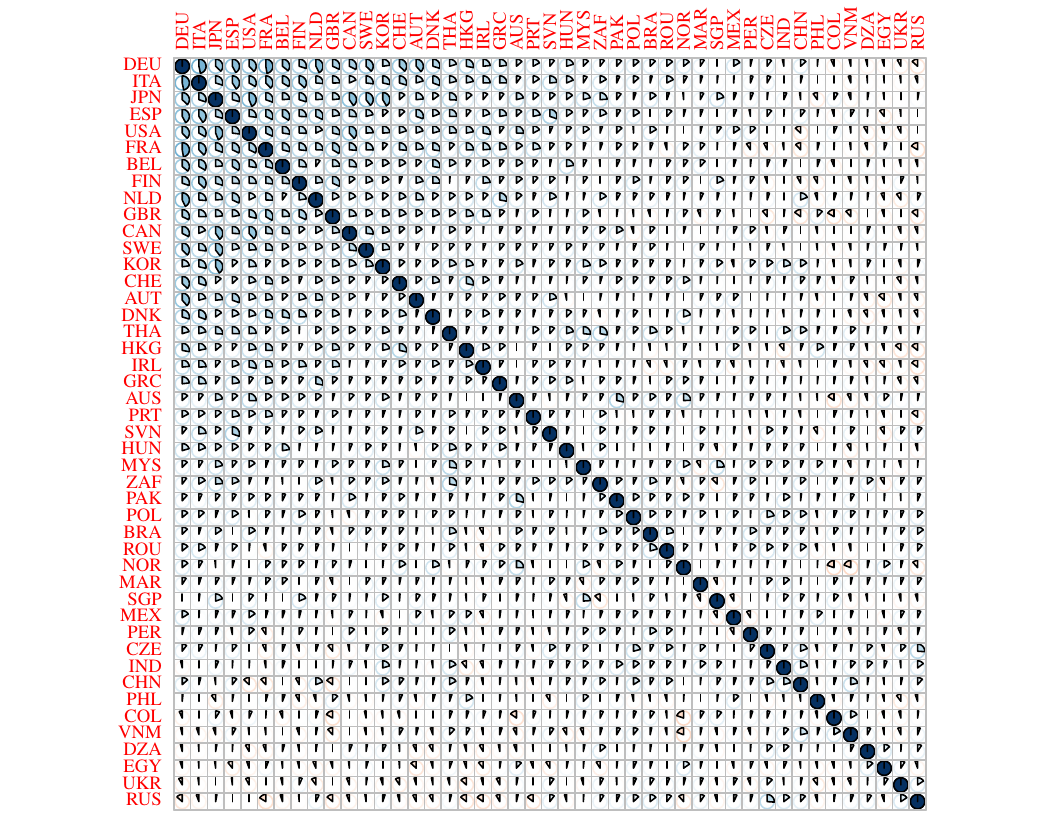}
    \caption{Partial trace estimate of the export correlation matrix  $\text{diag}(\Sigma_1)^{-1/2}\Sigma_1\text{diag}(\Sigma_1)^{-1/2}$.}
    \label{fig:Export}
\end{figure}

\begin{figure}[h]
    \centering
    \includegraphics[width = \textwidth, height = 0.6\textwidth]{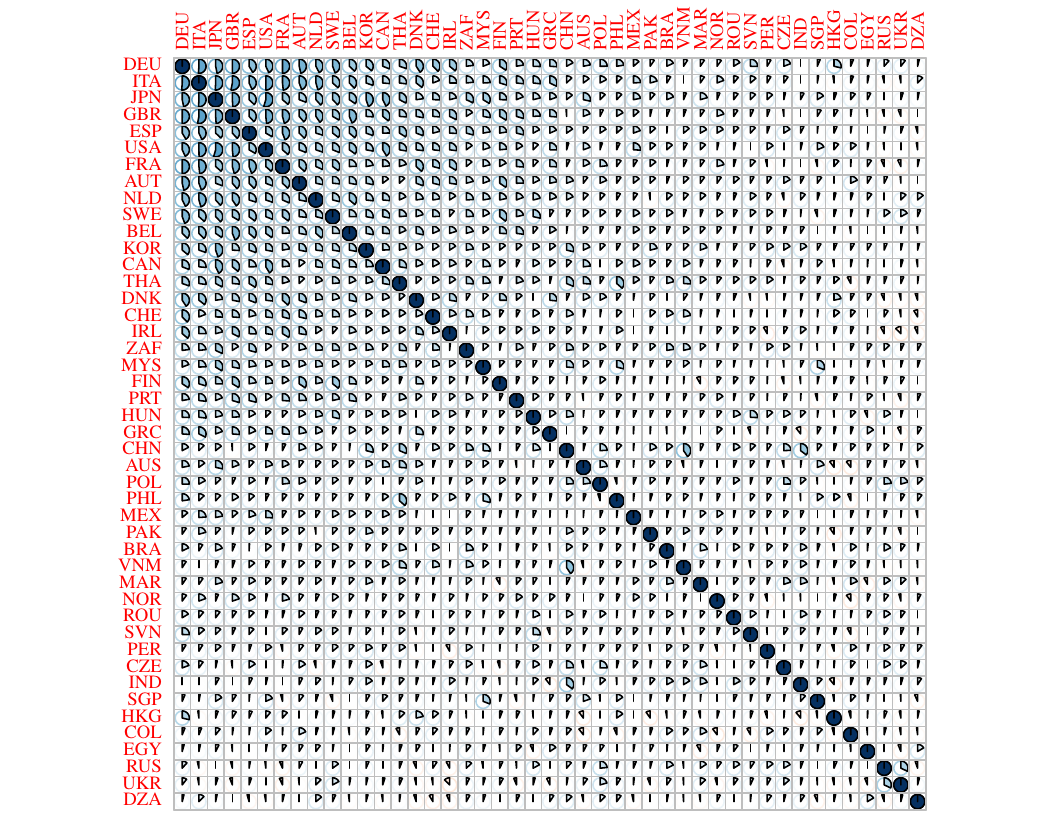}
    \caption{Partial trace estimate of the import correlation matrix  $\text{diag}(\Sigma_2)^{-1/2}\Sigma_2\text{diag}(\Sigma_2)^{-1/2}$.}
    \label{fig:Import}
\end{figure}

\begin{figure}[h]
    \centering
    \includegraphics[width = \textwidth, height = 0.6\textwidth]{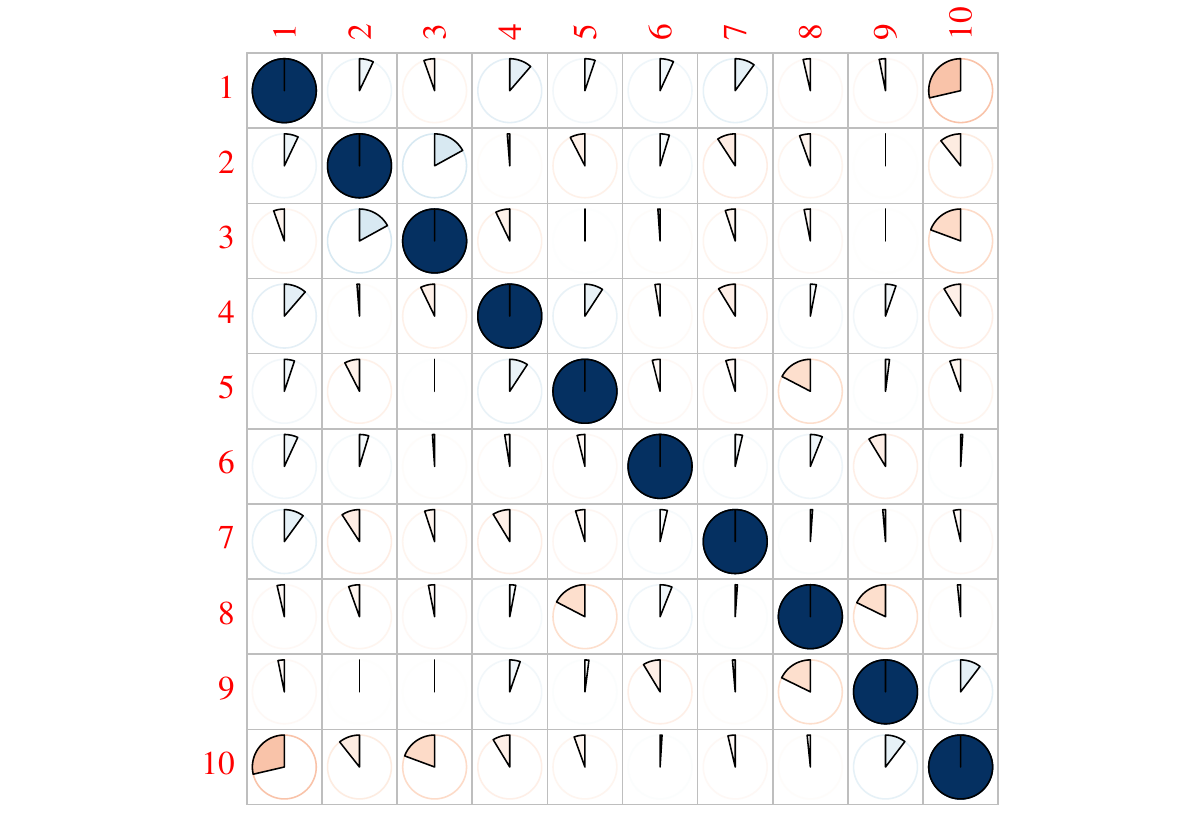}
    \caption{Partial trace estimate of the time correlation matrix $\text{diag}(\Sigma_3)^{-1/2}\Sigma_3\text{diag}(\Sigma_3)^{-1/2}$.}
    \label{fig:Time}
\end{figure}

\begin{figure}[h]
    \centering
    \includegraphics[width = .9\textwidth,height = 0.4\textwidth]{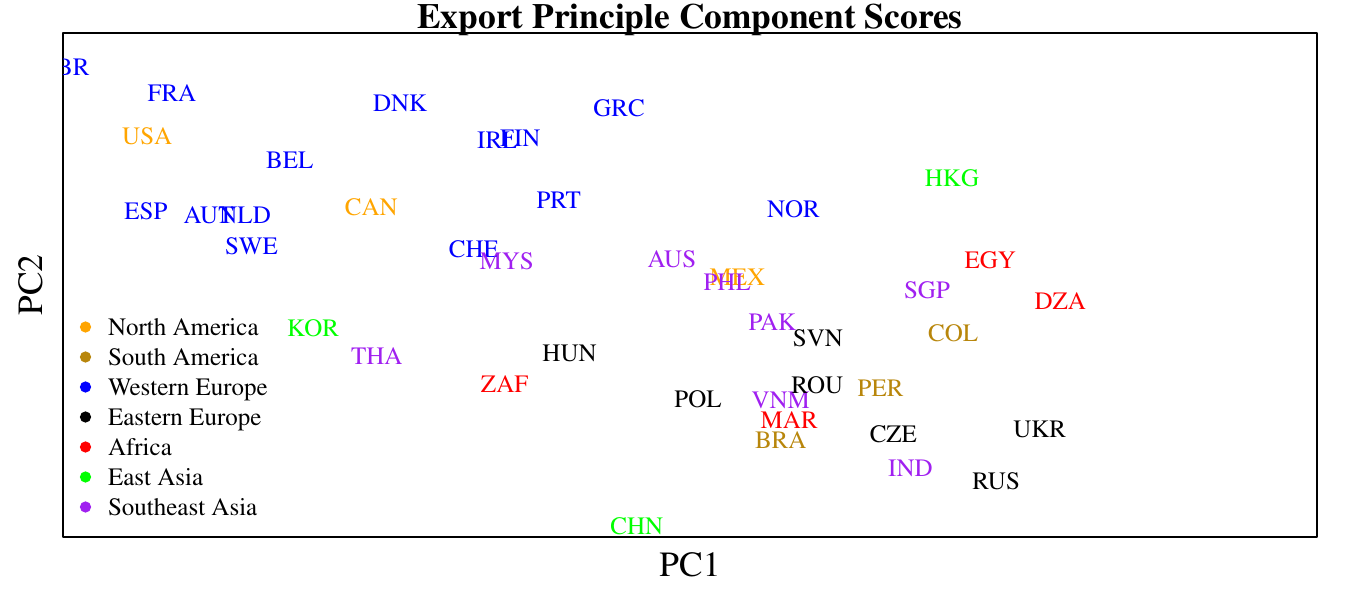}
    \caption{The first two principle component scores of the correlation matrix associated with $\hat{\Sigma}_1$.}
    \label{fig:PCExportScores}
\end{figure}

\begin{figure}[h]
    \centering
    \includegraphics[width = .8\textwidth,height = .55\textwidth]{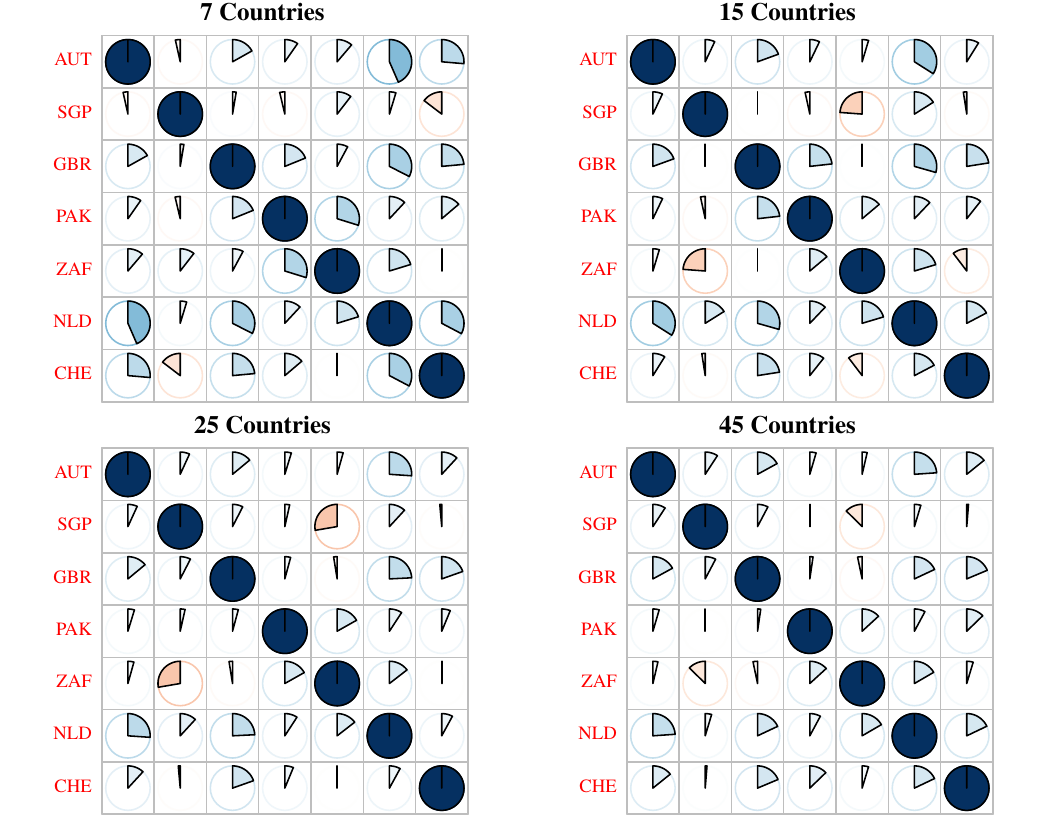}
    \caption{A second example of Figure \ref{fig:CorplotHoldOut} with $7$ different starting countries}
    \label{fig:CorplotHoldOut2}
\end{figure}

\end{appendix}

\end{document}